\documentclass[11pt]{amsart}
\allowdisplaybreaks[4]
\usepackage{tabularx}
\usepackage{amssymb}
\usepackage{amsfonts}
\usepackage{graphicx}
\usepackage{epstopdf}
\usepackage{dcolumn}
\usepackage{amsmath}
\usepackage{enumerate}
\usepackage{latexsym,bm}
\usepackage{slashed}
\usepackage{float}
\usepackage{cite}
\usepackage{CJK}
\usepackage[all,cmtip]{xy}
\usepackage{appendix}
\usepackage{ulem}
\usepackage[colorlinks,
linkcolor=black,
citecolor=red
]{hyperref}

\DeclareMathOperator{\Span}{Span}

\DeclareMathOperator{\Supp}{Supp}

\DeclareMathOperator{\CW}{CW}
\DeclareMathOperator{\adm}{adm}
\DeclareMathOperator{\cell}{cell}
\DeclareMathOperator{\NE}{NE}
\DeclareMathOperator{\NF}{NF}

\DeclareMathOperator{\dist}{dist}

\DeclareMathOperator{\id}{Id}

\DeclareMathOperator{\im}{im}

\DeclareMathOperator{\GLMY}{GLMY}

\newcommand{\Z}{\mathbb Z}

\newcommand{\R}{\mathbb R}

\theoremstyle{plain} %\newtheorem{Cor}{Corollary}[section]
\newtheorem{theorem}{Theorem}[section]
\newtheorem{lemma}[theorem]{Lemma}
\newtheorem{lem/def}[theorem]{Lemma/Definition}
\newtheorem{cor/def}[theorem]{Corollary/Definition}
\newtheorem{proposition}[theorem]{Proposition}
\newtheorem{corollary}[theorem]{Corollary}
\newtheorem{claim}[theorem]{Claim}
\newtheorem{conjecture}[theorem]{Conjecture}

\newtheorem{question}[theorem]{Question}
\theoremstyle{definition}
\newtheorem{definition}[theorem]{Definition}
\newtheorem{example}[theorem]{Example}

\newtheorem{remark}[theorem]{Remark}
\numberwithin{equation}{section}

\newenvironment{acknowledgement}{\smallskip{\sc Acknowledgments.}\rm}{\smallskip}

\allowdisplaybreaks

\begin{document}

\title{The Cellular Homology of Digraphs}
\author{Xinxing Tang and Shing-Tung Yau}

  \address{
X. Tang: Beijing Institute of Mathematical Sciences and Applications, Beijing, China;
}
\email{tangxinxing@bimsa.cn}
  \address{
S.-T. Yau: Yau Mathematical Sciences Center, Tsinghua University, and Beijing Institute of Mathematical Sciences and Applications, Beijing, China;\\
Beijing Institute of Mathematical Sciences and Applications, Beijing, China;
}
\email{styau@tsinghua.edu.cn}

\maketitle

\begin{abstract} In \cite{TY}, we investigate the pair $(P, \Supp(P))$ of minimal path $P$ and its supporting sub-digraph $\Supp(P)$ in the path complex of a digraph $G$ under the strongly regular condition. In this paper, first, we consider the special minimal path $P$ specified by the admissible condition (Definition \ref{admpair}), which means that $(P,\Supp(P))$ admits a singular cubical realization. Based on such a subset, we systematically introduce the definitions of cellular chain complex associated to $G$ and prove the well-definedness. Then we study several properties of such cellular homologies. Finally, we present several intriguing examples as well as some important observations.
\end{abstract}

\tableofcontents

\section{Introduction}

Since 2012, Grigoryan, Lin, Muranov and Yau develop a homotopy-homology approach to digraphs and initiate several studies of the relationship between discrete topology (where the discrete topology is played by objects of various categories of digraphs or graphs) and continuous topology, which is now known as GLMY-theory. Compared with other previously studied notions of graph homologies (e.g. via cliques \cite{CYY}, or via Hochschild (co)homology \cite{Hochschild}), it boasts many commendable properties (\cite{DIMZ,GJM,GJMY,GLMY3,GLMY,GLMY2}):
\begin{enumerate}
  \item It ``generalises" the homology of simplicial complex.
  \item It satisfies the properties that are analogous to Eilenberg-Steenrod axioms.
  \item It is well linked to graph-theoretical operations:
  \begin{itemize}
    \item it respects the digraph suspension;
    \item the K\"{u}nneth formula holds for both Cartesian product of two digraphs and join of two digraphs.
  \end{itemize}
  \item There is an appropriate notion of fundamental group and covering digraphs.
  \item There is a dual cohomology theory of digraphs which coincides with the one developed by Dimakis and M\"{u}ller-Hoissen in \cite{DM1,DM} as a reduction of the ``universal differential algebra".
\end{enumerate}

In \cite{TY}, we study the structure of (a sub-complex of) path complex of a digraph $G$. More explicitly, we focus on the $\Z$-coefficient generators of the path complex \textit{under the strongly regular condition}\footnote{The ``strongly regular condition" and ``minimal path" are originally proposed in the Huang-Yau paper\cite{HY}. They also give the definition of supporting digraph, but it is a little different from our definition \ref{supp}.} (that is, we stipulate that each elementary path must comprise distinct vertices, see Definition \ref{strongregular}), study the so-called minimal path $P$ (see, Definition \ref{minimaldef}) and its supporting digraph $\Supp(P)$ (see, Definition \ref{supp}).

\begin{question} Why do we consider the strongly regular condition?
\end{question}

On the one hand, given a finite simple digraph $G$, the GLMY's regular condition allows the path of the form
$$e_{i_0i_1\cdots i_n},\quad i_j\neq i_{j+1},~j=0,1,\ldots,n-1.$$
Under such a condition, the path complex as well as the path homology may be of infinite degree. For example,
\begin{enumerate}
  \item For the circulant digraph $C_{5}^{1,2}$ in Subsection \ref{circ}:
  $$\dim\Omega_n^{\GLMY}(C_5^{1,2})>0,\quad n\in\mathbb{N}.$$
  See Remark \ref{infinitepc} for an explanation;
  \item For the digraph $G$ in Appendix A of \cite{TY}:
  $$\dim H_{3n+1}^{\GLMY}(G)>0,\quad n\in\mathbb{N}.$$
\end{enumerate}

Thus, by definition, the strongly regular requirement gives rise to a finite sub-complex of GLMY's path complex. This sub-complex, by virtue of its finiteness, becomes more amenable to computational analysis and theoretical exploration. But it also presents several limitations and challenges. For example, in general,
$$\partial f_*\neq f_*\partial,$$
where $f$ is the digraph map and the push-forward $f_*$ is defined in \ref{f*}. Of course, when we consider the category of acyclic digraphs, the distinction between the strongly regular condition and the GLMY's regular condition disappears, then the functoriality holds in such a narrow case.

On the other hand, it is based on the topological/geometric consideration. We hope the minimal $n$-path $P$ and its supporting digraph $\Supp(P)$ could play the role of the regular $n$-cell in the topological space. To propel this idea, we study many examples and dissect the structure of $P$ and $\Supp(P)$ in \cite{TY}. The key results in \cite{TY} are the structure theorem (see, Theorem \ref{structurethm}) and acyclic result for $\widetilde{H}_*(\Supp(P))$. Immediately, we obtain several applications based on the standard acyclic model theorem. However, it is indeed meritorious to highlight that, in the course of our investigations, we have also made the following noteworthy observations:
\begin{itemize}
  \item For any integer $n\geq3$, there are infinitely many kinds of minimal $n$-paths.
  \item The digraph $\Supp(P)$ may not be contractible (Definition \ref{homotopyequiv}), such as the exotic cube\footnote{GJM \cite{GJM} show that the exotic cube possesses a nontrivial singular cubical cycle in degree 2, and they also prove the homotopy invariance of singular cubical homology. Thus, the exotic cube is not contractible.} (see, Example \ref{exotic}). More explicitly, for each $n\geq3$, only finite kinds of $\Supp(P)$ exhibit the property of contractibility. Specifically, as the number of elementary path components within the minimal path $P$ increases, the contractibility of $\Supp(P)$ will disappear.
\end{itemize}

Contractibility stands as a cornerstone property in the topological spaces, significantly influencing further constructions and analyses.
For example, the fiber bundle is trivial over the contractible space.

Due to the aforementioned reasons, this paper endeavors to establish a robust and nuanced connection between digraphs and topological spaces by studying special minimal paths under an additional condition and by examining the $\CW$ complex and cellular homology of the digraph $G$.

\subsection{Main results}

The key results in this paper are as follows.

\vskip 0.2cm
(1) By introducing the admissible conditions, we formulate the definition of the cellular chain complex of digraphs and prove the well-definedness.

\begin{definition}[Definition \ref{admpair}, \ref{admrel}] Let $G$ be a digraph and $P$ be a minimal $n$-path. The pair $(P,\Supp(P))$ is called \textbf{admissible}, if there exists a digraph map
$$f:I^{\boxdot n}\longrightarrow G$$
such that
\begin{itemize}
  \item $f_*(\omega_n)=cP$, for some $c\in\mathbb{Z}\setminus\{0\};$%$f_*(\omega_n)=P$ or $f_*(\omega_n)=-P$;
  \item the image digraph of $f$ is $\Supp(P)$.
\end{itemize}
%The set of such admissible pairs are denoted by $\mathcal{P}_{\adm,n}(G)$.
The linear relations among the admissible paths are called the \textbf{admissible relations}.
\end{definition}
\vskip 0.2cm

\begin{definition}[Definition \ref{cellular}] Let $\mathcal{P}_{\adm,n}(G)$ be the set of length $n$ admissible pairs (up to a sign) in $G$. We define the $n$-th cellular chain of $G$ to be $$C_n(G;\R)=\left\{\sum_{P\in\mathcal{P}_{\adm,n}}c_PP~\bigg|~c_P\in \R\right\}\bigg/\left\{\text{admissible relations}\right\}.$$
Moreover, we define a linear operator $\partial^{\cell}:C_n(G;\R)\longrightarrow C_{n-1}(G;\R)$. For any $P\in\mathcal{P}_{\adm,n}(G)$,
$$\partial^{\cell}P=\sum_{Q\in\mathcal{P}_{\adm,n-1}(G)}[Q:P]Q,$$
where the coefficient $[Q:P]$ is defined as follows:
\begin{equation*}
[Q:P]=\left\{
        \begin{array}{ll}
          1, & \hbox{if $Q<\partial P$;} \\
          -1, & \hbox{if $-Q<\partial P$;} \\
          0, & \hbox{otherwise.}
        \end{array}
      \right.
\end{equation*}
\end{definition}

Through the meticulous and intricate dissection of the relation between $\partial^{\cell}$ in the cellular complex and $\partial$ in the path complex, we elucidate the well-defined nature of $\partial^{\cell}$. Moreover, we arrive at the definition of the cellular homology of $G$, which are denoted by $H^{\cell}_*(G;\R)$.

Our construction is related to the singular cubical homology \cite{GJM} (which is recalled in Subsection \ref{singularcube}). In particular, the cellular chain complex could be understood as the image complex from singular cubical complex to path complex under strongly regular condition. In particular, after discussing several relations between our cellular homologies and singular cubical homologies of digraphs, we propose the following conjecture:

\begin{conjecture}[Conjecture \ref{conj}] If the digraph $G$ is acyclic, then
$$H_*^{\cell}(G;\R)\cong H_*^c(G;\R).$$
In particular, the admissible relations among the admissible pairs are the boundary conditions in the singular cubical homology.
\end{conjecture}

\vskip 0.2cm

(2) We prove several properties of the cellular homologies associated to digraphs.
\begin{itemize}
  \item Functoriality and homotopy invariance in the category of acyclic digraphs.
  \item Acyclic result for $\widetilde{H}^{\cell}_*(\Supp(P);\R)$ for $P\in\mathcal{P}_{\adm}(G)$.
  \item K\"{u}nneth formula for Cartesian product of two digraphs.
\end{itemize}

In particular, after checking the contractibility of all the admissible minimal $3$-paths (up to global orientation of edges) in Appendix \ref{AppendixA}, we also propose the following conjecture.
\begin{conjecture}[Conjecture \ref{contractible}] Let $P\in\mathcal{P}_{\adm}(G)$, then $\Supp(P)$ is contractible.
\end{conjecture}

\subsection{Organization}

In Section \ref{2}, we recall the necessary definitions and results in \cite{TY}. In Section \ref{3}, we initiate our construction process of $\CW$ complex and cellular chain group of a digraph $G$. In Subsection \ref{singularcube}, we recall the singular cubical complex of a digraph in \cite{GJM} and their relation with path complex. In Subsection \ref{3.2}, we introduce the definitions of admissible paths and admissible relations. The various examples are meticulously chosen to highlight the subtle yet significant distinctions between the traditional definitions and our novel approach. In Subsection \ref{3.3}, we define the cellular homology of a digraph and conduct a series of fundamental calculations. In subsection \ref{properties}, we analyze several properties of such cellular homology. In Section \ref{4}, we compute the cellular homologies of several interesting digraphs, such as special circulant digraphs and Johnson digraphs. In Section \ref{5}, we proceed with discussion, garnering a series of pivotal insights. In Appendix \ref{AppendixA}, we enumerate all the admissible minimal $3$-paths up to global orientation of edges and construct the corresponding deformation retractions.

\begin{acknowledgement} The authors thank Xin Fu, Sergei Ivanov, Yuri Muranov and Fedor Pavutnitskiy for useful discussions. The authors extend their gratitude to the anonymous reviewer for the meticulous perusal and the insightful suggestions. The work of X.T. is partially supported by research funding at BIMSA.
\end{acknowledgement}

\bigskip

\section{Path Complex and Minimal Path}\label{2}

In \cite{TY}, we analyze the structure of $\Z$-coefficient path complex $\Omega_*(G;\mathbb{Z})$ of a digraph $G$ via the minimal path $P$ and the corresponding supporting digraph $\Supp(P)$. In this section, we will recall necessary results on the path complex and the corresponding $\mathbb{Z}$-structure.

\subsection{Path complex}

Let $V$ be a finite set. For any $n\in\mathbb{N}$, an elementary $n$-path is any ordered sequence $i_0,\ldots, i_n$ of $n+1$ vertices of $V$, write it as $e_{i_0\ldots i_n}$. Let $\Lambda_n(V;\mathbb{Z})$ be the $\Z$-module generated by all elementary $n$-paths. We also call $n$ the length of the path in $\Lambda_n(V;\mathbb{Z})$.

The $\Z$-homomorphism $\partial:\Lambda_n(V;\mathbb{Z})\rightarrow \Lambda_{n-1}(V;\mathbb{Z})$ is defined via the generators:
$$\partial e_{i_0\ldots i_n}=\sum_{j=0}^n(-1)^je_{i_0\ldots \widehat{i_j}\ldots i_n}.$$
Clearly, $\partial\left(\Lambda_n(V;\mathbb{Z})\right)\subset\Lambda_{n-1}(V;\mathbb{Z})$, and $\partial^2=0$.

Now we will consider a sub-complex of $(\Lambda_*(V;\mathbb{Z}),\partial)$ which is constrained by a digraph $G$. More explicitly, let $G=(V,E)$ be a finite simple digraph\footnote{Here simple condition means that there are no multiple edges and self loops. In this paper, a digraph means a finite simple digraph.}, where $V$ is the set of vertices and $E$ is the set of directed edges. First, we have the following definition.
\begin{definition}[\cite{HY}]\label{strongregular}
Let $G=(V,E)$ be a finite simple digraph.
\begin{enumerate}
  \item An elementary path $e_{i_0i_1\ldots i_n}$ is called \textbf{allowed}, if each directed edge $i_{k-1}\rightarrow i_k$ belongs to $E$, $k=1,2\ldots, n$.
  \item An elementary path $e_{i_0i_1\ldots i_n}$ is called \textbf{strongly regular} (abbreviated as \textbf{s-regular}), if $i_0, i_1,\ldots,i_n$ are all distinct.
\end{enumerate}
Furthermore, we denote the set of s-regular allowed elementary $n$-paths in $G$ by $E_n(G)$.
\end{definition}

Let $\mathcal{A}_n(G;\Z)$ be the free $\Z$-module generated by the set $E_n(G)$. Note that, $\mathcal{A}_n(G;\Z)=0$, if $n\geq\big|V(G)\big|$.

In general, the graded module $\mathcal{A}_*(G;\Z)$ is not preserved under $\partial$. Thus we consider the submodule $\Omega_n(G;\Z)$ of $\partial$-invariant s-regular allowed $n$-paths. That is
$$\Omega_n(G;\Z)=\{p\in \mathcal{A}_{n}(G;\Z)| \partial p\in \mathcal{A}_{n-1}(G;\Z)\}.$$
%We call $k$ the length of the path.

\begin{definition}\label{pathcomplex} The complex $(\Omega_*(G;\Z), \partial)$ is called the path complex of a digraph $G$ with $\Z$-coefficients. For each $n\geq0$,
$$H_n(G,\Z):=\frac{\ker(\partial:\Omega_n(G;\Z)\rightarrow \Omega_{n-1}(G;\Z))}{\im(\partial:\Omega_{n+1}(G;\Z)\rightarrow \Omega_{n}(G;\Z))},$$
it is called the $n$-th path homology of $G$.
\end{definition}

\begin{remark}[GLMY's regular condition and GLMY's path complex]\label{regular}
The s-regular condition in Definition \ref{strongregular} is indeed stronger than GLMY's regular condition in \cite{GLMY3}. In \cite{GLMY3}, the set of
regular $n$-paths is defined as a quotient space $\mathcal{R}_n(V)$:
$$\mathcal{R}_n(V)=\Lambda_n(V)/I_n(V),$$
where $I_n(V)$ is the submodule generated by the irregular path $e_{i_0i_1\ldots i_n}$, that is, $i_{k-1}=i_k$ for some $k=1,\ldots,n$. One can easily check that $(\mathcal{R}_*(V),\partial)$ forms a quotient chain complex with the induced boundary operator, we still denote by $\partial$.

GLMY \cite{GLMY3} consider the submodule $\mathcal{A}_n^{\GLMY}(G)\subset\mathcal{R}_n(G)$ and $\partial$-invariant submodule $\Omega_n^{\GLMY}(G)$, then give the definition of path complex of a digraph $G$. We have,
\begin{itemize}
  \item The path complex in Definition \ref{pathcomplex} is a sub-complex of the GLMY's path complex. There are many examples that both path homologies are isomorphic. In particular, if the digraph $G$ is acyclic, that is, G does not contain the monotone cycle of the form
  $$v_0\rightarrow v_1\rightarrow v_2\rightarrow\cdots\rightarrow v_n\rightarrow v_0,$$
then the s-regular condition is the same as the regular condition. It follows that the corresponding homologies are the same.
  \item However, there also exist counterexamples that the two kinds of path homologies are not the same. Such an example could be found in Appendix A of \cite{TY}.
\end{itemize}
\end{remark}

In what follows, for simplicity, we will work with the s-regular condition.

\subsection{Minimal path}

For any $P=\sum_{p=1}^mc_pe_p\in \Omega_n(G;\Z)$ with $e_p\in E_n(G)$, we define
$$w(P)=\sum_{p=1}^m |c_p|,$$
and call it the width of the path $P$.

\begin{definition}[\cite{HY}]\label{minimaldef} Let $P=\sum_{p=1}^mc_pe_p\in\Omega_n(G;\Z)$, if there exists
$$P'=\sum_{p=1}^md_pe_p\in\Omega_n(G;\Z)\setminus\{0\},$$
satisfying
\begin{enumerate}
  \item $|c_p-d_p|\leq |c_p|$ and $|d_p|\leq |c_p|$ for each $p=1,\ldots,m$;
  \item $w(P')<w(P)$;
\end{enumerate}
then we call $P'$ is smaller than $P$, and denote by $P'<P$. If there does not exist such $P'$, we call $P$ a minimal $n$-path.
\end{definition}

We have the following immediate observations on the minimal path:
\begin{itemize}
  \item If $P$ is minimal, $-P$ is also minimal.
  \item If $P'<P$, then $P-P'\in\Omega_n(G;\Z)$ is also smaller than $P$.
  \item Minimal $0$-path is just represented by a single point (up to a sign); minimal $1$-path is just represented by $e_{ij}$ (up to a sign), $(i\rightarrow j)\in E$.
  \item Any element in $\Omega_n(G;\Z)$ is a $\Z$-linear combination of minimal paths.
  \item Any minimal path is a $\Z$-linear combination of s-regular allowed elementary paths with the same starting vertex (denote by $S$) and ending vertex (denote by $E$).
\end{itemize}

The path complex, especially the minimal path, depends on the digraph. To understand its structure, we introduce the following definitions.

\begin{definition}[\cite{TY}]\label{supp} For each minimal path $P$ in the digraph $G$, we define $\Supp(P)$ to be the minimal subgraph of $G$ such that $P\in\Omega_*(\Supp(P))$.
\end{definition}

One can obtain the sub-digraph $\Supp(P)$ as follows:
\begin{enumerate}
  \item First, we express $P$ uniquely in terms of the elementary paths:
        $$P=\sum_{p=1}^lc_pe_p,\quad e_p\in E_n(G;\Z),~ c_p\in\Z\setminus\{0\}.$$
        Then we get a digraph $\Supp(P)^{\mbox{pre}}$ given by $\{e_p\}_{p=1}^l$.
  \item Second, we express $\partial P$ uniquely in terms of the elementary paths:
        $$\partial P=\sum_{q=1}^md_qe_q,\quad e_q\in E_{n-1}(G;\Z),~d_q\in\Z\setminus\{0\}.$$
        Then we add necessary edges from $\{e_q\}_{q=1}^m$ to $\Supp(P)^{\mbox{pre}}$ and obtain $\Supp(P)$.
\end{enumerate}

For convenience, we also introduce two ``distance functions" on $V(P)$. First, for each s-regular allowed elementary path $e_I=e_{i_0i_1\ldots i_n}$, we have the position function
$$f_I:V(e_I)~\longrightarrow~\mathbb{N},\quad f_I(i_a)=a.$$
Then, for a minimal $n$-path $P$ with starting vertex $S$ and ending vertex $E$, we define
\begin{align*}
&d_S:V(P)\rightarrow \mathbb{N},\quad d_S(v)=\min\{f_I(v)|e_I\text{ is a component of }P\};\\
&d_E:V(P)\rightarrow \mathbb{N},\quad d_E(v)=\min\{n-f_I(v)|e_I\text{ is a component of }P\}.
\end{align*}
One can understand $d_S$ and $d_E$ as the distance from $v$ to $S$ and $E$ in $P$ (not in $\Supp(P)$) respectively.

Now let us give several examples on the minimal paths and the corresponding supporting digraphs as follows.

\begin{example}\label{exsuppdef} In the following 3-simplex digraph $G$, $P=e_{0123}$ is minimal, and its supporting subdigraph is given by
\begin{figure}[H]
	\centering
	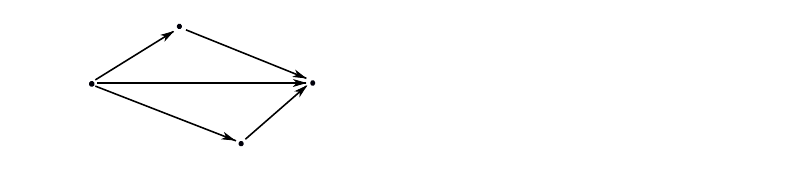
	%\caption[]{$\Supp(e_{0123})$}
	%\label{Fig:Q1}
\end{figure}
By definition, we have
\begin{align*}
&d_S(0)=0, \quad d_S(1)=1,\quad d_S(2)=2,\quad d_S(3)=3;\\
&d_E(0)=3, \quad d_E(1)=2,\quad d_E(2)=1,\quad d_E(3)=0.
\end{align*}
\end{example}

\begin{example}\label{EXfSfE1} The minimal digraph which supports
$$P=e_{0136}-e_{0156}+e_{0456}+e_{0246}-e_{0236}$$
is given by
\begin{figure}[H]
	\centering
	%% Creator: Inkscape 1.0.1 (3bc2e813f5, 2020-09-07), www.inkscape.org
%% Accompanies image file 'EXfSfE1.pdf' (pdf, eps, ps)
\begingroup%
  \makeatletter%
  \providecommand\color[2][]{%
    \errmessage{(Inkscape) Color is used for the text in Inkscape, but the package 'color.sty' is not loaded}%
    \renewcommand\color[2][]{}%
  }%
  \providecommand\transparent[1]{%
    \errmessage{(Inkscape) Transparency is used (non-zero) for the text in Inkscape, but the package 'transparent.sty' is not loaded}%
    \renewcommand\transparent[1]{}%
  }%
  \providecommand\rotatebox[2]{#2}%
  \newcommand*\fsize{\dimexpr\f@size pt\relax}%
  \newcommand*\lineheight[1]{\fontsize{\fsize}{#1\fsize}\selectfont}%
  \ifx\svgwidth\undefined%
    \setlength{\unitlength}{200.57693842bp}%
    \ifx\svgscale\undefined%
      \relax%
    \else%
      \setlength{\unitlength}{\unitlength * \real{\svgscale}}%
    \fi%
  \else%
    \setlength{\unitlength}{\svgwidth}%
  \fi%
  \global\let\svgwidth\undefined%
  \global\let\svgscale\undefined%
  \makeatother%
  \begin{picture}(1,0.4268118)%
    \lineheight{1}%
    \setlength\tabcolsep{0pt}%
    \put(0,0){\includegraphics[width=\unitlength,page=1]{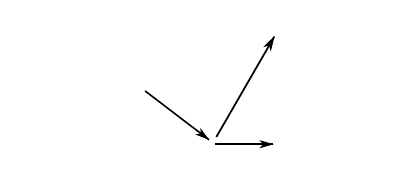}}%
    \put(0.28485588,0.19341972){\makebox(0,0)[lt]{\lineheight{1.25}\smash{\begin{tabular}[t]{l}$0$\end{tabular}}}}%
    \put(0.48873363,0.38519485){\makebox(0,0)[lt]{\lineheight{1.25}\smash{\begin{tabular}[t]{l}$1$\end{tabular}}}}%
    \put(0.48461972,0.00693325){\makebox(0,0)[lt]{\lineheight{1.25}\smash{\begin{tabular}[t]{l}$2$\end{tabular}}}}%
    \put(0.64651195,0.38756623){\makebox(0,0)[lt]{\lineheight{1.25}\smash{\begin{tabular}[t]{l}$3$\end{tabular}}}}%
    \put(0.65266385,0.00538687){\makebox(0,0)[lt]{\lineheight{1.25}\smash{\begin{tabular}[t]{l}$4$\end{tabular}}}}%
    \put(0.76289923,0.14730498){\makebox(0,0)[lt]{\lineheight{1.25}\smash{\begin{tabular}[t]{l}$5$\end{tabular}}}}%
    \put(0,0){\includegraphics[width=\unitlength,page=2]{EXfSfE1.pdf}}%
    \put(-0.00094961,0.19680724){\makebox(0,0)[lt]{\lineheight{1.25}\smash{\begin{tabular}[t]{l}$\Supp(P)=$\end{tabular}}}}%
    \put(0,0){\includegraphics[width=\unitlength,page=3]{EXfSfE1.pdf}}%
    \put(0.94546104,0.18469101){\makebox(0,0)[lt]{\lineheight{1.25}\smash{\begin{tabular}[t]{l}$6$\end{tabular}}}}%
  \end{picture}%
\endgroup%

	%\caption[]{$\Supp(e_{0123})$}
	%\label{Fig:Q1}
\end{figure}
By definition, we have
\begin{align*}
&d_S(0)=0,~~ d_S(1)=1,~~ d_S(2)=1,~~ d_S(3)=2,~~ d_S(4)=1,~~ d_S(5)=2,~~d_S(6)=3;\\
&d_E(0)=3,~~ d_E(1)=2,~~ d_E(2)=2,~~ d_E(3)=1,~~ d_E(4)=1,~~ d_E(5)=1,~~d_E(6)=0.
\end{align*}
\end{example}

\begin{example}\label{more} The minimal digraph which supports
\begin{align*}
P=~&e_{S06E}-e_{S16E}+e_{S17E}-e_{S27E}+e_{S28E}-e_{S38E}\\
   &+e_{S39E}-e_{S49E}+e_{S4(10)E}-e_{S5(10)E}.
\end{align*}
is given by
\begin{figure}[H]
	\centering
	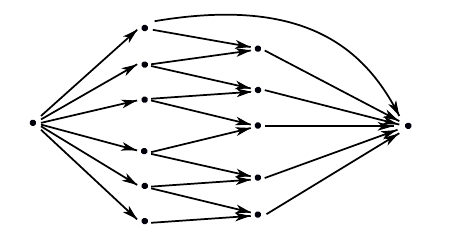
	%\caption[]{$\Supp(e_{0123})$}
	%\label{Fig:Q1}
\end{figure}
\end{example}

\begin{example}\label{length4}
Let us consider the following path
\begin{align*}
P=~&e_{S159E}-e_{S169E}+e_{S269E}\\
   &+e_{S16(10)E}-e_{S26(10)E}+e_{S27(10)E}-e_{S37(10)E}\\
   &-e_{S27(11)E}+e_{S37(11)E}-e_{S38(11)E}+e_{S48(11)E}.
\end{align*}
It could be a minimal $4$-path with the following supporting digraph
\begin{figure}[H]
	\centering
	\input{length4.pdf_tex}
	%\caption[]{$\Supp(e_{0123})$}
	%\label{Fig:Q1}
\end{figure}
\end{example}

We have the following structure theorem for $P$ and $\Supp(P)$. The readers can understand it in terms of the above examples.

\begin{theorem}[\cite{TY} Structure Theorem]\label{structurethm} Let $P\in\Omega_n(G;\Z)$ be a minimal path with the starting vertex $S$ and ending vertex $E$, $\Supp(P)$ be its supporting digraph and $d_S$, $d_E$ be the functions defined above.

(1) Let $S_1=d_S^{-1}(1)$ and $E_1=d_E^{-1}(1)$. Then
$P$ is a linear combination of s-regular allowed elementary paths from $S$ to $E$, with coefficients being either 1 or -1. And in $\Supp(P)$,
\begin{equation*}
\begin{aligned}
\partial P=\sum_{\alpha\in E_1}P_{S,n-1,\alpha}+\sum_{\beta\in S_1}P_{\beta,n-1,E}+\sum_{k\in I_P} P_{S,n-1,E}^k,
\end{aligned}
\end{equation*}
where
\begin{itemize}
  \item  $P_{S,n-1,\alpha},P_{\beta,n-1,E},P_{S,n-1,E}^k\in\Omega_{n-1}(\Supp(P);\Z)$, and moreover
  \begin{itemize}
          \item $P_{S,n-1,\alpha}$ is the minimal $(n-1)$-path starting from $S$ and ending with $\alpha$;
          \item $P_{\beta,n-1,E}$ is the minimal $(n-1)$-path starting from $\beta$ and ending with $E$;
          \item $P_{S,n-1,E}^k$ is the minimal $(n-1)$-path starting from $S$ and ending with $E$.
        \end{itemize}
  \item Such $P_{S,n-1,\alpha}$, $P_{\beta,n-1,E}$ are unique (up to a sign) in $\Supp(P)$ for each $\alpha\in E_1$, $\beta\in S_1$.
  \item The set $I_P$ in the last summand depends on $P$, and $|I_P|\leq 1$. That is, there exists at most one (up to a sign) minimal $(n-1)$-path in $\Supp(P)$ starting from $S$ and ending with $E$.
  %For each $k$ in the last summand, $P_{S,n-1,E}^k$ in $\partial P$ is unique (up to a sign) in $\Supp(P)$.
\end{itemize}

(2) For any $v\in d_E^{-1}(k)\cap\Supp(P)$, in $\Supp(P)$,
\begin{itemize}
  \item There is a unique minimal $(n-k)$-path (up to a sign) starting from $S$ and ending with $v$, denoted by $P_{S,n-k,v}$, as well as only one minimal $k$-path (up to a sign) starting from $v$ and ending with $E$, denoted by $P_{v,k,E}$.
  \item There is at most one minimal $(n-k-1)$-path (up to a sign) starting from $S$ and ending with $v$, denoted by $P_{S,n-k-1,v}$, as well as at most one minimal $(k-1)$-path starting from $v$ and ending with $E$, denoted by $P_{v,k-1,E}$.
  %\item There is no minimal path of length $<k-1$ starting from $v$ and ending with $E$.
\end{itemize}

(3) Each minimal $2$-path with fixed starting and ending vertices in $\Supp(P)$ ($n\geq2$) is unique.
\end{theorem}

In \cite{TY}, we also prove the following acyclic result.
\begin{theorem}[\cite{TY}]\label{acyclic} Let $P$ be a minimal path of $G$. Then
$$\widetilde{H}_*(\Supp(P))=0.$$
\end{theorem}
We switch to the complex with field (characteristic 0) coefficient, for example, $\R$, then the basis of $\Omega_*(G;\R)$ consisting of minimal paths gives the acyclic model of the path complex. Finally, we discuss several canonical applications of such an acyclic model in the path complex.

Meanwhile, we have the following two observations:
\begin{itemize}
  \item For any integer $n\geq3$, there are infinitely many minimal $n$-paths.
  \item Even if $\widetilde{H}_*(\Supp(P))=0$, the digraph $\Supp(P)$ may not be contractible. More explicitly, for each $n\geq3$, only finite kinds of $\Supp(P)$ have trivial cellular homologies defined in Section \ref{3}, then they are possibly contractible by the homotopy invariance \ref{homotopyinv}.
\end{itemize}

In the next section, we will consider a sub-complex of path complex, where certain minimal paths are specified by one more condition.
\bigskip

\section{The Cellular Homology of Digraphs}\label{3}

In this section, we will consider the $\CW$ complex associated to a digraph by specifying certain pairs (which we call admissible pairs) of minimal path and its supporting digraph. %It is isomorphic to some sub-complex of the path complex as a vector space. Roughly speaking, such a vector space is some kind of intersection of the path complex and singular cubical complex of a digraph $G$.
The admissible paths may be linearly dependent as paths. Our definition of $\CW$ complex depends on the choice of maximal linearly independent subset of admissible paths.
To circumvent the intricacies involved in the choice of a basis, we define the cellular chain group of the digraph as the quotient of the free group generated by so-called admissible pairs, modulo the linear relations imposed thereupon.

\subsection{The singular cubical complex of digraphs}\label{singularcube}

Inspired by the GJM's singular cubical homology theory, we introduce the admissible condition and further construction of celluar complex of digraphs. Now let us recall the corresponding necessary definitions and results for singular cubical complex of digraphs. For more details, please refer to \cite{GJM}.

\begin{definition}[Cartesian product of digraphs] For two digraphs $G=(V(G), E(G))$ and $H=(V(H), E(H))$, we define the Cartesian product $G\boxdot H$ as a digraph with the set of vertices $V(G)\times V(H)$ and with the set of edges as follows: for $v,v'\in V(G)$ and $w,w'\in V(H)$, we have
$(v,w)\rightarrow (v',w')$ in $G\boxdot H$ if and only if
$$\text{either $\{v=v',~w\rightarrow w'\}$, or $\{v\rightarrow v',~w=w'\}$}.$$
\end{definition}

Let $I$ be the digraph $I=0~\bullet~\longrightarrow~1~\bullet$. In the following discussion, for $n\in\mathbb{N}$, we will focus on the $n$-cube
$$I^{\boxdot n}=\underbrace{I\boxdot I\boxdot\cdots\boxdot I}_{n},~n\geq1;\quad I^{\boxdot 0}=\{0\}.$$
Where the vertices are represented by the binary sequence. Any of its $(n-1)$-face is of the form
$$I^{\boxdot k}\boxdot\{\epsilon\}\boxdot I^{\boxdot(n-1-k)},\quad 0\leq k\leq n-1,~\epsilon=0,1.$$
There are also natural inclusions $F_{j\epsilon}:I^{\boxdot(n-1)}\rightarrow I^{\boxdot n}$, for $1\leq j\leq n$,
\begin{equation}\label{n-1_n_inclusion}
F_{j\epsilon}(c_1,\ldots,c_{n-1})=(c_1,\ldots,c_{j-1},\epsilon,c_j,\ldots,c_{n-1}),\quad n\geq2,
\end{equation}
and $F_{1\epsilon}(0)=\epsilon$ for $n=1$.

%Denote by $E_n$ the set of allowed elementary $n$-paths between $(0,\ldots,0)$ and $(1,\ldots,1)$ in $I^{\boxdot n}$. Moreover, define the following allowed $n$-path in $I^{\boxdot n}$:
%$$\omega_n=\sum_{e_{\alpha}\in E_n}(-1)^{\sigma_n(\alpha)}e_{\alpha},$$
%where
%\begin{itemize}
%  \item $\alpha=\alpha_0~\rightarrow~\alpha_1~\rightarrow\cdots\rightarrow~\alpha_n$ with the vertex $\alpha_k=(\alpha_{k,i})_{i=1}^n$;
%  \item $\sigma(\alpha)$ is the number of inversions in the sequence $\{\sum_{i=1}^n2^i\alpha_{k,i}\}$.
%\end{itemize}

%$\alpha=\alpha_0~\rightarrow~\alpha_1~\rightarrow\cdots\rightarrow~\alpha_n$ with the vertex $\alpha_k=(\alpha_{k,i})_{i=1}^n$ and $\sigma(\alpha)$ is the number of inversions in the sequence $\{\sum_{i=1}^n2^i\alpha_{k,i}\}$.

\begin{definition}[Digraph map] A morphism from a digraph $G=(V(G),E(G))$ to a digraph $H=(V(H),E(H))$ is a map $f:V(G)\rightarrow V(H)$ such that for any edge $v\rightarrow w$ in $G$, we have
$$f(v)\rightarrow f(w)\quad\text{or}\quad f(v)=f(w) \text{ in }H.$$
The requirement is also denoted by $f(v)\overset{\rightarrow}{=}f(w)$. We also call such a morphism to be a digraph map, and denote by $f:G\rightarrow H$.
\end{definition}

A singular $n$-cube in a digraph $G$ is a digraph map $\phi:I^{\boxdot n}\rightarrow G$. Of course, the singular $n$-cubes could be degenerate enough, since $I^{\boxdot n}$ maps to any single point, leading to a singular $n$-cube. Such a map captures less information of a digraph $G$. To obtain the more explicit information of a digraph via the singular cubes, we want to module the degenerate singular cubes, which are defined as follows.

For $n\geq1$ and $1\leq j\leq n$, consider the natural projection $T^j:I^{\boxdot n}\rightarrow I^{\boxdot(n-1)}$: For $n=1$,
$$T^1(0)=0,\quad T^1(1)=0.$$
For $n\geq2$,
$$T^j(i_1,\ldots,i_n)=(i_1,\ldots,i_{j-1},i_{j+1},\ldots,i_n).$$
%Such a map $T^j$ is called the projection on the $j$-face.
\begin{definition}[\cite{GJM}] We call the singular $n$-cube $\phi:I^{\boxdot n}\rightarrow G$ degenerate if there exist $j\in\{1,2,\ldots,n\}$ and a singular $(n-1)$-cube $\psi:I^{\boxdot(n-1)}\rightarrow G$, such that
$$\phi=\psi\circ T^j.$$
\end{definition}

For $n\geq0$, we denote by $Q_n(G)$ the free module generated by the singular $n$-cubes. And for $n\geq1$, we denote by $B_n(G)$ the free submodule generated by all degenerate $n$-cubes, $B_0(G)=0$.

For $n\geq1$, one can define a homomorphism $\partial^c:Q_n(G)\rightarrow Q_{n-1}(G)$ via the inclusions $\eqref{n-1_n_inclusion}$ by
$$\partial^c{\phi}=\sum_{j=1}^n(-1)^j\left(\phi\circ F_{j0}-\phi\circ F_{j1}\right)$$
and $\partial^c=0$ for $n=0$.

\begin{proposition}[\cite{GJM}, \cite{HW}] The homomorphism $\partial^c$ defines a boundary operator, that is $\left(\partial^c\right)^2=0$. Meanwhile, $\partial^c\big(B_n(G)\big)\subset B_{n-1}(G)$.
\end{proposition}

Following from this proposition, $(B_*(G),\partial^c)\subset (Q_*(G),\partial^c)$ is a sub-complex. Hence, the quotient complex
$$\left(\Omega_*^c(G)=Q_*(G)/B_*(G),\partial^c\right)$$
is well-defined, and it is called the (normalized) singular cubical chain complex of the digraph $G$. The associated homology is called the singular cubical homology of $G$, and denoted by $H_*^c(G)$.\\

There is a chain map from the singular cubical complex of $G$ to GLMY's path complex of $G$.

First, there is a standard generator $\omega_n$ of $\Omega_n(I^{\boxdot n};\Z)$. Denote by $E_n$ the set of allowed elementary $n$-paths between $(0,\ldots,0)$ and $(1,\ldots,1)$ in $I^{\boxdot n}$. Then $\omega_n\footnote{There is another way to define $\omega_n$ by using the cross product $\times$ (Definition \ref{cross}), $$\omega_n=\underbrace{e_{01}\times\cdots\times e_{01}}_n.$$}$ is defined to be
$$\omega_n=\sum_{e_{\alpha}\in E_n}(-1)^{\sigma_n(\alpha)}e_{\alpha},$$
where
\begin{itemize}
  \item $\alpha=\alpha_0~\rightarrow~\alpha_1~\rightarrow\cdots\rightarrow~\alpha_n$ with the vertex $\alpha_k=(\alpha_{k,i})_{i=1}^n$;
  \item $\sigma(\alpha)$ is the number of inversions in the sequence $\{\sum_{i=1}^n2^i\alpha_{k,i}\}$.
\end{itemize}

Second, one can define a homomorphism $\tau_n:\Omega_n^c(G)\rightarrow\Omega_n(G)$ on the singular $n$-cube $\phi$ by
$$\tau_n(\phi):=\phi_*(\omega_n),$$
where $\phi_*$ is a linear map defined on the elementary paths by
\begin{equation}\label{phi*}
\phi_*(e_{i_0i_1\ldots i_n})=\left\{
                      \begin{array}{ll}
                        e_{\phi(i_0)\phi(i_1)\ldots \phi(i_n)}, & \hbox{if all $\phi(i_j)\neq\phi(i_{j+1})$;} \\
                        0, & \hbox{otherswise.}
                      \end{array}
                    \right.
\end{equation}
Then one has
\begin{proposition}[\cite{GJM}]\label{singularcubpath} For $n\geq0$, and $\phi\in\Omega_n^c(G)$,
$$\partial\circ\tau_n(\phi)=\tau_{n-1}\circ\partial^c(\phi).$$
\end{proposition}

\subsection{Admissible pair and admissible relations}\label{3.2}

\subsubsection{Admissible pair}

Let $\widetilde{\mathcal{P}}_n(G)$ be the set of minimal $n$-paths in $G$. If $P\in\widetilde{\mathcal{P}}_n(G)$, so is $-P$. Furthermore, let $\mathcal{P}_n(G)$ be the subset of $\widetilde{\mathcal{P}}_n(G)$ containing exactly one of $P$ and $-P$. %More explicitly, for any $P_{n,\alpha}\in\widetilde{\mathcal{P}}_n(G)$, we have
%\begin{itemize}
%  \item $P_{n,\alpha}\in\mathcal{P}_n(G)$ and $-P\notin\mathcal{P}_n(G)$; or
%  \item $-P_{n,\alpha}\in\mathcal{P}_n(G)$ and $P\notin\mathcal{P}_n(G)$.
%\end{itemize}
In the following, we write
$$\mathcal{P}_n=\mathcal{P}_n(G)=\{P_{n,\alpha}\}_{\alpha},\quad \mathcal{P}=\bigcup_{n\geq0}\mathcal{P}_n.$$
And call them the set of minimal $n$-paths and the set of minimal paths (up to a sign), respectively.

In the following, we will investigate a subset of $\mathcal{P}$ which is related to the singular cube.

\vskip 0.2cm
First, we modify the definition $f_*$ for the digraph map $f:G\rightarrow H$ to work with our situation.

\begin{definition}\label{f*} Let $G$ and $H$ be two digraphs, and $f:G\rightarrow H$ be a digraph map. Then we can define a $\Z$-linear map
$$f_*:\mathcal{A}_*(G;\Z)\rightarrow\mathcal{A}_*(H;\Z)$$
by
\begin{equation*}
f_*(e_{i_0i_1\ldots i_k})=\left\{
                      \begin{array}{ll}
                        e_{f(i_0)f(i_1)\ldots f(i_k)}, & \hbox{if all $f(i_j)$'s are dinstict;} \\
                        0, & \hbox{otherswise.}
                      \end{array}
                    \right.
\end{equation*}
\end{definition}

\begin{remark} The original definition \eqref{phi*} of $f_*$ induces a chain map between GLMY's path complexes. For our restricted case, if $H$ is acyclic, then $f_*$ induces a chain map $f_*:\Omega_*(G)\rightarrow\Omega_*(H)$.
\end{remark}

With the help of the induced map $f_*$, we give our key definition.

\begin{definition}\label{admpair} Let $G$ be a digraph. For $P\in\mathcal{P}_n(G)$, we call the pair $(P,\Supp(P))$ \textbf{admissible}, if there exists a digraph map
$$f:I^{\boxdot n}\longrightarrow G$$
such that
\begin{enumerate}
  \item $f_*(\omega_n)=cP$ for some $c\in\Z\setminus\{0\}$;
  \item the image digraph of $f$ is $\Supp(P)$.
\end{enumerate}
We also say that $(P,\Supp(P))$ admits a singular cubical realization.
\end{definition}

Let $\mathcal{P}_{\adm}=\mathcal{P}_{\adm}(G)$ be the set of the admissible pairs in $G$, that is
$$\mathcal{P}_{\adm}=\left\{(P,\Supp(P))\text{ is admissible. }\big|~P\in\mathcal{P}\right\}.$$
Since $\Supp(P)$ is determined by $P$, we also regard $\mathcal{P}_{\adm}$ as a subset of $\mathcal{P}$, and write it as
$$\mathcal{P}_{\adm}=\left\{P\in\mathcal{P}~\big|~(P,\Supp(P))\text{ is admissible.}\right\}.$$

Now let us give several examples to understand the definition.

\begin{example}[Triangle $K_3$ and Cycle $C_3$]\label{K3} Let us consider the following digraphs $K_3$ and $C_3$:
\begin{figure}[H]
	\centering
	%% Creator: Inkscape 1.0.1 (3bc2e813f5, 2020-09-07), www.inkscape.org
%% Accompanies image file 'K3C3.pdf' (pdf, eps, ps)
\begingroup%
  \makeatletter%
  \providecommand\color[2][]{%
    \errmessage{(Inkscape) Color is used for the text in Inkscape, but the package 'color.sty' is not loaded}%
    \renewcommand\color[2][]{}%
  }%
  \providecommand\transparent[1]{%
    \errmessage{(Inkscape) Transparency is used (non-zero) for the text in Inkscape, but the package 'transparent.sty' is not loaded}%
    \renewcommand\transparent[1]{}%
  }%
  \providecommand\rotatebox[2]{#2}%
  \newcommand*\fsize{\dimexpr\f@size pt\relax}%
  \newcommand*\lineheight[1]{\fontsize{\fsize}{#1\fsize}\selectfont}%
  \ifx\svgwidth\undefined%
    \setlength{\unitlength}{260.05594569bp}%
    \ifx\svgscale\undefined%
      \relax%
    \else%
      \setlength{\unitlength}{\unitlength * \real{\svgscale}}%
    \fi%
  \else%
    \setlength{\unitlength}{\svgwidth}%
  \fi%
  \global\let\svgwidth\undefined%
  \global\let\svgscale\undefined%
  \makeatother%
  \begin{picture}(1,0.20052505)%
    \lineheight{1}%
    \setlength\tabcolsep{0pt}%
    \put(0,0){\includegraphics[width=\unitlength,page=1]{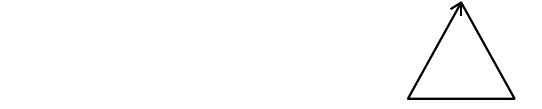}}%
    \put(0.70467791,0.00468185){\makebox(0,0)[lt]{\lineheight{1.25}\smash{\begin{tabular}[t]{l}$0$\end{tabular}}}}%
    \put(0.88694441,0.18591952){\makebox(0,0)[lt]{\lineheight{1.25}\smash{\begin{tabular}[t]{l}$1$\end{tabular}}}}%
    \put(0.97279564,0.00565701){\makebox(0,0)[lt]{\lineheight{1.25}\smash{\begin{tabular}[t]{l}$2$\end{tabular}}}}%
    \put(0.60697958,0.08817613){\makebox(0,0)[lt]{\lineheight{1.25}\smash{\begin{tabular}[t]{l}$C_3=$\end{tabular}}}}%
    \put(0,0){\includegraphics[width=\unitlength,page=2]{K3C3.pdf}}%
    \put(0.09571804,0.00189054){\makebox(0,0)[lt]{\lineheight{1.25}\smash{\begin{tabular}[t]{l}$0$\end{tabular}}}}%
    \put(0.27798458,0.18312822){\makebox(0,0)[lt]{\lineheight{1.25}\smash{\begin{tabular}[t]{l}$1$\end{tabular}}}}%
    \put(0.36383581,0.00286571){\makebox(0,0)[lt]{\lineheight{1.25}\smash{\begin{tabular}[t]{l}$2$\end{tabular}}}}%
    \put(-0.00198029,0.08538469){\makebox(0,0)[lt]{\lineheight{1.25}\smash{\begin{tabular}[t]{l}$K_3=$\end{tabular}}}}%
    \put(0,0){\includegraphics[width=\unitlength,page=3]{K3C3.pdf}}%
  \end{picture}%
\endgroup%

	%\caption[]{$\Supp(P)$}
	%\label{Fig:Q1}
\end{figure}
For the triangle digraph $K_3$, the set of minimal paths is
$$\mathcal{P}=\{e_0,~e_1,~e_2;~e_{01},~ e_{02}, ~e_{12};~ e_{012}\}.$$
Obviously, the admissible set $\mathcal{P}_{\adm}$ is the same as $\mathcal{P}$.

In particular, the admissible pair $(e_{012},K_3=\Supp(e_{012}))$ could be realized by the following digraph map $f: I^{\boxdot2}\rightarrow K_3$.
$$f((0,0))=0,\quad f((0,1))=1,\quad f((1,0))=f((1,1))=2,$$
or alternatively,
$$f((0,0))=f((1,0))=0,\quad f((0,1))=1,\quad f((1,1))=2.$$
For the 3-cycle digraph $C_3$, we have
$$\mathcal{P}_{\adm}=\mathcal{P}=\{e_0,~e_1,~e_2;~e_{01},~ e_{12},~ e_{20}\}.$$
\end{example}

\begin{example}[Square $S$ and Cycle $C_4$] Let us consider the following digraphs $S$ and $C_4$:
\begin{figure}[H]
	\centering
	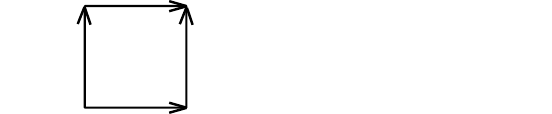
	%\caption[]{$\Supp(P)$}
	%\label{Fig:Q1}
\end{figure}
For the square digraph $S$, $S\cong I^{\boxdot 2}$. Obviously, we have
$$\mathcal{P}_{\adm}(S)=\mathcal{P}(S)=\{e_0,~e_1,~e_2,~e_3;~ e_{01},~ e_{02},~ e_{13},~ e_{23};~ e_{013}-e_{023}\}.$$
For the 4-cycle digraph $C_4$, we have
$$\mathcal{P}_{\adm}(C_4)=\mathcal{P}(C_4)=\{e_0,~e_1,~e_2,~e_3;~ e_{01},~ e_{13},~ e_{32},~ e_{20}\}.$$
\end{example}

\begin{example} For the supporting digraph $\Supp(e_{0123})$ (see Example \ref{exsuppdef}), the admissible set $\mathcal{P}_{\adm}$ is the same as $\mathcal{P}$. That is,
\begin{align*}
\mathcal{P}_{\adm}=\mathcal{P}=\{&e_0,~e_1,~e_2,~e_3;~ e_{01},~ e_{02},~ e_{12},~ e_{13}, ~e_{23};\\
 &e_{012},~e_{013}-e_{023},~e_{123};~e_{0123}\}
\end{align*}
In particular, the admissible pair $(e_{0123}, \Supp(e_{0123}))$ could be realized by $f:I^{\boxdot3}\longrightarrow\Supp(e_{0123})$:
\begin{align*}
&f((0,0,0))=0,\quad f((1,0,0))=1,\quad f((0,1,0))=2,\quad f((0,0,1))=2,\\
&f((1,1,0))=2,\quad f((1,0,1))=3,\quad f((0,1,1))=2,\quad f((1,1,1))=3.
\end{align*}
\end{example}

\begin{example}[Example \ref{EXfSfE1} revisited]\label{NE5} For the supporting digraph in Example \ref{EXfSfE1}, that is,
$$P=e_{0136}-e_{0156}+e_{0456}+e_{0246}-e_{0236}.$$
\begin{figure}[H]
	\centering
	
	%\caption[]{$\Supp(e_{0123})$}
	%\label{Fig:Q1}
\end{figure}
The admissible set is the same as the set of minimal paths. That is,
\begin{align*}
\mathcal{P}_{\adm}=\mathcal{P}=\{&e_0,~e_1,~e_2,~e_3,~e_4,~e_5,~ e_6;\\
              &e_{01},~e_{02},~e_{04},~e_{13},~e_{15},~e_{23},~e_{24},~e_{36},~e_{45},~e_{46},~e_{56};\\
              &e_{013}-e_{023},~e_{015}-e_{045},~e_{024},~e_{136}-e_{156},~e_{236}-e_{246},~e_{456};\\
              &e_{0136}-e_{0156}+e_{0456}+e_{0246}-e_{0236}\}
\end{align*}
In particular, the admissible pair $(P, \Supp(P))$ could be realized by $f:I^{\boxdot 3}\rightarrow \Supp(P)$:
\begin{align*}
&f((0,0,0))=0,\quad f((1,0,0))=1,\quad f((0,1,0))=2,\quad f((0,0,1))=4,\\
&f((1,1,0))=3,\quad f((1,0,1))=5,\quad f((0,1,1))=4,\quad f((1,1,1))=6.
\end{align*}
\end{example}

\begin{example}\label{NE32} For the following digraph $G=\Supp(e_{0135}-e_{0235}+e_{0245})$, we have $\mathcal{P}_{\adm}=\mathcal{P}$.
\begin{figure}[H]
	\centering
	%% Creator: Inkscape 1.0.1 (3bc2e813f5, 2020-09-07), www.inkscape.org
%% Accompanies image file 'EX315.pdf' (pdf, eps, ps)
\begingroup%
  \makeatletter%
  \providecommand\color[2][]{%
    \errmessage{(Inkscape) Color is used for the text in Inkscape, but the package 'color.sty' is not loaded}%
    \renewcommand\color[2][]{}%
  }%
  \providecommand\transparent[1]{%
    \errmessage{(Inkscape) Transparency is used (non-zero) for the text in Inkscape, but the package 'transparent.sty' is not loaded}%
    \renewcommand\transparent[1]{}%
  }%
  \providecommand\rotatebox[2]{#2}%
  \newcommand*\fsize{\dimexpr\f@size pt\relax}%
  \newcommand*\lineheight[1]{\fontsize{\fsize}{#1\fsize}\selectfont}%
  \ifx\svgwidth\undefined%
    \setlength{\unitlength}{202.50255194bp}%
    \ifx\svgscale\undefined%
      \relax%
    \else%
      \setlength{\unitlength}{\unitlength * \real{\svgscale}}%
    \fi%
  \else%
    \setlength{\unitlength}{\svgwidth}%
  \fi%
  \global\let\svgwidth\undefined%
  \global\let\svgscale\undefined%
  \makeatother%
  \begin{picture}(1,0.29128394)%
    \lineheight{1}%
    \setlength\tabcolsep{0pt}%
    \put(0,0){\includegraphics[width=\unitlength,page=1]{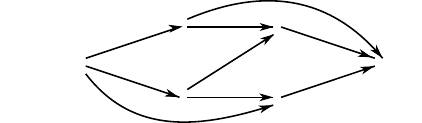}}%
    \put(0.13380566,0.12165521){\makebox(0,0)[lt]{\lineheight{1.25}\smash{\begin{tabular}[t]{l}$0$\end{tabular}}}}%
    \put(0.37790773,0.25379493){\makebox(0,0)[lt]{\lineheight{1.25}\smash{\begin{tabular}[t]{l}$1$\end{tabular}}}}%
    \put(0.39874329,0.09341465){\makebox(0,0)[lt]{\lineheight{1.25}\smash{\begin{tabular}[t]{l}$2$\end{tabular}}}}%
    \put(0.65527228,0.17228895){\makebox(0,0)[lt]{\lineheight{1.25}\smash{\begin{tabular}[t]{l}$3$\end{tabular}}}}%
    \put(0.65391179,0.00563322){\makebox(0,0)[lt]{\lineheight{1.25}\smash{\begin{tabular}[t]{l}$4$\end{tabular}}}}%
    \put(0.93471551,0.1222493){\makebox(0,0)[lt]{\lineheight{1.25}\smash{\begin{tabular}[t]{l}$5$\end{tabular}}}}%
    \put(0,0){\includegraphics[width=\unitlength,page=2]{EX315.pdf}}%
    \put(-0.00343812,0.12120995){\makebox(0,0)[lt]{\lineheight{1.25}\smash{\begin{tabular}[t]{l}$G=$\end{tabular}}}}%
  \end{picture}%
\endgroup%

	%\caption[]{$\Supp(e_{0123})$}
	%\label{Fig:Q1}
\end{figure}
That is,
\begin{align*}
\mathcal{P}_{\adm}=\mathcal{P}=\{&e_0,~e_1,~e_2,~e_3,~e_4,~e_5;\\
              &e_{01},~e_{02},~e_{04},~e_{13},~e_{15},~e_{23},~e_{24},~e_{35},~e_{45};\\
              &e_{013}-e_{023},~e_{015}-e_{045},~e_{024},~e_{135},~e_{235}-e_{245};\\
              &e_{0135}-e_{0235}+e_{0245}\}.
\end{align*}
In particular, we have the digraph map $f: I^{\boxdot 3}\longrightarrow \Supp(e_{0135}-e_{0235}+e_{0245})$:
\begin{align*}
&f((0,0,0))=0,\quad f((1,0,0))=1,\quad f((0,1,0))=2,\quad f((0,0,1))=4,\\
&f((1,1,0))=3,\quad f((1,0,1))=5,\quad f((0,1,1))=4,\quad f((1,1,1))=5.
\end{align*}
\end{example}

\begin{example}\label{NE46} For the following digraph $G=\Supp(e_{0136}-e_{0146}-e_{0236}+e_{0256})$, we have $\mathcal{P}_{\adm}=\mathcal{P}$.
\begin{figure}[H]
	\centering
	%% Creator: Inkscape 1.0.1 (3bc2e813f5, 2020-09-07), www.inkscape.org
%% Accompanies image file 'EX316.pdf' (pdf, eps, ps)
\begingroup%
  \makeatletter%
  \providecommand\color[2][]{%
    \errmessage{(Inkscape) Color is used for the text in Inkscape, but the package 'color.sty' is not loaded}%
    \renewcommand\color[2][]{}%
  }%
  \providecommand\transparent[1]{%
    \errmessage{(Inkscape) Transparency is used (non-zero) for the text in Inkscape, but the package 'transparent.sty' is not loaded}%
    \renewcommand\transparent[1]{}%
  }%
  \providecommand\rotatebox[2]{#2}%
  \newcommand*\fsize{\dimexpr\f@size pt\relax}%
  \newcommand*\lineheight[1]{\fontsize{\fsize}{#1\fsize}\selectfont}%
  \ifx\svgwidth\undefined%
    \setlength{\unitlength}{211.43943246bp}%
    \ifx\svgscale\undefined%
      \relax%
    \else%
      \setlength{\unitlength}{\unitlength * \real{\svgscale}}%
    \fi%
  \else%
    \setlength{\unitlength}{\svgwidth}%
  \fi%
  \global\let\svgwidth\undefined%
  \global\let\svgscale\undefined%
  \makeatother%
  \begin{picture}(1,0.35952127)%
    \lineheight{1}%
    \setlength\tabcolsep{0pt}%
    \put(0,0){\includegraphics[width=\unitlength,page=1]{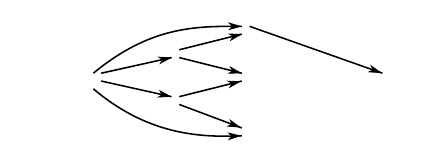}}%
    \put(0.58228007,0.16822231){\makebox(0,0)[lt]{\lineheight{1.25}\smash{\begin{tabular}[t]{l}$3$\end{tabular}}}}%
    \put(0,0){\includegraphics[width=\unitlength,page=2]{EX316.pdf}}%
    \put(0.16337957,0.18965666){\rotatebox{-180}{\makebox(0,0)[lt]{\lineheight{1.25}\smash{\begin{tabular}[t]{l}$0$\end{tabular}}}}}%
    \put(0.37246092,0.18113263){\makebox(0,0)[lt]{\lineheight{1.25}\smash{\begin{tabular}[t]{l}$1$\end{tabular}}}}%
    \put(0.38371257,0.08129551){\makebox(0,0)[lt]{\lineheight{1.25}\smash{\begin{tabular}[t]{l}$2$\end{tabular}}}}%
    \put(0.54624681,0.32217071){\makebox(0,0)[lt]{\lineheight{1.25}\smash{\begin{tabular}[t]{l}$4$\end{tabular}}}}%
    \put(0.54046794,0.00336046){\makebox(0,0)[lt]{\lineheight{1.25}\smash{\begin{tabular}[t]{l}$5$\end{tabular}}}}%
    \put(0.91185319,0.16126236){\makebox(0,0)[lt]{\lineheight{1.25}\smash{\begin{tabular}[t]{l}$6$\end{tabular}}}}%
    \put(0,0){\includegraphics[width=\unitlength,page=3]{EX316.pdf}}%
    \put(-0.00258299,0.16577442){\makebox(0,0)[lt]{\lineheight{1.25}\smash{\begin{tabular}[t]{l}$G=$\end{tabular}}}}%
  \end{picture}%
\endgroup%

	%\caption[]{$\Supp(e_{0123})$}
	%\label{Fig:Q1}
\end{figure}
That is,
\begin{align*}
\mathcal{P}_{\adm}=\mathcal{P}=\{&e_0,~e_1,~e_2,~e_3,~e_4,~e_5,~e_6;\\
              &e_{01},~e_{02},~e_{04},~e_{05},~e_{13},~e_{14},~e_{23},~e_{25},~e_{36},~e_{46},~e_{56};\\
              &e_{013}-e_{023},~e_{046}-e_{056},~e_{014},~e_{025},~e_{136}-e_{146},~e_{236}-e_{256};\\
              &e_{0136}-e_{0146}-e_{0236}+e_{0256}\}.
\end{align*}
In particular, we have the digraph map $f: I^{\boxdot 3}\longrightarrow \Supp(e_{0136}-e_{0146}-e_{0236}+e_{0256})$:
\begin{align*}
&f((0,0,0))=0,\quad f((1,0,0))=1,\quad f((0,1,0))=0,\quad f((0,0,1))=2,\\
&f((1,1,0))=4,\quad f((1,0,1))=3,\quad f((0,1,1))=5,\quad f((1,1,1))=6.
\end{align*}
\end{example}

In the above examples, $\mathcal{P}_{\adm}=\mathcal{P}$. In particular, we see that there are many different admissible pairs. In what follows, we will give several examples in which $\mathcal{P}_{\adm}$ is a proper subset of $\mathcal{P}$.

\begin{example}[exotic cube]\label{exotic} The minimal digraph which supports
$$P=e_{0158}-e_{0258}+e_{0268}-e_{0368}+e_{0378}-e_{0478}$$
is given by the following digraph, which is known as the exotic cube.
\begin{figure}[H]
	\centering
	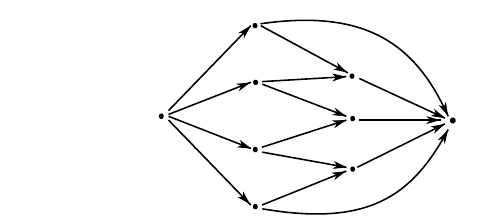
	%\caption[]{$\Supp(e_{0123})$}
	%\label{Fig:Q1}
\end{figure}
The set of minimal paths is given by
\begin{align*}
\mathcal{P}=\{&e_0, ~e_1,~e_2,~e_3,~e_4,~e_5,~e_6,~e_7,~e_8; \\
&e_{01}, ~e_{02},~ e_{03}, ~e_{04}, ~e_{15},~ e_{18},~ e_{25},~ e_{26},~ e_{36}, ~e_{37},~ e_{47},~ e_{48},~ e_{58},~ e_{68},~ e_{78};\\
&e_{015}-e_{025},~ e_{026}-e_{036},~ e_{037}-e_{047},~ e_{018}-e_{048},~ e_{158},~ e_{258}-e_{268},~ e_{368}-e_{378},~ e_{478};\\
&e_{0158}-e_{0258}+e_{0268}-e_{0368}+e_{0378}-e_{0478}\}.
\end{align*}

Note that the exotic cube has one more vertex than $I^{\boxdot 3}$, thus it can not be an image digraph of $I^{\boxdot 3}$. Thus,
$$\mathcal{P}_{\adm}=\mathcal{P}\setminus\{P\}.$$
\end{example}

\begin{remark}
One can find more examples as the exotic cube, such as, Example \ref{more}. In fact, although there are various kinds of admissible pairs, for each $n\in\mathbb{N}_+$, the total numbers are finite, due to the second condition in Definition \ref{admpair}. In particular, we conjecture that each supporting digraph of the admissible path is contractible in Subsection \ref{properties}.
\end{remark}

Let us further study the supporting digraphs of two minimal $4$-paths, $\Supp(e_{01234})$ and Example \ref{length4}.

\begin{example} The digraph $\Supp(e_{01234})$ is as follows:
\begin{figure}[H]
	\centering
	%% Creator: Inkscape 1.0.1 (3bc2e813f5, 2020-09-07), www.inkscape.org
%% Accompanies image file '01234.pdf' (pdf, eps, ps)
\begingroup%
  \makeatletter%
  \providecommand\color[2][]{%
    \errmessage{(Inkscape) Color is used for the text in Inkscape, but the package 'color.sty' is not loaded}%
    \renewcommand\color[2][]{}%
  }%
  \providecommand\transparent[1]{%
    \errmessage{(Inkscape) Transparency is used (non-zero) for the text in Inkscape, but the package 'transparent.sty' is not loaded}%
    \renewcommand\transparent[1]{}%
  }%
  \providecommand\rotatebox[2]{#2}%
  \newcommand*\fsize{\dimexpr\f@size pt\relax}%
  \newcommand*\lineheight[1]{\fontsize{\fsize}{#1\fsize}\selectfont}%
  \ifx\svgwidth\undefined%
    \setlength{\unitlength}{251.93884422bp}%
    \ifx\svgscale\undefined%
      \relax%
    \else%
      \setlength{\unitlength}{\unitlength * \real{\svgscale}}%
    \fi%
  \else%
    \setlength{\unitlength}{\svgwidth}%
  \fi%
  \global\let\svgwidth\undefined%
  \global\let\svgscale\undefined%
  \makeatother%
  \begin{picture}(1,0.20999051)%
    \lineheight{1}%
    \setlength\tabcolsep{0pt}%
    \put(0,0){\includegraphics[width=\unitlength,page=1]{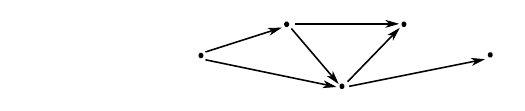}}%
    \put(0.32035071,0.09072258){\makebox(0,0)[lt]{\lineheight{1.25}\smash{\begin{tabular}[t]{l}$0$\end{tabular}}}}%
    \put(0.53074603,0.19357358){\makebox(0,0)[lt]{\lineheight{1.25}\smash{\begin{tabular}[t]{l}$1$\end{tabular}}}}%
    \put(0.63694939,0.00195147){\makebox(0,0)[lt]{\lineheight{1.25}\smash{\begin{tabular}[t]{l}$2$\end{tabular}}}}%
    \put(0.7636693,0.19491444){\makebox(0,0)[lt]{\lineheight{1.25}\smash{\begin{tabular}[t]{l}$3$\end{tabular}}}}%
    \put(0.97191917,0.09353856){\makebox(0,0)[lt]{\lineheight{1.25}\smash{\begin{tabular}[t]{l}$4$\end{tabular}}}}%
    \put(0,0){\includegraphics[width=\unitlength,page=2]{01234.pdf}}%
    \put(-0.00110383,0.09200434){\makebox(0,0)[lt]{\lineheight{1.25}\smash{\begin{tabular}[t]{l}$\Supp(e_{01234})=$\end{tabular}}}}%
  \end{picture}%
\endgroup%

	%\caption[]{$\Supp(e_{0123})$}
	%\label{Fig:Q1}
\end{figure}
Similarly, we have
\begin{align*}
\mathcal{P}_{\adm}=\mathcal{P}=\{&e_0,~e_1,~e_2,~e_3,~e_4;~e_{01},~e_{02},~e_{12},~e_{13},~e_{23},~e_{24},~e_{34};\\
&e_{012},~e_{013}-e_{023},~e_{123},~e_{124}-e_{134},~e_{234};\\
&e_{0123},~e_{1234},~e_{0124}-e_{0134}+e_{0234};~e_{01234}\}.
\end{align*}
In particular, we have the two digraph maps realizing $e_{0124}-e_{0134}+e_{0234}$ and $e_{01234}$ respectively:
\begin{align*}
&f:I^{\boxdot 4}\longrightarrow \Supp(e_{01234}),\\
&g:I^{\boxdot 3}\longrightarrow \Supp(e_{0124}-e_{0134}+e_{0234}),
\end{align*}
which are given by
\begin{align*}
&f((0,0,0,0))=0,\quad f((1,0,0,0))=1,\quad f((0,1,0,0))=1,\quad f((0,0,1,0))=2,\\
&f((0,0,0,1))=2,\quad f((1,1,0,0))=3,\quad f((1,0,1,0))=2,\quad f((1,0,0,1))=3,\\
&f((0,1,1,0))=3,\quad f((0,1,0,1))=2,\quad f((0,0,1,1))=3,\quad f((1,1,1,0))=4,\\
&f((1,1,0,1))=3,\quad f((1,0,1,1))=4,\quad f((0,1,1,1))=3,\quad f((1,1,1,1))=4;\\
&g=f\big|_{I^{\boxdot3}\boxdot\{0\}}.
\end{align*}
Note that in this example,
$$\Supp(e_{0124}-e_{0134}+e_{0234})=\Supp(e_{01234}),$$
but we regard $\left(e_{0124}-e_{0134}+e_{0234},~ \Supp(e_{0124}-e_{0134}+e_{0234})\right)$ as a length $3$ admissible pair, while
$\left(e_{01234}, ~\Supp(e_{01234})\right)$ as a length $4$ admissible pair. The former path is a minimal face component of the latter one. %Moreover, one can compute that the only non-trivial CW homology is $H_0^{\CW}(\Supp(e_{01234});\R)\cong\R$.
\end{example}

\begin{remark} In general, the existence of the digraph map $f:I^{\boxdot n}\longrightarrow \Supp(P_{n,\alpha})$ in our definition is not unique, such as the singular cubical realizations of $(e_{012},K_3)$ in Example \ref{K3}, there may exist different ``degenerations" on faces. Meanwhile, $I^{\boxdot n}$ has a non-trivial automorphism group, which also leads to the non-uniqueness. %For example, there are many digraph maps realizing $\Supp(e_{012\cdots n})$ from $I^{\boxdot n}$.
\end{remark}

\begin{example}[Example \ref{length4} revisited]\label{length4rev} Let us look at the digraph $\Supp(P)$ in Example \ref{length4}. According to the structure theorem \ref{structurethm}, we have
\begin{align*}
\mathcal{P}_3=\{& e_{S159}-e_{S169}+e_{S269},~e_{S16(10)}-e_{S26(10)}+e_{S27(10)}-e_{S37(10)},\\
&e_{S27(11)}-e_{S37(11)}+e_{S28(11)}-e_{S48(11)};\\
&e_{159E}-e_{169E}+e_{16(10)E},~e_{269E}-e_{26(10)E}+e_{27(10)E}-e_{27(11)E},\\
&e_{37(10)E}-e_{37(11)E}+e_{38(11)E},~e_{48(11)E};\\
&P_{S,3,E}:=e_{S15E}-e_{S1(10)E}+e_{S2(11)E}-e_{S29E}-e_{S38E}+e_{S3(10)E}+e_{S48E}-e_{S4(11)E}+e_{S59E}\}.
\end{align*}
According to the above examples, all the minimal $3$-paths in $\mathcal{P}_3$ are admissible except the last one $P_{S,3,E}$, since it has 11 vertices. Thus,
$$\mathcal{P}_{\adm,3}=\mathcal{P}_3\setminus\{P_{S,3,E}\}.$$

Furthermore, $P\in\mathcal{P}_4$ is also not admissible. The reason is as follows: if $(P,\Supp(P))$ is realizable by the singular $4$-cube $f$, then $f_*(\omega_n)=P$ implies
\begin{align*}
& f((0,0,0,0))=S,\quad f((1,1,1,1))=E; \\
& f:\{(1,0,0,0),(0,1,0,0),(0,0,1,0),(0,0,0,1)\}\rightarrowtail \{1,2,3,4\}.
\end{align*}
We claim that the edge $S\rightarrow 5$ can not be in the image digraph of $f$. If $S\rightarrow 5$ is realizable by $f$, let us analyze the possible pre-images of $S$:
\begin{itemize}
  \item $f((0,0,0,0))=S$. Note that there are only four directed edges starting from $(0,0,0,0)$ in $I^{\boxdot 4}$, which have been mapped to $S\rightarrow 1$, $S\rightarrow 2$, $S\rightarrow 3$ and $S\rightarrow 4$. Then $S\rightarrow 5$ can not be realized as the pre-image of $f$ via the edges starting from $(0,0,0,0)$.
  \item $f((\epsilon_1,\epsilon_2,\epsilon_3,\epsilon_4))=S$, for some $\epsilon_i$ with two of them being $1$ and the other two of them being $0$. Then the image digraph has the edge $1\rightarrow S$, or $2\rightarrow S$, or $3\rightarrow S$, or $4\rightarrow S$. Contraction.
  \item $f((\epsilon_1,\epsilon_2,\epsilon_3,\epsilon_4))=S$, for some $\epsilon_i$ with three of them being $1$ and the remaining of them being $0$. Then the image digraph has the edge $S\rightarrow E$. Contraction.
\end{itemize}
Thus, $\Supp(P)$ could not be the image digraph of any singular $4$-cube. We have
$$\mathcal{P}_{\adm}=\mathcal{P}\setminus\{P, P_{S,3,E}\}.$$
\end{example}

\subsubsection{Admissible relations}

\begin{definition}\label{admrel} For $P_1$, $P_2\in\mathcal{P}_{\adm}(G)$, if one of $\pm P_1\pm P_2$ belongs to $\mathcal{P}_{\adm}(G)$, denote it by $P_3$, that is,
$$(P_3,\Supp(P_3))\in\mathcal{P}_{\adm}(G),$$
%if $P_3:=P_1+P_2$ (resp. $P_4:=P_1-P_2$) $\in\mathcal{P}_{\adm}(G)$, more explicitly,
%$$(P_1+P_2,\Supp(P_1+P_2))~(\text{resp. } (P_1-P_2,\Supp(P_1-P_2))) \in\mathcal{P}_{\adm}(G),$$
%$$(P_3,\Supp(P_3))~(\text{resp. } (P_4,\Supp(P_4))) \in\mathcal{P}_{\adm}(G),$$
we say that $P_1$ and $P_2$ are \textbf{summable} in $\mathcal{P}_{\adm}(G)$. If $P_3=P_1+P_2$, we also write the relation as
$$(P_3,\Supp(P_3))=(P_1,\Supp(P_1))+(P_2,\Supp(P_2)).$$
We call the linear relations among $\mathcal{P}_{\adm}(G)$ the \textbf{admissible relations}.

%\begin{itemize}
%  \item $P_1$ and $P_2$ are summable in $\mathcal{P}_{\adm}(G)$;
%  \item %$P_1$, $P_2$ and $P_1+P_2$ (or $P_1-P_2$) are linearly independent in $\mathcal{P}_{\adm}$.
%\end{itemize}

\end{definition}

\begin{example}\label{k-square} For any integer $k\geq 3$, we consider the following digraph $S_k$, with
$$V(S_k)=\{S,E,1,2,\ldots,k\},\qquad E(S_k)=\{S~\rightarrow~i\}_{i=1}^k\cup\{i~\rightarrow~E\}_{i=1}^k.$$
\begin{figure}[H]
	\centering
	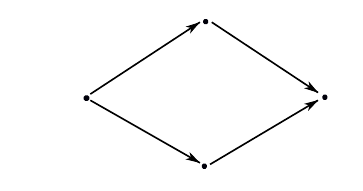
	%\caption[]{$\Supp(e_{0123})$}
	%\label{Fig:Q1}
\end{figure}
It is easy to see that
\begin{align*}
\mathcal{P}_{\adm}=\mathcal{P}=\{&e_S,~e_E,~e_i,~i=1,\ldots,k;\\
              &e_{Si}, ~e_{iE}, ~i=1,\ldots,k;\\
              &e_{SiE}-e_{SjE},~i,j=1,\ldots,k, \text{and }i\neq j\text{ (up to a sign)}\}
\end{align*}
For distinct $a,b,c\in\{1,\ldots,k\}$, $e_{SaE}-e_{SbE}$ and $e_{SbE}-e_{ScE}$ are summable and
$$(e_{SaE}-e_{SbE})+(e_{SbE}-e_{ScE})=e_{SaE}-e_{ScE}\in\mathcal{P}_{\adm}.$$
\end{example}

\begin{example}\label{summable} For the following digraph,
\begin{figure}[H]
	\centering
	%% Creator: Inkscape 1.0.1 (3bc2e813f5, 2020-09-07), www.inkscape.org
%% Accompanies image file 'minrel1.pdf' (pdf, eps, ps)
\begingroup%
  \makeatletter%
  \providecommand\color[2][]{%
    \errmessage{(Inkscape) Color is used for the text in Inkscape, but the package 'color.sty' is not loaded}%
    \renewcommand\color[2][]{}%
  }%
  \providecommand\transparent[1]{%
    \errmessage{(Inkscape) Transparency is used (non-zero) for the text in Inkscape, but the package 'transparent.sty' is not loaded}%
    \renewcommand\transparent[1]{}%
  }%
  \providecommand\rotatebox[2]{#2}%
  \newcommand*\fsize{\dimexpr\f@size pt\relax}%
  \newcommand*\lineheight[1]{\fontsize{\fsize}{#1\fsize}\selectfont}%
  \ifx\svgwidth\undefined%
    \setlength{\unitlength}{170.88345625bp}%
    \ifx\svgscale\undefined%
      \relax%
    \else%
      \setlength{\unitlength}{\unitlength * \real{\svgscale}}%
    \fi%
  \else%
    \setlength{\unitlength}{\svgwidth}%
  \fi%
  \global\let\svgwidth\undefined%
  \global\let\svgscale\undefined%
  \makeatother%
  \begin{picture}(1,0.43487891)%
    \lineheight{1}%
    \setlength\tabcolsep{0pt}%
    \put(0,0){\includegraphics[width=\unitlength,page=1]{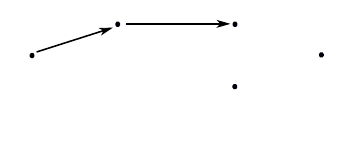}}%
    \put(-0.00202826,0.2590384){\makebox(0,0)[lt]{\lineheight{1.25}\smash{\begin{tabular}[t]{l}$0$\end{tabular}}}}%
    \put(0.30816419,0.41067485){\makebox(0,0)[lt]{\lineheight{1.25}\smash{\begin{tabular}[t]{l}$1$\end{tabular}}}}%
    \put(0.65157022,0.41265173){\makebox(0,0)[lt]{\lineheight{1.25}\smash{\begin{tabular}[t]{l}$3$\end{tabular}}}}%
    \put(0.95859955,0.26319008){\makebox(0,0)[lt]{\lineheight{1.25}\smash{\begin{tabular}[t]{l}$6$\end{tabular}}}}%
    \put(0,0){\includegraphics[width=\unitlength,page=2]{minrel1.pdf}}%
    \put(0.30905074,0.23649841){\makebox(0,0)[lt]{\lineheight{1.25}\smash{\begin{tabular}[t]{l}$2$\end{tabular}}}}%
    \put(0.64886866,0.23706708){\makebox(0,0)[lt]{\lineheight{1.25}\smash{\begin{tabular}[t]{l}$4$\end{tabular}}}}%
    \put(0.47053018,0.07157103){\makebox(0,0)[lt]{\lineheight{1.25}\smash{\begin{tabular}[t]{l}$5$\end{tabular}}}}%
  \end{picture}%
\endgroup%

	%\caption[]{$\Supp(e_{0123})$}
	%\label{Fig:}
\end{figure}
The admissible set is the same as the set of minimal paths. In particular, the length 3 part is given by
$$\mathcal{P}_{\adm,3}=\mathcal{P}_3=\{e_{0136}-e_{0236}+e_{0256},~e_{0136}-e_{0146}+e_{0246}-e_{0236},~e_{0146}-e_{0246}+e_{0256}\}.$$
Clearly, any two of them are summable in $\mathcal{P}_{\adm,3}$.
\end{example}

\begin{example} For the following digraph,
\begin{figure}[H]
	\centering
	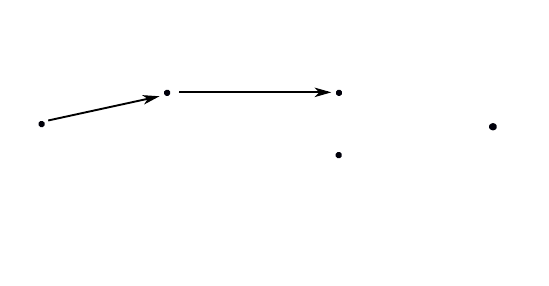
	%\caption[]{$\Supp(e_{0123})$}
	%\label{Fig:}
\end{figure}
Similarly, we look at the length 3 part.
\begin{align*}
\mathcal{P}_3=\{& e_{S05E}-e_{S15E}+e_{S16E},~e_{S16E}-e_{S17E}+e_{S27E},\\
                & e_{S27E}-e_{S37E}+e_{S38E},~e_{S38E}-e_{S39E}+e_{S49E},\\
                & e_{S05E}-e_{S15E}+e_{S17E}-e_{S27E},~e_{S16E}-e_{S17E}+e_{S37E}-e_{S38E},\\
                & e_{S27E}-e_{S37E}+e_{S39E}-e_{S49E},\\
                & e_{S05E}-e_{S15E}+e_{S17E}-e_{S37E}+e_{S38E},\\
                & e_{S16E}-e_{S17E}+e_{S37E}-e_{S39E}+e_{S49E},\\
                & e_{S05E}-e_{S15E}+e_{S17E}-e_{S37E}+e_{S39E}-e_{S49E}
\},
\end{align*}
Note that $e_{S05E}-e_{S15E}+e_{S17E}-e_{S37E}+e_{S38E}$ is not admissible, since its supporting digraph has 13 edges, more than those of 3-cube $I^{\boxdot 3}$. One can apply the same argument to the final two paths in $\mathcal{P}_3$. Then we have
\begin{align*}
\mathcal{P}_{\adm,3}=\{& e_{S05E}-e_{S15E}+e_{S16E},~e_{S16E}-e_{S17E}+e_{S27E},\\
                       &e_{S27E}-e_{S37E}+e_{S38E},~e_{S38E}-e_{S39E}+e_{S49E},\\
                       &e_{S05E}-e_{S15E}+e_{S17E}-e_{S27E},~e_{S16E}-e_{S17E}+e_{S37E}-e_{S38E},\\
                       &e_{S27E}-e_{S37E}+e_{S39E}-e_{S49E}
\}.
\end{align*}
In particular, note that
\begin{itemize}
  \item $e_{S05E}-e_{S15E}+e_{S16E}$ and $e_{S16E}-e_{S17E}+e_{S27E}$ are summable in both $\mathcal{P}_3$ and $\mathcal{P}_{\adm,3}$;
  \item $e_{S05E}-e_{S15E}+e_{S16E}$ and $e_{S16E}-e_{S17E}+e_{S37E}-e_{S38E}$ is summable in $\mathcal{P}_3$, but not summable in $\mathcal{P}_{\adm,3}$.
\end{itemize}
\end{example}

\subsubsection{CW complex of the digraph}

Now we follow Whitehead's definition of CW complex of a Hausdorff topological space, and give our definition of CW complex of a digraph $G$.

\begin{definition}\label{CWG} Let $G$ be a digraph, and $\mathcal{P}_{\adm}=\bigcup_{n\geq0}\mathcal{P}_{\adm,n}=\bigcup_{n\geq0}\{P_{n,\alpha}\}_{\alpha}$ be the set of admissible minimal paths in $G$. Then the pair $\left(G,\mathcal{P}'_{\adm}=\bigcup_{n\geq0}\{P_{n,\alpha'}\}_{\alpha'}\subset\mathcal{P}_{\adm}\right)$ is called a CW complex of $G$ if the following three conditions are satisfied:
\begin{enumerate}
  \item $G=\bigcup_{P_{n,\alpha'}\in\mathcal{P}'_{\adm}}\Supp(P_{n,\alpha'})$;
  \item For each $n\geq0$, there are not admissible relations among $\{P_{n,\alpha'}\}_{\alpha'}$;
  \item $\mathcal{P}'_{\adm}$ is the maximal subset of $\mathcal{P}_{\adm}$ satisfying (1) (2).
\end{enumerate}
\end{definition}

\begin{remark}By Definition \ref{CWG} (2) (3), the $\CW$ complex associated to $G$ depends on the choice of the basis due to the possible existence of admissible relations.
\end{remark}

\subsection{Cellular homology of a digraph}\label{3.3}

By the standard argument, there exists a cellular chain complex associated to each $\CW$ complex. But, to circumvent the intricacies involved in the choice of a basis, we define the associated cellular chain complex to be a quotient complex as follows. For simplicity, in the sequel, we work with coefficient $\R$ or other fields with characteristic $0$.

\begin{definition}[Cellular chain of a digraph $G$]\label{cellular} Let $G$ be a digraph and $\mathcal{P}_{\adm}$ be the set of admissible minimal paths of $G$. For all $n\in\mathbb{N}$, the $\R$-module $C_n(G;\R)$, called the $n$-th cellular chain group of $G$ with coefficient $\R$, is generated by the admissible $n$-minimal paths quotient by the admissible relations
$$C_n(G;\R)=\left\{\sum_{P\in\mathcal{P}_{\adm,n}}c_PP~\bigg|~c_P\in \R\right\}\bigg/\left\{\text{admissible relations}\right\}.$$
Moreover, %we define an operator $\partial^{\cell}:C_n(G;K)\longrightarrow C_{n-1}(G;K)$. For any $P\in\mathcal{P}_{\adm,n}$,
for any $P\in\mathcal{P}_{\adm,n}$, and $n\geq1$, we define $\partial^{\cell}P$ as follows:
$$\partial^{\cell}P=\sum_{Q\in\mathcal{P}_{\adm,n-1}}[Q:P]Q,$$
where the number $[Q:P]$ is defined as follows:
\begin{equation*}
[Q:P]=\left\{
        \begin{array}{ll}
          1, & \hbox{if $Q<\partial P$;} \\
          -1, & \hbox{if $-Q<\partial P$;} \\
          0, & \hbox{otherwise.}
        \end{array}
      \right.
\end{equation*}
\end{definition}
We also call $Q$ to be the real face of $P$ if $[Q:P]\neq0$. For $P\in\mathcal{P}_{\adm,0}$, $\partial^{\cell}P:=0$; or we make the convention that $\mathcal{P}_{\adm,-1}=0$, $C_{-1}(G;\R)=0$.

It extends $\R$-linearly to the graded free group generated by $P\in\mathcal{P}_{\adm}(G)$.

\begin{proposition}\label{welldef} The operator $\partial^{\cell}$ can be descended to the quotient group $C_*(G;\R)$, still denoted by $\partial^{\cell}$. That is, we have a well-defined operator
$$\partial^{\cell}:C_n(G;\R)\longrightarrow C_{n-1}(G;\R).$$
Moreover, it is a boundary operator: $\left(\partial^{\cell}\right)^2=0$.
\end{proposition}

The proof is based on the following several lemmas.

\begin{lemma}\label{elementarylem} For $P\in\mathcal{P}_{\adm,n}(G)$. In the sub-digraph $\Supp(P)\subset G$, the boundary operator
$$\partial^{\cell}: C_n(\Supp(P);\R)\longrightarrow C_{n-1}(\Supp(P);\R)$$
coincides with $\partial$ of path complex of $\Supp(P)$. That is, we have
$$\partial^{\cell}P=\partial P=\sum_{\alpha\in E_1}P_{S,n-1,\alpha}+\sum_{\beta\in S_1}P_{\beta,n-1,E}+\sum_{k\in I_P}P_{S,n-1,E}^k.$$
\end{lemma}

\begin{proof} For any $P\in \mathcal{P}_{\adm,n}$ with starting vertex $S$ and ending vertex $E$, by the structure theorem \ref{structurethm} for $\Supp(P)$, we have the decomposition in $\Supp(P)$,
\begin{equation}\label{decomp1}
\partial P=\sum_{\alpha\in E_1}P_{S,n-1,\alpha}+\sum_{\beta\in S_1}P_{\beta,n-1,E}+\sum_{k\in I_P}P_{S,n-1,E}^k.
\end{equation}
And the minimal face components of right hand side are all the real faces of $P$ in $\Supp(P)$ with the right sign. It remains to check that each of them has a singular cubical realization, i.e.
$$P_{S,n-1,\alpha},~P_{\beta,n-1,E},~P_{S,n-1,E}^k\in\mathcal{P}_{\adm,n-1}(\Supp(P)).$$
Since $P$ is admissible, there exists a digraph map $f:I^{\boxdot n}\rightarrow \Supp(P)$, such that
$$f_*(\omega_n)=cP,\quad\text{for some }c\in\Z\setminus\{0\}.$$
In particular, by the length reason, there must be
$$f((0,0,\ldots,0))=S,\quad f((1,1,\ldots,1))=E.$$
By Proposition \ref{singularcubpath}\footnote{This proposition works here since $\Supp(P)$ is acyclic.}, we have
\begin{equation}\label{decomp2}
\begin{aligned}
\partial P=&~\partial\circ f_*(\omega_n)=\partial\circ\tau_n(f)=\tau_{n-1}\circ\partial^c(f)=(\partial^c(f))_*(\omega_{n-1})\\
          =&~\sum_{j=1}^n(-1)^j(f\circ F_{j0})_*(\omega_{n-1})-\sum_{j=1}^n(-1)^j(f\circ F_{j1})_*(\omega_{n-1}).
\end{aligned}
\end{equation}
Now let us compare the two decompositions \eqref{decomp1} and \eqref{decomp2}. First, by definition of the induced homomorphism $(f\circ F_{j\epsilon})_*$, each $(f\circ F_{j\epsilon})_*(\omega_{n-1})\in\Omega_{n-1}(\Supp(P);\Z)$. Then we can classify it according to the starting and ending vertices.

More explicitly, note that for the minimal path $\omega_n$ in $I^{\boxdot n}$, we have
\begin{align*}
&S(\omega_n)=(0,0,\ldots,0),\qquad E(\omega_n)=(1,1,\ldots,1),\\
&S_1(\omega_n)=\{(1,0,0,\ldots,0),~(0,1,0,\ldots,0),~\cdots,~(0,\ldots,0,0,1)\}\\
&E_1(\omega_n)=\{(0,1,1,\ldots,1),~(1,0,1,\ldots,1),~\cdots,~(1,\ldots,1,1,0)\}.
\end{align*}
Furthermore, let
\begin{align*}
&E_{1,\alpha}=\{j_{\alpha}|0 \text{ appears in the $j_{\alpha}$-th position of } f^{-1}(\alpha)\cap E_1(\omega_n)\},\\
&S_{1,\beta}=\{k_{\beta}|1 \text{ appears in the $k_{\beta}$-th position of }f^{-1}(\beta)\cap S_1(\omega_n)\}.
\end{align*}
Then by \eqref{decomp1}, we have
\begin{align*}
&P_{S,n-1,\alpha}=\sum_{j_{\alpha}\in E_{1,\alpha}}(-1)^{j_{\alpha}}(f\circ F_{j_{\alpha}0})_*(\omega_{n-1}),\\
&P_{\beta,n-1,E}=\sum_{k_{\beta}\in S_{1,\beta}}(-1)^{k_{\beta}}(f\circ F_{k_{\beta}1})_*(\omega_{n-1}),\\
&P_{S,n-1,E}^k=\partial P-\left(\sum_{\alpha\in E_1}\sum_{j_{\alpha}\in E_{1,\alpha}}(-1)^{j_{\alpha}}(f\circ F_{j_{\alpha}0})_*(\omega_{n-1})-\right.\\
&\qquad\qquad\qquad\left.\sum_{\beta\in S_1}\sum_{k_{\beta}\in S_{1,\beta}}(-1)^{k_{\beta}}(f\circ F_{k_{\beta}1})_*(\omega_{n-1})\right).
\end{align*}
Since $\Supp(P)$ is an acyclic digraph, then each $(f\circ F_{j\epsilon})_*(\omega_{n-1})\in\Omega_{n-1}(\Supp(P))$. By the decomposition \eqref{decomp1}, and the unique result in the structure theorem, there exist non-zero integers $\{c_{\alpha}\}_{\alpha\in E_1}$, $\{c_{\beta}\}_{\beta\in S_1}$, $\{d_k\}_{k\in I_P}$ such that
\begin{align*}
&c_{\alpha} P_{S,n-1,\alpha}=(f\circ F_{j_{\alpha}0})_*(\omega_{n-1}),\text{ for some }j_{\alpha},\\
&c_{\beta} P_{\beta,n-1,E}=(f\circ F_{k_{\beta}1})_*(\omega_{n-1}),\text{ for some }k_{\beta},\\
&d_k P_{S,n-1,E}^k=(f\circ F_{j\epsilon})_*(\omega_{n-1}),\text{ for some } j \text{ and some }\epsilon.
\end{align*}

It remains to check the images of the above chosen $f\circ F_{j_{\alpha}0}$, $f\circ F_{k_{\beta}1}$ and $f\circ F_{j\epsilon}$ are $\Supp(P_{S,n-1},\alpha)$, $\Supp(P_{\beta,n-1,E}$ and $\Supp(P_{S,n-1,E}^k)$ respectively.

Since $(f\circ F_{j_{\alpha}0})_*(\omega_{n-1})=c_{\alpha}P_{S,n-1,\alpha}\in\Omega_{n-1}(\Supp(P))$, by the minimal condition of $\Supp(P_{S,n-1,\alpha})$, we have
$$\Supp(P_{S,n-1,\alpha})\subset\im(f\circ F_{j_{\alpha}0}).$$
Now assume $\Supp(P_{S,n-1,\alpha})\neq\im(f\circ F_{j_{\alpha}0})$, then applying the structure theorem to $\Supp(P)$ and $\Supp(P_{S,n-a,\alpha})$, we see that there exists
$$v\in V(\im(f\circ F_{j_{\alpha}0}))\setminus V(\Supp(P_{S,n-1,\alpha})).$$

Let us choose one of the pre-images $(f\circ F_{j_{\alpha}0})^{-1}(v)$ in $I^{\boxdot (n-1)}$, and denote it by
$$(\epsilon_1^v, \epsilon_2^v, \ldots, \epsilon_{n-1}^v).$$
Then we consider the elementary $(n-1)$-path in $I^{\boxdot(n-1)}$ passing through $(\epsilon_1^v, \epsilon_2^v, \ldots, \epsilon_{n-1}^v)$. Its image under $f\circ F_{j_{\alpha}0}$ is an elementary $k$-path ($k\leq (n-1)$) from $S$ to $\alpha$ passing through $v$. By applying the structure theorem of $\Supp(P)$ and $\Supp(P_{S,n-1,\alpha})$ again, we see that all the elementary $k$-paths ($k\leq (n-1)$) from $S$ to $\alpha$ in $\Supp(P)$ lie in $\Supp(P_{S,n-1,\alpha})$. Thus, $v\in V(\Supp(P_{S,n-1,\alpha}))$, contradiction.

The arguments for $\Supp(P_{\beta,n-1,E})$ and $\Supp(P_{S,n-1,E}^k)$ are the same. In summary,
$$P_{S,n-1,\alpha},~P_{\beta,n-1,E},~P_{S,n-1,E}^k\in\mathcal{P}_{\adm,n-1}(\Supp(P)).$$
\end{proof}

Note that each face component of $P$ in $\Supp(P)$ is a linear combination of s-regular allowed elementary paths with coefficients $1$ or $-1$. But it may not be minimal in $G$. Thus, to prove the coincidence of $\partial^{\cell}$ and $\partial$ in $G$, first we need to prove if some face component of $P$ is not minimal in $G$, then any of its minimal components have singular cubical realizations. More explicitly, we need to prove the following lemma.
\begin{lemma}\label{Keylemma} Let $P$ be a linear combination of s-regular allowed elementary paths with coefficients $1$ or $-1$ and there exists a digraph map
$$f:I^{\boxdot n}\longrightarrow G$$
such that $f_*(\omega_n)=cP$ for some $c\in\mathbb{Z}\setminus\{0\}$, and the image digraph is $\Supp(P)$. Then for any $P'\in\mathcal{P}_n(G)$ with $P'\leq P$, we have $P'\in\mathcal{P}_{\adm,n}(G)$.
\end{lemma}

Before the proof, let us give two examples to explain it.

\begin{example}\label{basicex} Let $G$ be the following digraph and $P=e_{013}-e_{023}\in\Omega_2(G)$,
\begin{figure}[H]
	\centering
	%% Creator: Inkscape 1.0.1 (3bc2e813f5, 2020-09-07), www.inkscape.org
%% Accompanies image file 'KS.pdf' (pdf, eps, ps)
\begingroup%
  \makeatletter%
  \providecommand\color[2][]{%
    \errmessage{(Inkscape) Color is used for the text in Inkscape, but the package 'color.sty' is not loaded}%
    \renewcommand\color[2][]{}%
  }%
  \providecommand\transparent[1]{%
    \errmessage{(Inkscape) Transparency is used (non-zero) for the text in Inkscape, but the package 'transparent.sty' is not loaded}%
    \renewcommand\transparent[1]{}%
  }%
  \providecommand\rotatebox[2]{#2}%
  \newcommand*\fsize{\dimexpr\f@size pt\relax}%
  \newcommand*\lineheight[1]{\fontsize{\fsize}{#1\fsize}\selectfont}%
  \ifx\svgwidth\undefined%
    \setlength{\unitlength}{104.55151055bp}%
    \ifx\svgscale\undefined%
      \relax%
    \else%
      \setlength{\unitlength}{\unitlength * \real{\svgscale}}%
    \fi%
  \else%
    \setlength{\unitlength}{\svgwidth}%
  \fi%
  \global\let\svgwidth\undefined%
  \global\let\svgscale\undefined%
  \makeatother%
  \begin{picture}(1,0.523072)%
    \lineheight{1}%
    \setlength\tabcolsep{0pt}%
    \put(0,0){\includegraphics[width=\unitlength,page=1]{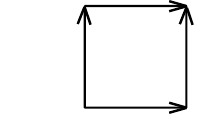}}%
    \put(0.24697374,0.00641143){\makebox(0,0)[lt]{\lineheight{1.25}\smash{\begin{tabular}[t]{l}$0$\end{tabular}}}}%
    \put(0.24212725,0.46073292){\makebox(0,0)[lt]{\lineheight{1.25}\smash{\begin{tabular}[t]{l}$1$\end{tabular}}}}%
    \put(0.93233332,0.45847124){\makebox(0,0)[lt]{\lineheight{1.25}\smash{\begin{tabular}[t]{l}$3$\end{tabular}}}}%
    \put(0.93168731,0.00931989){\makebox(0,0)[lt]{\lineheight{1.25}\smash{\begin{tabular}[t]{l}$2$\end{tabular}}}}%
    \put(-0.00492566,0.22203108){\makebox(0,0)[lt]{\lineheight{1.25}\smash{\begin{tabular}[t]{l}$G=$\end{tabular}}}}%
    \put(0,0){\includegraphics[width=\unitlength,page=2]{KS.pdf}}%
  \end{picture}%
\endgroup%

	%\caption[]{$\Supp(e_{0123})$}
	%\label{Fig:}
\end{figure}
\noindent we have the digraph map $f:I^{\boxdot 2}\longrightarrow G$,
$$f((0,0))=0,\quad f((1,0))=1,\quad f((0,1))=2,\quad f((1,1))=3.$$
Then $f_*(\omega_2)=P$, the image digraph of $f$ is the square given by $\Supp(P)$.

Note that $P$ is not minimal in $G$ because of the edge $0\rightarrow 3$. The paths $P_1=e_{013}$ and $P_2=-e_{023}$ are minimal in $G$, and we have the digraph maps $g_1, g_2:I^{\boxdot 2}\longrightarrow G$ respectively
\begin{align*}
&g_1((0,0))=0,\quad g_1((0,1))=0,\quad g_1((1,0))=1,\quad g_1((1,1))=3;\\
&g_2((0,0))=0,\quad g_2((0,1))=2,\quad g_2((1,0))=0,\quad g_2((1,1))=3.
\end{align*}
\end{example}

\begin{example} In Example 3.11, we give the singular cubical realization of the pair $(P=e_{0135}-e_{0235}+e_{0245},\Supp(P))$. But in the following bigger digraph $G$,
\begin{figure}[H]
	\centering
	%% Creator: Inkscape 1.0.1 (3bc2e813f5, 2020-09-07), www.inkscape.org
%% Accompanies image file 'PP1P2.pdf' (pdf, eps, ps)
\begingroup%
  \makeatletter%
  \providecommand\color[2][]{%
    \errmessage{(Inkscape) Color is used for the text in Inkscape, but the package 'color.sty' is not loaded}%
    \renewcommand\color[2][]{}%
  }%
  \providecommand\transparent[1]{%
    \errmessage{(Inkscape) Transparency is used (non-zero) for the text in Inkscape, but the package 'transparent.sty' is not loaded}%
    \renewcommand\transparent[1]{}%
  }%
  \providecommand\rotatebox[2]{#2}%
  \newcommand*\fsize{\dimexpr\f@size pt\relax}%
  \newcommand*\lineheight[1]{\fontsize{\fsize}{#1\fsize}\selectfont}%
  \ifx\svgwidth\undefined%
    \setlength{\unitlength}{174.98948826bp}%
    \ifx\svgscale\undefined%
      \relax%
    \else%
      \setlength{\unitlength}{\unitlength * \real{\svgscale}}%
    \fi%
  \else%
    \setlength{\unitlength}{\svgwidth}%
  \fi%
  \global\let\svgwidth\undefined%
  \global\let\svgscale\undefined%
  \makeatother%
  \begin{picture}(1,0.33708163)%
    \lineheight{1}%
    \setlength\tabcolsep{0pt}%
    \put(0,0){\includegraphics[width=\unitlength,page=1]{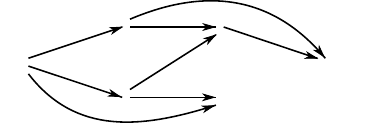}}%
    \put(-0.00238346,0.14078257){\makebox(0,0)[lt]{\lineheight{1.25}\smash{\begin{tabular}[t]{l}$0$\end{tabular}}}}%
    \put(0.27237559,0.29369817){\makebox(0,0)[lt]{\lineheight{1.25}\smash{\begin{tabular}[t]{l}$1$\end{tabular}}}}%
    \put(0.30420972,0.10810195){\makebox(0,0)[lt]{\lineheight{1.25}\smash{\begin{tabular}[t]{l}$2$\end{tabular}}}}%
    \put(0.60107175,0.19937726){\makebox(0,0)[lt]{\lineheight{1.25}\smash{\begin{tabular}[t]{l}$3$\end{tabular}}}}%
    \put(0.59949747,0.00651886){\makebox(0,0)[lt]{\lineheight{1.25}\smash{\begin{tabular}[t]{l}$4$\end{tabular}}}}%
    \put(0.92445107,0.14147016){\makebox(0,0)[lt]{\lineheight{1.25}\smash{\begin{tabular}[t]{l}$5$\end{tabular}}}}%
    \put(0,0){\includegraphics[width=\unitlength,page=2]{PP1P2.pdf}}%
  \end{picture}%
\endgroup%

	%\caption[]{$\Supp(e_{0123})$}
	%\label{Fig:}
\end{figure}
we have
$$P_1=e_{0135}-e_{0235}<P,\quad P_2=e_{0245}<P.$$
And we also have the digraph maps
$$g_1:I^{\boxdot 3}\longrightarrow \Supp(P_1),\quad g_2:I^{\boxdot 3}\longrightarrow \Supp(P_2)$$
which are given by
\begin{align*}
&g_1((0,0,0))=0,\quad g_1((1,0,0))=1,\quad g_1((0,1,0))=2,\quad g_1((0,0,1))=2;\\
&g_1((1,1,0))=3,\quad g_1((1,0,1))=5,\quad g_1((0,1,1))=5,\quad g_1((1,1,1))=5;\\
&g_2((0,0,0))=0,\quad g_2((1,0,0))=2,\quad g_2((0,1,0))=2,\quad g_2((0,0,1))=4;\\
&g_2((1,1,0))=2,\quad g_2((1,0,1))=5,\quad g_2((0,1,1))=4,\quad g_2((1,1,1))=5.\\
\end{align*}
\end{example}

\begin{proof}[Proof of Lemma \ref{Keylemma}] To prove this lemma, we need to construct a singular cubical realization $g:I^{\boxdot n}\longrightarrow \Supp(P')$. The construction is based on the following two steps:
\begin{enumerate}
  \item Construct a map $p:\Supp(P)\rightarrow\Supp(P')$, which is not necessarily a digraph map, but satisfying
  $$p_*(P)=P',~\mod\text{ (non-allowed paths)}.$$
  \item Modify the singular cubical realization $f:I^{\boxdot n}\rightarrow \Supp(P)$ to $f':I^{\boxdot n}\rightarrow \Supp(P)$ such that $f'\circ p$ gives the singular cubical realization of $(P',\Supp(P'))$. More explicitly,
  \begin{itemize}
    \item $f'\circ p: I^{\boxdot n}\rightarrow \Supp(P')$ is a digraph map;
    \item $(f'\circ p)_*(\omega_n)=cP'$ for some $c\in\Z\setminus\{0\}$ and the image digraph of $f'\circ p$ is $\Supp(P')$.
  \end{itemize}
\end{enumerate}
The construction is done by induction on the number of elementary path components of $P$.

If $P$ has only one elementary path components, that is, $P=e_{i_0i_1\cdots i_n}$. Then, the condition $P'\leq P$ implies that $P'=P$. Thus, we can take
$$p=\id_{\Supp(P)},\quad f'=f.$$

\noindent\textbf{Hypothesis:} Assume the pair of maps $(p,f')$ can be constructed for $P$ with $\leq (m-1)$ elementary path components.

Now let us consider $P$ with $m$ ($m\geq2$) elementary path components. First, we can write $P$ as follows:
$$P=\sum_{\alpha_k}P_{S,n-k,\alpha_k}\alpha_{k-1}\alpha_{k-2}\cdots\alpha_1E,\quad\text{for some }k=1,\ldots,n-1,$$
where
\begin{itemize}
  \item $\alpha_k\in d_E^{-1}(k)$, $\alpha_0=E$;
  \item $\sum_{\alpha_k}P_{S,n-k,\alpha_k}$ has $m$ elementary path components and $|d_E^{-1}(k)|\geq2$.
\end{itemize}
For example, the following digraph is the supporting digraph of the path (with 4 elementary path components)
$$P=e_{01467}-e_{01367}+e_{02367}-e_{02567}.$$
\begin{figure}[H]
	\centering
	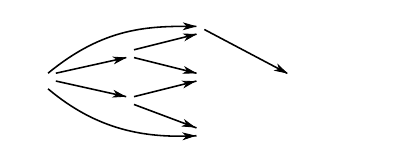
	%\caption[]{$\Supp(e_{0123})$}
	%\label{Fig:Q1}
\end{figure}
We can also write $P$ as
$$P=[(e_{014})+(e_{023}-e_{013})+(e_{025})]67, \quad \text{with } d_E^{-1}(2)=\{3,4,5\}.$$

Then, we apply the technique in the proof of Theorem 4.1 in \cite{TY} to $P_{S,n-k,\alpha_k}\alpha_{k-1}$ (abbreviated as $P_{\alpha_{k-1}}$). That is, we add some edges
$$\alpha_{k+1}\rightarrow \alpha_{k-1},\quad \text{for some } \alpha_{k+1}\in d_E^{-1}(k+1)$$
to obtain a larger digraph $\Supp(P)_1$ and a good splitting of $P_{\alpha_{k-1}}$:
$$P_{\alpha_{k-1}}=P_{\alpha_{k-1},1}+P_{\alpha_{k-1},2},$$
where
\begin{itemize}
  \item $P_{\alpha_{k-1},1}, P_{\alpha_{k-1},2}\in\mathcal{P}_{n-k+1}(\Supp(P)_1)$, and $P_{\alpha_{k-1},1}, P_{\alpha_{k-1},2}<P_{\alpha_{k-1}}$ in $\Supp(P)_1$;
  \item $P_{\alpha_{k-1},1}$ is of the form
$$P_{\alpha_{k-1},1}=\sum_{\alpha_{k+1}}P_{S,n-k-1,\alpha_{k+1}}\alpha_k\alpha_{k-1},\quad \text{i.e. } d_E^{-1}(k)\cap V(P_{\alpha_{k-1},1})=1.$$
\end{itemize}

Let us set
\begin{align*}
&P_1:=P_{\alpha_{k-1},1}\alpha_{k-2}\cdots\alpha_1E=\sum_{\alpha_{k+1}}P_{S,n-k-1,\alpha_{k+1}}\alpha_k\alpha_{k-1}\alpha_{k-2}\cdots\alpha_1E;\\
&P_1':=P-P_1.
\end{align*}
By the above construction, we have the following two observations:
\begin{description}
  \item[(1)] $P_1,P_1'\in\mathcal{P}_n(\Supp(P)_1)$ and $P_1,P_1'<P$ in $\Supp(P)_1$;
  \item[(2)] $|V(P_1)\setminus V(P_1')|\geq 1$.
\end{description}
%\begin{enumerate}
%  \item $P_1,P_2\in\mathcal{P}_n(\Supp(P)_1)$ and $P_1,P_2<P$ in $\Supp(P)_1$;
%  \item $|V(P_1)\setminus V(P_2)|\geq 1$.
%\end{enumerate}

If $|V(P_1)\setminus V(P_1')|>1$, we continue performing the operation to $P_1$ in $\Supp(P)$, and obtain the new digraph $\Supp(P)_2 (\supset\Supp(P)_1\supset\Supp(P))$ and $P_2$ such that
\begin{description}
  \item[(3)] $P_2\in\mathcal{P}_n(\Supp(P)_2)$ is of the form
  $$P_2=\sum_{\alpha_{k+2}}P_{S,n-k-2,\alpha_{k+2}}\alpha_{k+1}\alpha_{k}\alpha_{k-1}\cdots\alpha_1E<P_1;$$
  \item[(4)] $1\leq|V(P_2)\setminus(V(P_1-P_2)\cup V(P_1'))|=|V(P_2)\setminus V(P-P_2)|<|V(P_1)\setminus V(P_1')|$.
\end{description}

Note that, if $P\in\mathcal{P}_n(G)$ has more than 2 elementary path components, then for each elementary path component $e_{i_0i_1\ldots i_n}$ in $P$, there exists at least one path component $e_{j_0j_1\ldots j_n}$ in $P$ such that
$$|V(e_{i_0i_1\ldots i_n})\setminus V(e_{j_0j_1\ldots j_n})|=1.$$
Thus, if $|V(P_2)\setminus(V(P-P_2))|>1$, after several operations, one can obtain a much larger digraph $\Supp(P)_l$ and $P_l$ such that
\begin{description}
  \item[(5)] $P_l\in\mathcal{P}_n(\Supp(P)_l)$ is of the form
  $$P_l=P_{S,n-k-l,\alpha_{k+l}}\alpha_{k+l-1}\alpha_{k+1}\alpha_{k}\cdots\alpha_1E<\cdots<P_2<P_1,\quad l>2$$
  \item[(6)] $|V(P_l)\setminus V(P-P_l)|=1$, and what is more important is that
  $$V(P_l)\setminus V(P-P_l)=\{\alpha_{k+l-1}\}.$$
\end{description}

Without loss of generality, we assume that $|V(P_1)\setminus V(P_1')|=1$ and by our notation,
$$V(P_1)\setminus V(P_1')=\alpha_k.$$
(Or one can work with $P_l$ and remove the redundant new edges to make the path $P-P_l$ minimal.)

Now we need to deal with the following two kinds of local digraphs around $\alpha_k$ in $\widehat{\Supp(P)}$:
\begin{figure}[H]
	\centering
	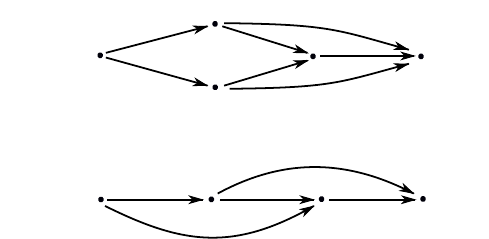
	%\caption[]{$\Supp(e_{0123})$}
	%\label{Fig:Q1}
\end{figure}

where, by our construction in Type I subgraph, one of the edges $\alpha_{k+1}^i\rightarrow \alpha_{k-1}$ ($i=1,2$) is added to obtain $P_1$. And in Type II subgraph, the edge $\alpha_{k+1}^3\rightarrow\alpha_{k-1}$ is new.

Now we define the following map
\begin{equation}\label{p1}
\begin{aligned}
p_1:~&\Supp(P)\rightarrow \Supp(P_1')\\
&p_1(\alpha_k)=\alpha_{k-1},\\
&p_1(v)=v,~\text{ for }v\in V(P)\setminus\{\alpha_k\}.
\end{aligned}
\end{equation}

If $\widehat{\Supp(P)}$ only contains Type I local digraph around $\alpha_k$, we have,
\begin{itemize}
  \item $p_1$ is a digraph map;
  \item $(p_1)_*(P)=P_1'$ and its image digraph is $\Supp(P_1')$.
\end{itemize}
Then the composition $f\circ p_1$ gives us the singular cubical realization of $P_1'(<P)$.

If $\widehat{\Supp(P)}$ only contains Type II local digraph around $\alpha_k$, $p_1$ is not a digraph map, since
$$\alpha_{k+2}^{2}\rightarrow\alpha_k,\quad \alpha_{k+2}^3\nrightarrow \alpha_{k-1}.$$
But we still have
$$(p_1)_*(P)=P_1',~\mod\text{ (non-allowed paths)}.$$

To obtain the singular cubical realization of $P_1'$, let us modify $f:I^{\boxdot n}\rightarrow \Supp(P)$ as follows. Note that the triangle with vertices $$\{\alpha_{k+2}^2,~\alpha_{k+1}^3,~ \alpha_k\}$$
must be realized by a digraph map from the square. Without loss of generality, for convenience, assume it is realized by $f\big|_{(\epsilon_1,\ldots,\epsilon_{n-2})\boxdot I^{\boxdot 2}}$. That is, let us write $\epsilon'=(\epsilon_1,\ldots,\epsilon_{n-2})$, the restriction map is given by exactly one of the following four cases.
\begin{align*}
&f((\epsilon',0,0))=\alpha_{k+2}^2,\quad f((\epsilon',1,0))=\alpha_{k+2}^2,\\
&f((\epsilon',0,1))=\alpha_{k+1}^3,\quad f((\epsilon',1,1))=\alpha_k;
\end{align*}
or
\begin{align*}
&f((\epsilon',0,0))=\alpha_{k+2}^2,\quad f((\epsilon',1,0))=\alpha_k,\\
&f((\epsilon',0,1))=\alpha_{k+1}^3,\quad f((\epsilon',1,1))=\alpha_k;
\end{align*}
or the other two maps induced by the automorphism of the square.

For all the four cases, let us modify $f$ to be
\begin{align*}
f':&~I^{\boxdot n}\rightarrow \Supp(P)\\
&f'((\epsilon',1,0))=f'((\epsilon',0,1))=\alpha_{k+1}^3;\\
&f'(v)=f(v),\quad \text{for }v\neq (\epsilon',1,0), (\epsilon',0,1).
\end{align*}
If there are other edges $a\rightarrow b$ in $I^{\boxdot n}$ mapping to $\alpha_{k+2}^2\rightarrow \alpha_k$ via $f$, similarly, modify it locally to $f'$:
$$f'(a)=f(a)=\alpha_{k+2}^2,\quad f'(b)=\alpha_{k+1}^3.$$
Then, one can easily check that $f'\circ p_1$ gives the singular cubical realization of $P_1'$. That is,
\begin{itemize}
  \item $f'\circ p_1:I^{\boxdot n}\rightarrow \Supp(P_1')$ is a digraph map;
  \item $(f'\circ p_1)_*(\omega_n)=cP_1'$, where $c$ is the same as above.
\end{itemize}

Now if $\widehat{\Supp(P)}$ contain both Type I and Type II local digraphs around $\alpha_k$, there must be $\alpha_{k+2}^1\neq\alpha_{k+2}^2$, then $f'$ will have no effect on Type I case. Thus, the modification $f'$ still works.

In summary, we construct a singular cubical realization of $P_1'$ which has elementary path components less than $m$. By Hypothesis, for any $\tilde{P}_1'\in\mathcal{P}_n(G)$, $\tilde{P}_1'\leq P_1'$, we have $\tilde{P}_1'\in\mathcal{P}_{\adm,n}(G)$. Now any $P'\in\mathcal{P}_n(G)$ with $P'<P$, it could be first realized as some $P_1'$ or a smaller path, and then we obtain a singular cubical realization for $P'$. We are done.
\end{proof}

To complete the proof of the coincidence of $\partial^{\cell}$ and $\partial$, we require the ensuing lemma on unique decomposition.

\begin{lemma}\label{unique} Let $P\in\mathcal{P}(G)$ and $R$ be a minimal path (of arbitrary length) in $\Supp(P)$, if $R$ is not minimal in $G$, it admits the following unique decomposition (up to the summation order) in $G$:
\begin{equation}\label{decompunique2}
R=R_1+R_2+\cdots+R_k,\quad R_i\in\mathcal{P}(G),~R_i<R,~i=1,2,\ldots,k.
\end{equation}
In other words, there do not exist admissible relations among $\{R_i\in\mathcal{P}_n(G)|R_i<R\}$.
\end{lemma}

\begin{proof} The proof is done by induction on the length of $R$. For the case of length $0$ and length-$1$, it's obvious; and the length-2 result follows from the structure theorem \ref{structurethm} (3). Assume the result holds for the case of length $(n-1)$. Now, let $R$ be a length $n$ minimal path in $\Supp(P)$. First we have
$$R=\sum_{\alpha\in E_1}R_{\alpha}E,$$
where $R_{\alpha}\in\mathcal{P}_{n-1}(\Supp(P))$ with the ending vertex $\alpha$.
By hypothesis, we can decompose $R_{\alpha}$ in $G$ and write $R$ as
$$R=\sum_{\alpha\in E_1}\sum_{i\in I}R_{\alpha,i}E.$$
Now if $R$ admits two different decompositions up to the summation order
\begin{equation}\label{twodecomp}
\begin{aligned}
&R=R_1+R_2+\cdots+R_k,\qquad R_i\in\mathcal{P}_n(G),~R_i<R,~i=1,2,\ldots,k;\\
&R=R_1'+R_2'+\cdots+R_l',\qquad  R_j'\in\mathcal{P}_n(G),~R_j'<R,~j=1,2,\ldots,l.
\end{aligned}
\end{equation}
Without loss of generality, assume that $R_1\neq R_j'$, for any $j=1,2,\ldots,l$. And we can expand $R_1$ as
$$R_1=\sum_{\alpha\in A\subset E_1}\sum_{i\in I'\subset I}R_{\alpha,i}E.$$
Since the decomposition for $R_{\alpha}$ in $G$ is unique, by \eqref{twodecomp}, for any $\alpha\in A$ and $i\in I_1$, $R_{\alpha,i}E$ is also a component in some $R_j'$.
We can also expand $R_j'$ as
$$R_1=\sum_{\alpha\in A'\subset E_1}\sum_{i\in I''\subset I}R_{\alpha,i}E$$

Assume that there exists $\beta$ in $R_{\alpha,i}$, such that
$$\beta\rightarrow \alpha, \quad \beta\nrightarrow E.$$
(Otherwise $R_{\alpha,i}E\in\mathcal{P}_{n}$, the minimal conditions for $R_1$ and $R_j'$ imply that $R_1=R_j'=R_{\alpha,i}$, contradiction.)
Then the condition $R_1, R_j\in\mathcal{P}_n(G)$ implies that there exist path components $R_{\alpha',i'}$ in $R_1$ and $R_{\alpha'',i''}$ in $R_j'$ such that
$$\beta\rightarrow \alpha',\quad \beta\rightarrow \alpha''.$$
If $\alpha'\neq\alpha''$, then $R_1$, $R_j<R$ and $R$ is a minimal path in $\Supp(P)$, which imply that $\Supp(P)$ contains the following subgraph
\begin{figure}[H]
\vskip -0.2cm
	\centering
	%% Creator: Inkscape 1.0.1 (3bc2e813f5, 2020-09-07), www.inkscape.org
%% Accompanies image file 'subdigraph.pdf' (pdf, eps, ps)
\begingroup%
  \makeatletter%
  \providecommand\color[2][]{%
    \errmessage{(Inkscape) Color is used for the text in Inkscape, but the package 'color.sty' is not loaded}%
    \renewcommand\color[2][]{}%
  }%
  \providecommand\transparent[1]{%
    \errmessage{(Inkscape) Transparency is used (non-zero) for the text in Inkscape, but the package 'transparent.sty' is not loaded}%
    \renewcommand\transparent[1]{}%
  }%
  \providecommand\rotatebox[2]{#2}%
  \newcommand*\fsize{\dimexpr\f@size pt\relax}%
  \newcommand*\lineheight[1]{\fontsize{\fsize}{#1\fsize}\selectfont}%
  \ifx\svgwidth\undefined%
    \setlength{\unitlength}{146.83839573bp}%
    \ifx\svgscale\undefined%
      \relax%
    \else%
      \setlength{\unitlength}{\unitlength * \real{\svgscale}}%
    \fi%
  \else%
    \setlength{\unitlength}{\svgwidth}%
  \fi%
  \global\let\svgwidth\undefined%
  \global\let\svgscale\undefined%
  \makeatother%
  \begin{picture}(1,0.4940469)%
    \lineheight{1}%
    \setlength\tabcolsep{0pt}%
    \put(0,0){\includegraphics[width=\unitlength,page=1]{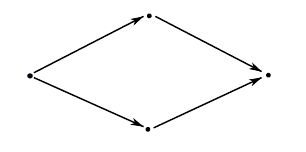}}%
    \put(-0.00387292,0.22819892){\makebox(0,0)[lt]{\lineheight{1.25}\smash{\begin{tabular}[t]{l}$\beta$\end{tabular}}}}%
    \put(0.43911474,0.46217468){\makebox(0,0)[lt]{\lineheight{1.25}\smash{\begin{tabular}[t]{l}$\alpha$\end{tabular}}}}%
    \put(0.43610737,0.00833215){\makebox(0,0)[lt]{\lineheight{1.25}\smash{\begin{tabular}[t]{l}$\alpha''$\end{tabular}}}}%
    \put(0.90612346,0.22498683){\makebox(0,0)[lt]{\lineheight{1.25}\smash{\begin{tabular}[t]{l}$E$\end{tabular}}}}%
    \put(0,0){\includegraphics[width=\unitlength,page=2]{subdigraph.pdf}}%
    \put(0.44293967,0.26831076){\makebox(0,0)[lt]{\lineheight{1.25}\smash{\begin{tabular}[t]{l}$\alpha'$\end{tabular}}}}%
  \end{picture}%
\endgroup%

	%\caption[]{$\Supp(e_{0123})$}
	%\label{Fig:Q1}
\end{figure}
\vskip 0.1cm
It contradicts with the structure of $\Supp(P)$. Then $\alpha'=\alpha''$. Furthermore, the unique decompositions of all $R_{\alpha}$'s imply that $i''=i'$. Thus,
$$R_{\alpha',i'}=R_{\alpha'',i''}.$$
Continuing the discussion, one can obtain two $\partial$-invariant completions $\widetilde{R_1}$, $\widetilde{R_j'}$ of $R_{\alpha,i}+R_{\alpha',i'}$ in $R_1$ and $R_j'$. Applying the same argument, we obtain
$$R_1=\widetilde{R_1}=\widetilde{R_j'}=R_j',$$
contradiction. Thus, the required decomposition is unique.
\end{proof}

\begin{corollary}\label{coincidence} For any $n\geq1$, $P\in\mathcal{P}_{\adm,n}(G)$, $\partial P$ admits a representation which coincides with $\partial^{\cell}P$.
\end{corollary}

\begin{proof} First, the decomposition of $\partial P$ in $\Supp(P)$ is unique. We write $\partial P$ as
$$\partial P=Q_1+Q_2+\cdots+Q_k, \quad Q_i\in\mathcal{P}_{\adm,n-1}(\Supp(P)).$$
Now we consider all the paths $Q_i$ in $G$. If $Q_i$ is not minimal in $G$, by Lemma \ref{unique}, we have the unique decomposition for $Q_i$ as follows:
\begin{equation}\label{decompunique}
Q_i=\sum_{r=1}^{k_i}Q_{ir},\quad Q_{ir}\in\mathcal{P}_{n-1}(G),~Q_{ir}<Q_i.
\end{equation}
By Lemma \ref{Keylemma}, we have $Q_{ir}\in\mathcal{P}_{\adm,n-1}(G)$.

According to the definition of $\partial^{\cell}P$, the uniqueness of the decomposition implies that the admissible paths $\{Q_{ir}\}_{i,r}$ are all the summands in $\partial^{\cell}P$. It follows that
$$\partial^{\cell}P=\sum_{i,r}Q_{ir}=\partial P.\qedhere$$
\end{proof}

\begin{proof}[Proof of Proposition \ref{welldef}]
Assume that $P_1,\ldots,P_k\in\mathcal{P}_{\adm,n}(G)$ satisfy the admissible (linear) relations, which we write it as
\begin{equation}\label{adm}
P_k=\sum_{i=1}^{k-1}c_iP_i,\quad c_i\in\mathbb{R}.
\end{equation}
First, by Corollary \ref{coincidence}, for each $P_i\in\mathcal{P}_{\adm,n}(G)$ in the relation, we have
$$\partial^{\cell}P_i=\sum_{r}Q_{ir}=\partial P_i.$$
Then on the one hand, by $\R$-linear extension,
$$\partial^{\cell}\left(\sum_{i=1}^{k-1}c_iP_i\right)=\sum_{i=1}^{k-1}c_i\partial^{\cell}P_i=\sum_{i=1}^{k-1}c_i\sum_{r}Q_{ir}=\partial\left(\sum_{i=1}^{k-1}c_iP_i\right);$$
on the other hand,
$$\partial^{\cell}P_k=\sum_rQ_{kr}=\partial P_k.$$
It follows from the condition $Q_{ir}\in\mathcal{P}_{\adm,n-1}(G)$ and the relation \eqref{adm} that
$$\partial P_k=\partial\left(\sum_{i=1}^{k-1}c_iP_i\right)$$
imply the admissible relations among $\{Q_{ir}\}_{i=1,\ldots,k;r}$
$$\partial^{\cell} P_k=\sum_rQ_{kr}=\sum_{i=1}^{k-1}c_i\sum_{r}Q_{ir}=\partial^{\cell}\left(\sum_{i=1}^{k-1}c_iP_i\right).$$
Thus, we obtain a well-defined operator, still denoted by $\partial^{\cell}$, on the cellular chain group
$$\partial^{\cell}:C_n(G;\R)\longrightarrow C_{n-1}(G;\R).$$
Its relation with $\partial$ of path complex implies that $\partial^{\cell}$ is a boundary operator.
\end{proof}
%\noindent\textbf{Step 2.} Note that each face component of $P$ in $\Supp(P)$ is a linear combination of allowed elementary paths with coefficients $1$ or $-1$. But it may not be minimal in $G$.
%\begin{itemize}
%  \item If all of them are minimal in $G$, then $P$ is preserved by $\partial^{\CW}$ in $G$.
%  \item If some face component is not minimal in $G$, we will show that any of its minimal component has a singular cubical realization. More explicitly, we need to prove the following claim.
%\end{itemize}
%\noindent\textbf{Claim:} If $P\in\Omega_n(G)$ is a linear combination of allowed elementary paths with coefficients $1$ or $-1$, and there exists a digraph map
%$$f:I^{\boxdot n}\longrightarrow G,$$
%such that $f_*(\omega_n)=P$, and the image digraph of $f$ is $\Supp(P)$.
%If $P'\in\Omega_n(G)$ and $P'<P$, then there exists a digraph map
%$$g:I^{\boxdot n}\longrightarrow G,$$
%such that
%\begin{enumerate}
%  \item $g\big|_{f^{-1}(V(P'))}=f\big|_{f^{-1}(V(P'))}$;
%  \item $g_*(\omega_n)=P'$, and the image digraph of $g$ is $\Supp(P')$.
%\end{enumerate}
%
%
%Now assume the claim holds for any length $n-1$ element satisfying the above conditions. Now for $P$ in the claim,

\begin{definition} Let $G$ be a digraph. For any $n\in\mathbb{N}$, the $n$-th cellular homology of $G$ is defined to be
$$H_n^{\cell}(G;\R)=\frac{\ker\{\partial_n^{\cell}: C_n(G;\R)\longrightarrow C_{n-1}(G;\R)\}}{\im\{\partial_{n+1}^{\cell}: C_{n+1}(G;\R)\longrightarrow C_{n}(G;\R)\}}.$$
\end{definition}

\begin{example}\label{simpleex} We can compute the cellular homologies of the above several digraphs.
\begin{enumerate}
  \item For the $k$-square digraph $S_k$ in Example \ref{k-square}, we have
\begin{equation*}
H_n^{\cell}(S_k;\R)\cong\left\{
                      \begin{array}{ll}
                        \R, & \hbox{$n=0$;} \\
                        0, & \hbox{others.}
                      \end{array}
                    \right.
\end{equation*}
  \item For the exotic cube, we have
\begin{equation*}
H_n^{\cell}(\mbox{exotic cube};\R)\cong\left\{
                      \begin{array}{ll}
                        \R, & \hbox{$n=0,2$;} \\
                        0, & \hbox{others.}
                      \end{array}
                    \right.
\end{equation*}
  \item For the supporting digraph of the path $P$,
\begin{align*}
P=~&e_{S159E}-e_{S169E}+e_{S269E}\\
   &+e_{S16(10)E}-e_{S26(10)E}+e_{S27(10)E}-e_{S37(10)E}\\
   &-e_{S27(11)E}+e_{S37(11)E}-e_{S38(11)E}+e_{S48(11)E},
\end{align*}
i.e. the digraph in Example \ref{length4}, according to the analysis of its admissible paths in $\Supp(P)$, i.e. Example \ref{length4rev}, we have
\begin{equation*}
H_n^{\cell}(\Supp(P);\R)\cong\left\{
                      \begin{array}{ll}
                        \R, & \hbox{$n=0,2$;} \\
                        0, & \hbox{others.}
                      \end{array}
                    \right.
\end{equation*}
\end{enumerate}
\end{example}

\vskip 0.2cm

\subsection{Properties of cellular homology of digraphs}\label{properties}

In this subsection, we will discuss several properties of our cellular homology of digraphs, which include functoriality and homotopy invariance in the category of acyclic digraphs, acyclic model and K\"{u}nneth formula for Cartesian product of digraphs.

\subsubsection{Functoriality and Homotopy invariance}

Let us recall basic homotopy theory of digraphs developed in \cite{GLMY}.

Fix $n\geq0$. Denote by $I_n=(V(I_n),E(I_n))$ any digraph, with $V(I_n)=\{0,1,\ldots,n\}$ and $E(I_n)$ containing exactly one of the edges $i\rightarrow i+1$, $i+1\rightarrow i$ for any
$i=0,1,\ldots, n-1$. We call it the line digraph.

\begin{definition}[Homotopy]\label{homotopyequiv} Let $G$, $H$ be two digraphs.
\begin{enumerate}
  \item Two digraph maps $f,g:G\rightarrow H$
are called \textbf{homotopic} if there exists a line digraph $I_n$ for some $n\geq1$ and a digraph map $F:G\boxdot I_n\rightarrow H$, such that
$$F\big|_{G\boxdot\{0\}}=f\quad\text{and}\quad F\big|_{G\boxdot\{n\}}=g.$$
We shall write $f\overset{F}{\simeq} g$ (or simply $f\simeq g$) and call $F$ a $n$-step homotopy between $f$ and $g$.
  \item Two digraphs $G$ and $H$ are called \textbf{homotopy equivalent} if there exist digraph maps
  $$f:G\rightarrow H,\quad g:H\rightarrow G$$
  such that
  $$f\circ g\simeq \id_H,\quad g\circ f\simeq\id_G.$$
A digraph $G$ is called \textbf{contractible} if $G\simeq\{*\}$ where $\{*\}$ is a single vertex digraph.
\end{enumerate}

\end{definition}

\begin{theorem}[Homotopy invariance of path complex\cite{GLMY}] Let $G$, $H$ be two digraphs.
\begin{enumerate}
  \item Let $f\simeq g: G\rightarrow H$ be two homotopic digraph maps. Then these maps
induce the identical homomorphisms of path homologies of $G$ and $H$:
$$f_*=g_*:H_*(G)\rightarrow H_*(H).$$
  \item If the digraphs $G$ and $H$ are homotopy equivalent, then they have isomorphic path homologies.
\end{enumerate}
\end{theorem}

The proof is first done for one-step homotopy $f\overset{F}{\simeq} g$ and then by induction on the step of the homotopy. The key point there is to construct the chain homotopy
$$L_p:\Omega_p(G)\rightarrow\Omega_{p+1}(H),\quad p\geq0$$
such that
$$\partial L_p+L_{p-1}\partial=g_*-f_*.$$
The answer is given by
$$L_p(v)=F_*(v\times e_{01}),$$
where $\times$ is the cross product of paths, which will be recalled in Definition \ref{cross}.

\vskip 0.2cm
Now let us return to cellular homologies of digraphs in the category of acyclic digraphs.

\begin{theorem}[Funtoriality]\label{1} Let $G$, $H$ be two acyclic digraphs. If $f:G\rightarrow H$ is a digraph map, then it induces a chain map
$$f_*:C_*(G;\R)\rightarrow C_*(H;\R).$$
\end{theorem}

\begin{proof} Let $P\in\mathcal{P}_{\adm,n}(G)$, and $\phi: I^{\boxdot n}\rightarrow G$ be its singular cubical realization. Then
$$f_*(P)=c^{-1}f_*\phi_*(\omega_n)=c^{-1}(f\circ\phi)_*(\omega_n),\quad \text{for some }c\in\Z\setminus\{0\}.$$
It means that $f_*(P)$ admits a singular cubical realization. Furthermore, it follows from the proof of Lemma \ref{elementarylem} that, any $Q\in\mathcal{P}_n(H)$, with $Q\leq f_*(P)$, we have $Q\in\mathcal{P}_{\adm,n}(H)$. Thus, $f_*(P)\in C_*(H;\R)$.
\end{proof}

\begin{remark} It is worth to mention that the digraph map is not a cellular map, since the push-forward $f_*(P)$ of a minimal path is not necessarily minimal.
\end{remark}

\begin{lemma}\label{2} If $v\in\mathcal{P}_{\adm,p}(G)$, we have
$$v\times e_{01}\in\mathcal{P}_{\adm,p+1}(G\boxdot I).$$
\end{lemma}
\begin{proof} We will prove the general case in Lemma \ref{PtimesQ}, which states as follows: if $P\in\mathcal{P}_{\adm,p}(X)$, $Q\in\mathcal{P}_{\adm,q}(Y)$, then
$$P\times Q\in\mathcal{P}_{\adm,p+q}(X\boxdot Y).$$
\end{proof}

Following from Lemmas \ref{1}, \ref{2} and the same proof in the path complex, the above morphism $L_p$ gives the chain homotopy between $C_*(G;\R)$ and $C_*(H;\R)$. Thus we obtain the homotopy invariance of cellular homology of digraphs.

\begin{theorem}[Homotopy invariance]\label{homotopyinv} Let $G$, $H$ be two acyclic digraphs.
\begin{enumerate}
  \item Let $f\simeq g: G\rightarrow H$ be two homotopic digraph maps. Then
$$f_*=g_*:H_*^{\cell}(G;\R)\rightarrow H_*^{\cell}(H;\R).$$
  \item If the digraphs $G$ and $H$ are homotopy equivalent, then they have isomorphic cellular homologies.
\end{enumerate}
\end{theorem}

\subsubsection{Acyclic result for $\widetilde{H}^{\cell}_*(\Supp(P);\R)$, $P\in\mathcal{P}_{\adm}(G)$}
One can define the reduced cellular homologies $\tilde{H}_*^{\cell}(G;\R)$ of $G$ in the usual way by modifying $H_0^{\cell}(G;\R)$ to be
$$\widetilde{H}_0^{\cell}(G;\R)=\ker\{\partial^{\cell}_{0}: C^{\cell}_0(G;\R)\rightarrow \R\}/\im\{\partial_1^{\cell}: C_1^{\cell}(G;\R)\rightarrow C_0^{\cell}(G;\R)\}.$$
where
$$\partial^{\cell}_{0}: C^{\cell}_0(G;\R)\rightarrow \R,\quad \partial^{\cell}_0\left(\sum_{a\in V(G)}c_a a\right)=\sum_{a\in V(G)}c_a.$$

It follows from Theorem \ref{acyclic} and the relation between path complex and cellular chain complex associated to the digraphs, we obtain
\begin{proposition} For $P\in\mathcal{P}_{\adm}(G)$, we have
$$\widetilde{H}_*^{\cell}(\Supp(P);\R)=0.$$
\end{proposition}

We list all the admissible minimal $3$-paths in Appendix \ref{AppendixA}. In particular, all the corresponding supporting digraphs are contractible as we expect. Thus, we propose the following conjecture.
\begin{conjecture}\label{contractible} For any $n\geq0$ and $P\in\mathcal{P}_{\adm,n}(G)$, the digraph $\Supp(P)$ is contractible.
\end{conjecture}

\subsubsection{K\"{u}nneth formula}

In this subsection, we will study the cellular chain structure of the Cartesian product of digraphs.

\begin{definition}[Cross product]\label{cross}
Let $X$, $Y$ be two digraphs. For any s-regular allowed elementary $p$-path $e_x$ in $X$ and $q$-path $e_y$ in $Y$, the cross product $e_x\times e_y$ is defined as a $(p+q)$-path in $X\boxdot Y$ by
$$e_x\times e_y=\sum_{z\in\Sigma_{x,y}}(-1)^{L(z)}e_z,$$
where $\Sigma_{x,y}$  is the set of all stair-like paths $z$ on $X\boxdot Y$ whose projections on $X$ and $Y$ are respectively $x$ and $y$, $L(z)$ is the number of cells in $\mathbb{N}_+^2$ below the staircase $S(z)$, see the following picture\footnote{We refer to the picture in \cite{GMY}.} for an explanation.
\begin{figure}[H]
  \centering
  \includegraphics[width=14cm]{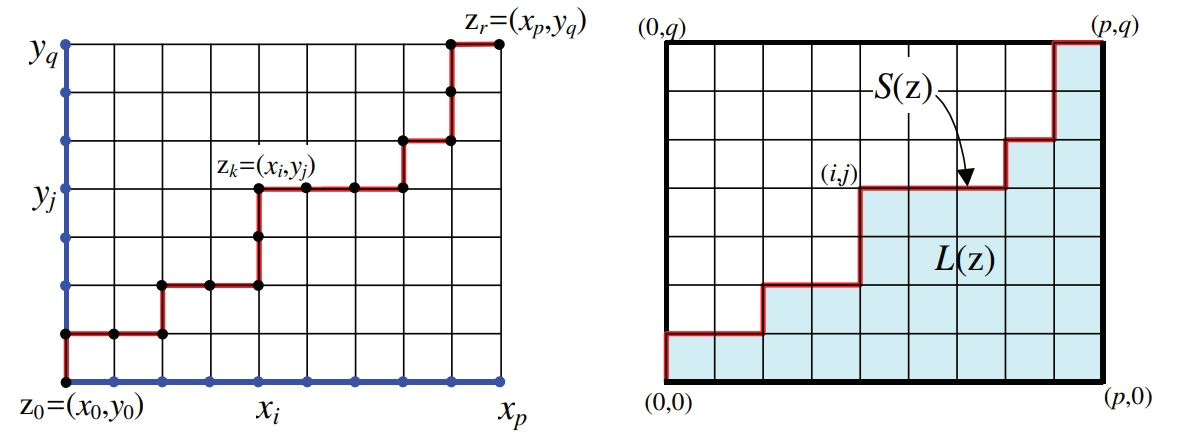}\\
  %\caption{}\label{}
\end{figure}
\end{definition}

\begin{example}\label{crossexample} Let $a,b\in V(X)$, $1,2\in V(Y)$, we have
\begin{enumerate}
  \item $e_a\times e_{12}=e_{(a1)(a2)}$, $e_{ab}\times e_1=e_{(a1)(b1)}$;
  \item $e_{ab}\times e_{12}=e_{(a1)(b1)(b2)}-e_{(a1)(a2)(b2)}$.
\end{enumerate}
\end{example}

The following results which are originally proven for GLMY path complex, also hold under the strongly regular condition.
\begin{lemma}[\cite{GMY}]\label{useful12} Let $X$, $Y$ be two digraphs.
\begin{enumerate}
  \item The cross product extends by linearity to the vector space of paths. For $u\in\Omega_p(X)$, and $v\in\Omega_q(Y)$, we have $u\times v\in\Omega_{p+q}(X\boxdot Y)$. The boundary operator $\partial$ is compatible with the cross product in the following way:
$$\partial(u\times v)=(\partial u)\times v+(-1)^pu\times(\partial v).$$
In particular, if $X=I^{\boxdot p}$ and $Y=I^{\boxdot q}$, then
$$\omega_p\times\omega_q=\omega_{p+q}\in\Omega_{p+q}(I^{\boxdot(p+q)};\Z).$$
  \item Any $\omega\in\Omega_*(X\boxdot Y)$ admits a representation
$$\omega=\sum_{e_y\in E_*(Y)}u^y\times e_y,\quad u^y\in\Omega_*(X),$$
where $E_*(Y)=\cup_{q\geq0}E_q(Y)$, with $E_q(Y)$ being the set of all s-regular allowed elementary $q$-paths in $Y$.
\end{enumerate}
\end{lemma}

\begin{remark}\label{useful12rem} In \cite{GMY}, the authors work with the coefficient being the field $\mathbb{K}$. But the above lemma also works for the $\Z$-coefficient. In particular,
any $\omega\in\Omega_*(X\boxdot Y;\Z)$ admits a representation
$$\omega=\sum_{e_y\in E_*(Y)}u^y\times e_y,\quad u^y\in\Omega_*(X;\Z),$$
%since all the non-zero $u^y\times e_y$'s are linearly independent for different $e_y$.
since for different $e_{y_1}$, $e_{y_2}\in E_*(Y)$, if we represent $u^{y_i}\times e_{y_i}$ in terms of elementary paths uniquely, then the corresponding two sets of elementary paths have empty intersections:
$$E(u^{y_1}\times e_{y_1})\cap E(u^{y_2}\times e_{y_2})=\emptyset.$$
\end{remark}

Here, we expect a more restricted result. That is
\begin{lemma}\label{PtimesQ} For any $P\in\mathcal{P}_{\adm,p}(X)$ and $Q\in\mathcal{P}_{\adm,q}(Y)$, we have
$$P\times Q\in\mathcal{P}_{\adm,p+q}(X\boxdot Y).$$
\end{lemma}
\begin{proof} By definition, we need to prove
\begin{enumerate}
  \item $P\times Q$ is minimal in $X\boxdot Y$;
  \item $P\times Q$ has a singular cubical realization.
\end{enumerate}

(1) Let $R=P\times Q\in\Omega_k(X\boxdot Y;\Z)$, $k=p+q$. Since $Q\in\mathcal{P}_{\adm,q}(Y)\subset\mathcal{P}_q(Y)$, then we can write $Q=\sum_{e_y\in E_q(Y)}\epsilon_ye_y$, where $\epsilon_y=0$, $1$ or $-1$. Thus,
$$R=P\times Q=\sum_{e_y\in E_q(Y)}P\times \epsilon_ye_y=\sum_{e_y\in E_q(Y)}\epsilon_yP\times e_y.$$

Assume that $R$ is not minimal, that is, there exists a path $\hat{R}\in\Omega_k(X\boxdot Y;\Z)$ such that $\hat{R}<R$.

By Lemma \ref{useful12}(2) and Remark \ref{useful12rem}, $\hat{R}$ admits a representation
$$\hat{R}=\sum_{e_y\in E_q(Y)}\hat{R}^y\times e_y,\quad \hat{R}^y\in\Omega_p(X;\Z).$$

For any paths $u$, $v$ (it is possible that $v=u$) and different elementary paths $e_{y_1}$, $e_{y_2}$,
$$E(u\times e_{y_1})\cap E(v\times e_{y_2})=\emptyset,$$
then $\hat{R}<R$ implies
$$\hat{R}^y\leq \epsilon_y P.$$
\begin{itemize}
  \item If $\epsilon_y=0$, then $\hat{R}^y=0$.
  \item If $\epsilon_y=\pm1$, since $P$ is minimal, then $\hat{R}^y=0$ or $\hat{R}^y=\epsilon_yP$.
\end{itemize}
Thus we write
\begin{equation*}
\hat{R}^y=\epsilon_y'P:=\left\{
                 \begin{array}{ll}
                 \epsilon_yP, & \hbox{if $\hat{R}^y=\epsilon_yP$;} \\
                        0, & \hbox{if $\hat{R}^y<\epsilon_yP$.}
                      \end{array}
                    \right.
\end{equation*}

Now, we can write $\hat{R}$ as follows:
$$\hat{R}=\sum_{e_y\in E_q(Y)}\hat{R}^y\times e_y=\sum_{e_y\in E_q(Y)}\epsilon_y'P\times e_y=\sum_{e_y\in E_q(Y)}P\times\epsilon_y'e_y.$$
Since $\hat{R}\in\Omega_k(X\boxdot Y;\Z)$, by Lemma \ref{useful12} (1), we have
$$\partial\hat{R}=(\partial P)\times \sum_{e_y\in E_q(Y)}\epsilon_y'e_y+(-1)^pP\times \partial\bigg(\sum_{e_y\in E_q(Y)}\epsilon_y'e_y\bigg)\in \mathcal{A}_{k-1}(X\boxdot Y,\Z).$$
Combining with the conditions $P\in\Omega_p(X;\Z)$ and $\sum_{e_y\in E_q(Y)}\epsilon_y'e_y\in \mathcal{A}_q(Y;\Z)$, we obtain
$$\partial\bigg(\sum_{e_y\in E_q(Y)}\epsilon_y' e_y\bigg)\in \mathcal{A}_{q-1}(Y;\Z).$$
Thus, $\sum_{e_y\in E_q(Y)}\epsilon_y'e_y\in \Omega_q(Y;\Z)$. It further implies that
$$\sum_{e_y\in E_q(Y)}\epsilon_y'e_y<Q,$$
contradiction. Thus $P\times Q$ is minimal.

\vskip 0.2cm

(2) Since $P$ and $Q$ are admissible, there exist singular cubical realizations
$$f:I^{\boxdot p}\longrightarrow \Supp(P)\subset X,\qquad g:I^{\boxdot q}\longrightarrow \Supp(Q)\subset Y.$$
Now we define a map
$$h:I^{\boxdot(p+q)}\longrightarrow X\boxdot Y$$
as follows:
$$h(a,b)=(f(a),g(b)), \quad a\in V(I^{\boxdot p}),~b\in V(I^{\boxdot q}).$$
It is easy to check that $h$ is a digraph map by definition of the Cartesian product. It remains to prove that $h$ gives a singular cubical realization of $P\times Q$.
\begin{claim} For any allowed elementary $p$-path $e_p=e_{i_0i_1\cdots i_p}$ in $I^{\boxdot p}$, $q$-path $e_q=e_{j_0j_1\cdots j_q}$ in $I^{\boxdot q}$, we have
$$h_*(e_p\times e_q)=f_*(e_p)\times g_*(e_q).$$
\end{claim}
By definition, we have
\begin{align*}
h_*(e_p\times e_q)%=~&h_*\left(\sum_{r\in\Sigma_{p,q}}(-1)^{L(r)}e_r\right)\\
                  =~&h_*\left(\sum_{(i_0,j_0)\cdots(i_p,j_q)\in\Sigma_{p,q}}(-1)^{L((i_0,j_0)\cdots(i_p,j_q))}e_{(i_0,j_0)\cdots(i_p,j_q)}\right)\\
                  =~&\sum_{(i_0,j_0)\cdots(i_p,j_q)\in\Sigma_{p,q}}(-1)^{L((i_0,j_0)\cdots(i_p,j_q))}e_{(f(i_0),g(j_0))\cdots(f(i_p),g(j_q))}.\\
                  %=~&\sum_{(f(i_0),g(j_0))\cdots(f(i_p),g(j_q))\in\Sigma_{f_*(p),g_*(q)}}(-1)^{L((f(i_0),g(j_0))\cdots(f(i_p),g(j_q)))}e_{(f(i_0),g(j_0))\cdots(f(i_p),g(j_q))}\\
                  %=~&e_{f(i_0)\cdots f(i_p)}\times e_{g(j_0)\cdots g(j_q)}\\
                  %=~&f_*(e_p)\times g_*(e_q).
\end{align*}
If there exist $i_k$ and $i_{k+1}$ such that $f(i_k)=f(i_{k+1})$, then
\begin{itemize}
  \item $f_*(e_p)=0$;
  \item $(f(i_k),g(j_l))=(f(i_{k+1}),g(j_l))$, for any $l=0,1\ldots,q$.
\end{itemize}
Thus,
$$h_*(e_p\times e_q)=0=f_*(e_p)\times g_*(e_q).$$
Similar argument holds for the case $g(j_l)=g(j_{l+1})$ for some $l$.

If $f_*(e_p)$ and $g_*(e_q)$ do not vanish, then there is a bijection between $\Sigma_{p,q}$ and $\Sigma_{f_*(p), f_*(q)}$ given by
$$(i_0,j_0)\cdots(i_p,j_q)\in\Sigma_{p,q}~\leftrightarrow~ (f(i_0),g(j_0))\cdots(f(i_p),g(j_q))\in\Sigma_{f_*(p),g_*(q)}$$
and moreover
$$L((i_0,j_0)\cdots(i_p,j_q))=L((f(i_0),g(j_0))\cdots(f(i_p),g(j_q))),$$
since the staircases on both sides have the same shape.
Thus, we can continue the computation
\begin{align*}
h_*(e_p\times e_q)=~&\sum_{(f(i_0),g(j_0))\cdots(f(i_p),g(j_q))\in\Sigma_{f_*(p),g_*(q)}}(-1)^{L((i_0,j_0)\cdots(i_p,j_q))}e_{(f(i_0),g(j_0))\cdots(f(i_p),g(j_q))}\\
                  =~&\sum_{(f(i_0),g(j_0))\cdots(f(i_p),g(j_q))\in\Sigma_{f_*(p),g_*(q)}}(-1)^{L((f(i_0),g(j_0))\cdots(f(i_p),g(j_q)))}e_{(f(i_0),g(j_0))\cdots(f(i_p),g(j_q))}\\
                  =~&e_{f(i_0)\cdots f(i_p)}\times e_{g(j_0)\cdots g(j_q)}\\
                  =~&f_*(e_p)\times g_*(e_q).
\end{align*}
Since $h_*$, $f_*$ and $g_*$ are linear maps, we get
$$h_*(\omega_{p+q})=h_*(\omega_p\times\omega_q)=f_*(\omega_p)\times g_*(\omega_q)=c_PP\times c_QQ,\quad c_P,c_Q\in\Z\setminus\{0\}.$$
Note that the image digraph of $h$ is given by $\Supp(P)\boxdot\Supp(Q)$ and
$$P\times Q\in\Omega_{p+q}(\Supp(P)\boxdot\Supp(Q)),$$
thus
$$\Supp(P\times Q)\subset\Supp(P)\boxdot\Supp(Q).$$
Conversely, any edge $e_{ij}\in\mathcal{A}_1(\Supp(P)\boxdot\Supp(Q))$, it is of the form
$$e_{(i_0,j_0)(i_0,j_1)}\quad\text{or}\quad e_{(i_0,j_0)(i_1,j_0)}.$$
Meanwhile,
\begin{align*}
&e_{(i_0,j_0)(i_0,j_1)}\in\mathcal{A}_1(\{i_0\}\boxdot\Supp(Q))=\mathcal{A}_1(\Supp(\{i_0\}\times Q))\subset\mathcal{A}_1(\Supp(P\times Q));\\
&e_{(i_0,j_0)(i_1,j_0)}\in\mathcal{A}_1(\Supp(P)\boxdot\{j_0\})=\mathcal{A}_1(\Supp(P\times \{j_0\}))\subset\mathcal{A}_1(\Supp(P\times Q)).
\end{align*}
Then
$$\Supp(P)\boxdot\Supp(Q)=\bigcup_{e_{ij}\in\mathcal{A}_1(\Supp(P)\boxdot\Supp(Q))}\{i\rightarrow j\}\subset\Supp(P\times Q).$$
Thus, we have
$$\Supp(P\times Q)=\Supp(P)\boxdot\Supp(Q).$$
\end{proof}

Combining with the above two lemmas, we obtain
\begin{corollary}\label{beforeKun} For $u\in \mathcal{P}_{\adm,p}(G)$ and $v\in \mathcal{P}_{\adm,q}(G)$, we have
$$\partial^{\cell}(u\times v)=\left(\partial^{\cell}u\right)\times v+(-1)^pu\times\left(\partial^{\cell}v\right).$$
\end{corollary}

\begin{theorem}[K\"{u}nneth formula] For any digraphs $X$ and $Y$, any $n\geq 0$, we have
$$C_n(X\boxdot Y;\R)\cong \bigoplus_{p\geq0,q\geq0;~p+q=n}C_p(X;\R)\otimes C_q(Y;\R).$$
Moreover, we have
$$H_*^{\cell}(X\boxdot Y;\R)\cong H_*^{\cell}(X;\R)\otimes H_*^{\cell}(Y;\R).$$
\end{theorem}

\begin{proof} According to the proof of Lemma \ref{PtimesQ}, any $R\in\mathcal{P}_n(X\boxdot Y)$ admits a decomposition
$$R=P\times Q, \quad\text{for some }P\in\mathcal{P}_{p}(X), ~Q\in\mathcal{P}_{q}(Y),~p+q=n,$$
and
$$\Supp(R)=\Supp(P)\boxdot\Supp(Q).$$

Now let us study the singular cubical realizations of $P$ and $Q$ for $R\in\mathcal{P}_{\adm,n}(X\boxdot Y)$. First, let
$$h=(h_1,h_2):I^{\boxdot n}\rightarrow\Supp(R)=\Supp(P)\boxdot\Supp(Q)$$
be a singular cubic realization of $R$, where
$$h_1: I^{\boxdot n}\rightarrow\Supp(P),\quad h_2:I^{\boxdot n}\rightarrow\Supp(Q).$$
Note that there are $C_{p+q}^p$ kinds of embedding digraph maps
$$i_a: I^{\boxdot p}\hookrightarrow I^{\boxdot (p+q)},\quad a=1,2,\ldots,C_{p+q}^p$$
with the remaining $q$ coordinates being $0$; corresponding the positions of the remaining $q$ coordinates determine another embedding digraph maps
$$j_a: I^{\boxdot q}\hookrightarrow I^{\boxdot (p+q)},\quad a=1,2,\ldots,C_{p+q}^p$$
with the remaining $p$ coordinates being $0$.

Via $i_a$, $j_a$ and $h$, we can define two digraph maps $f_a$ and $g_a$ as follows
\begin{align*}
&f_a=h_1\circ i_a: I^{\boxdot p}\hookrightarrow I^{\boxdot (p+q)}\rightarrow \Supp(P),\\
&g_a=h_2\circ j_a: I^{\boxdot q}\hookrightarrow I^{\boxdot (p+q)}\rightarrow \Supp(Q).
\end{align*}
Let us denote $S_P$ (resp. $S_Q$) and $E_P$ (resp. $E_Q$) to be the starting vertex and ending vertex of $P$ (resp. $Q$) respectively.

Next, let us analyze the properties of all $f_a$, $a=1,\ldots,C_{p+q}^p$ along the following steps. The argument for $g_a$ is the same.

\begin{enumerate}
  \item If $f_a(1,1,\ldots,1)\neq E_P$, then
  \begin{enumerate}
    \item either $(f_a)_*(\omega_p)=0$,
    \item or $(f_a)_*(\omega_p)$ gives another non-$P$ element in $\Omega_p(\Supp(P);\Z)$, it is impossible.
  \end{enumerate}
  \item Let us consider the elementary path components going through $(E_P,S_Q)$ in $R$, that is, the paths start from $(S_P,S_Q)$ and then go along the $X$-directions at the first $p$ steps. Then the condition $h_*(\omega_{p+q})=R$ means that there must exist $a\in\{1,\ldots,C_{p+q}^p\}$ such that
  \begin{equation}\label{f_a}
  f_a(\underbrace{1,1,\ldots,1}_{p})=E_P.
  \end{equation}
  \item We claim that there exists $a$ satisfying \eqref{f_a}, such that $(f_a)_*(\omega_p)\neq0$. Then the minimal condition for $P$ implies that $f_a$ gives a singular cubical realization for $P$. In fact, if all $(f_a)_*(\omega_p)=0$ for $a$ satisfying \eqref{f_a}, the vanishing here means the complete cancellations of elementary $p$-paths in $(f_a)_*(\omega_p)$. By comparing the relation between $f_1$ and $h$, it also implies that there are not elementary path components in $R$ going through $(E_P,S_Q)$, contradiction.
\end{enumerate}

In summary, we have
$$C_n(X\boxdot Y;\R)=\Span_{\R}\{u\times v~\big|~u\in\mathcal{P}_{\adm,p}(X),v\in\mathcal{P}_{\adm,q}(Y)\}/\{\text{admissible relations}\}.$$
Now, let us study the possibly admissible relations among $\mathcal{P}_{\adm,n}(X\boxdot Y)$. First, for $R_i\in\mathcal{P}_{\adm,n}(X\boxdot Y)$, write $R_i=P_i\times Q_i$, $P_i\in\mathcal{P}_{\adm,p}(X),Q_i\in\mathcal{P}_{\adm,q}(Y)$. Assume we have the admissible (linear) relation in $\mathcal{P}_{\adm,n}(X\boxdot Y)$, and write it as
\begin{equation}\label{boxdotlinear}
R_k=\sum_{i=1}^{k-1}a_iR_i,\quad a_i\in\mathbb{R}.
\end{equation}
%By the minimal condition for each $R_i$, we can assume it is a reduced relation that all $R_i$'s have the same starting vertex as well as the same ending vertex.
%Since $R_i=P_i\times Q_i$, and by definition,
%\begin{itemize}
%  \item the starting vertex $S$ of $R_i$ is $(S_{P_i}, S_{Q_i})$;
%  \item the ending vertex $E$ of $R_i$ is $(E_{P_i}, E_{Q_i})$.
%\end{itemize}
%The reduced admissible relation tells us
%\begin{align*}
%&S_{P_1}=S_{P_2}=\cdots=S_{P_k}=:S_P,\quad S_{Q_1}=S_{Q_2}=\cdots=S_{Q_k}:=S_Q;\\
%&E_{P_1}=E_{P_2}=\cdots=E_{P_k}=:E_P,\quad E_{Q_1}=E_{Q_2}=\cdots=E_{Q_k}:=E_Q.
%\end{align*}
%Let $P_i$ has length $p_i$, $i=1,\ldots,k$. Note that if $p_i\neq p_j$ for $i\neq j$, then $R_i=P_i\times Q_i$ doesn't have the same elementary path component as $P_j\times Q_j$. Thus, we arrive at a refined admissible relation from the original relation
%\begin{equation}\label{refine}
%R_k=l_r(R_1,\ldots, R_{l-1}),
%\end{equation}
%where
%\begin{itemize}
%  \item $R_i=P_i\times Q_i$ with starting vertex $(S_P,S_Q)$ and ending vertex $(E_P,E_Q)$, $i=1,\ldots,k$;
%  \item $p_1=p_2=\cdots=p_{l-1}=p_k=:p$, corresponding $Q_i$ has the same length $q=n-p$.
%\end{itemize}
Moreover, we expand $Q_i$ as the linear combinations of s-regular allowed elementary paths
$$Q_k=\sum \epsilon_qe_q,\quad \epsilon_q=1,\text{ or }-1.$$
Then we can pick up the term $P_i\times\epsilon_qe_q$ on both sides of \eqref{boxdotlinear}, and obtain a linear relation among $\{P_i\times e_q\}$.
%$$P_k\times\epsilon_qe_q=\sum_{i=1}^{k-1}a_iP_i\times \epsilon_qe_q$$
Similarly, we can expand $P_k$ and obtain a linear relation to express $\epsilon_pe_p\times Q_k$.

Thus, we have an isomorphism
$$C_n(X\boxdot Y;\R)\cong \bigoplus_{p\geq0,q\geq0;~p+q=n}C_p(X;\R)\otimes C_q(Y;\R)$$
which is given on the generators by
$$R=P\times Q\longmapsto P\otimes Q.$$
The remaining follow from Corollary \ref{beforeKun} and the standard argument.
\end{proof}

\bigskip

\section{More Examples}\label{4}

\subsection{The circulant graph}\label{circ}
The so-called circulant graphs (also known as loop networks), has been much studied. Let us recall the definition.

\begin{definition}\label{circulant}
Let $1\leq\gamma_1\leq\cdots\leq\gamma_d\leq\left[\frac{n}{2}\right]$ be positive integers. A circulant graph $C_{n}^{\gamma_1,\ldots,\gamma_d}$ is the $2d$-regular graph with $n$ vertices labelled by $0,1\ldots,n-1$ such that each vertex $\nu$ is connected to $\nu\pm\gamma_i$ $\mod n$ for all $i\in\{1,2,\cdots,d\}$.
\end{definition}

It is well known that each graph $C_n^{\Gamma}$ (where $\Gamma=\{1,\gamma_2,\ldots,\gamma_d\}$) is isomorphic to the $d$-dimensional discrete torus $\Z^d/\Lambda_{\Gamma}\Z^d$ with $\Lambda_{\Gamma}$ being the $d\times d$ matrix:
\begin{align*}
&(\Lambda_{\Gamma})_{11}=n,\quad \Lambda_{1i}=-\gamma_i,~i=2,3,\ldots,d\\
&(\Lambda_{\Gamma})_{ij}=\delta_{ij},\quad i=2,\ldots,d; ~j=1,\ldots, d.
\end{align*}

\begin{example}[$C_5^{1,2}$ and $C_7^{1,3}$] The circulant graphs $C_5^{1,2}$ and $C_7^{1,3}$ are given as follows:
\begin{figure}[H]
	\centering
	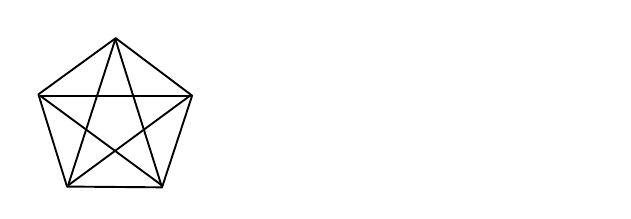
	%\caption[]{$\Supp(e_{0123})$}
	\label{Fig:The circulant graphs $C_5{1,2}$ and $C_7^{1,3}$}
\end{figure}
\end{example}

In our situation, we study the circulant graph as a digraph by considering the directions of edges: for each $i=1,\ldots,d$,
$$\nu-\gamma_i~\rightarrow~\nu~\rightarrow~\nu+\gamma_i\quad\mod n.$$
Then the above circulant graphs, which we still denote by $C_5^{1,2}$ and $C_7^{1,3}$, become
\begin{figure}[H]
	\centering
	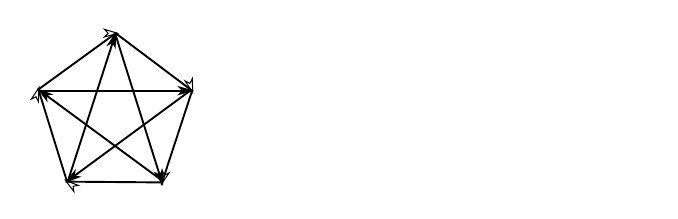
	%\caption[]{$\Supp(e_{0123})$}
	\label{Fig:The circulant digraphs $C_5{1,2}$ and $C_7^{1,3}$}
\end{figure}

Now, let us compute the cellular homologies of the digraphs $C_5^{1,2}$ and $C_7^{1,3}$.
\begin{example}[The digraph $C_5^{1,2}$]
For $n\in\mathbb{N}$, the sets of admissible paths $\mathcal{P}_{\adm,n}(C_5^{1,2})$, which are the same as the sets of minimal paths $\mathcal{P}_n(C_5^{1,2})$, are given by
\begin{align*}
&\mathcal{P}_{\adm,0}(C_5^{1,2})=\mathcal{P}_0(C_5^{1,2})=V(C_5^{1,2});\\
&\mathcal{P}_{\adm,1}(C_5^{1,2})=\mathcal{P}_1(C_5^{1,2})=E(C_5^{1,2});\\
&\mathcal{P}_{\adm,2}(C_5^{1,2})=\mathcal{P}_2(C_5^{1,2})=\{e_{012},~e_{123},~e_{234},~e_{340},~e_{401},~e_{013}-e_{023},~e_{124}-e_{134},\\
&\qquad\qquad\qquad\qquad\qquad\qquad e_{230}-e_{240},~e_{301}-e_{341},~e_{402}-e_{412}\};\\
&\mathcal{P}_{\adm,3}(C_5^{1,2})=\mathcal{P}_3(C_5^{1,2})=\{e_{0123},~e_{1234},~e_{2340},~e_{3401},~e_{4012},~e_{0124}-e_{0134}+e_{0234},\\
&\qquad\qquad\qquad\qquad\qquad\qquad e_{1230}-e_{1240}+e_{1340},~e_{2341}-e_{2301}+e_{2401},\\
&\qquad\qquad\qquad\qquad\qquad\qquad e_{3402}-e_{3412}+e_{3012},~e_{4013}-e_{4023}+e_{4123}\};\\
&\mathcal{P}_{\adm,4}(C_5^{1,2})=\mathcal{P}_4(C_5^{1,2})=\{e_{01234},~ e_{12340},~ e_{23401},~ e_{34012},~e_{40123}\};\\
&\mathcal{P}_{\adm,n}(C_5^{1,2})=\mathcal{P}_n(C_5^{1,2})=0,\quad n\geq5.
\end{align*}
One can check that $\left(C_5^{1,2},\cup_{n=0}^4\mathcal{P}_{n}\right)$ gives the path complex structure as well as CW structure of $C_5^{1,2}$, due to the absence of the admissible relations. After a careful calculation, we have
\begin{align*}
&H_0^{\cell}(C_5^{1,2};\R)=H_0(C_5^{1,2};\R)\cong\R;\\
&H_1^{\cell}(C_5^{1,2};\R)=H_1(C_5^{1,2};\R)=\Span_{\R}\bigg\{[e_{01}+e_{12}+e_{23}+e_{34}+e_{40}]\bigg\};\\
%&\qquad\qquad\qquad\qquad\qquad\qquad\qquad [e_{02}+e_{24}+e_{41}+e_{13}+e_{30}]\bigg\};\\
%&H_2^{\cell}(C_5^{1,2};\R)=H_2(C_5^{1,2};\R)=\Span_{\R}\bigg\{[(e_{301}-e_{341})+(e_{412}-e_{402})\\
%&\qquad\qquad\qquad\qquad\qquad\qquad\qquad +(e_{023}-e_{013})+(e_{134}-e_{124})+(e_{240}-e_{230})]\bigg\};\\
&H_n^{\cell}(C_5^{1,2};\R)=H_n(C_5^{1,2};\R)=0,\quad n\geq 2.
\end{align*}
\end{example}

\begin{remark}\label{infinitepc} (1) One can choose different representatives for $H_1^{\cell}(C_5^{1,2};\R)$. Actually,
\begin{align*}
[e_{01}+e_{12}+e_{23}+e_{34}+e_{40}]=~&[e_{01}+e_{13}+e_{30}]=[e_{12}+e_{24}+e_{41}]\\
=~&[e_{23}+e_{30}+e_{02}]=[e_{34}+e_{41}+e_{13}]\\
=~&[e_{40}+e_{02}+e_{24}];\\
[e_{02}+e_{24}+e_{41}+e_{13}+e_{30}]=~&2[e_{01}+e_{12}+e_{23}+e_{34}+e_{40}].
\end{align*}

It is easy to see that the following 2-chain
$$(e_{013}-e_{023})+(e_{124}-e_{134})+(e_{230}-e_{240})+(e_{341}-e_{301})+(e_{402}-e_{412})$$
is closed. But it is also the boundary of $e_{0123}+e_{1234}+e_{2340}+e_{3401}+e_{4012}$. One can further check the 3-closed chains are also trivial. From the point of view of discrete torus, $C_5^{1,2}$ could be regarded as a ``solid" discrete torus.

(2) If we consider GLMY's regular condition of path complex in Remark \ref{regular}, one can find that for any $k\in\mathbb{N}$,
\begin{align*}
&e_{{\underbrace{01234\cdots01234}_{k}0}},~ e_{{\underbrace{01234\cdots01234}_{k}01}},~ e_{{\underbrace{01234\cdots01234}_{k}012}}, \\ &e_{{\underbrace{01234\cdots01234}_{k}0123}}, ~e_{{\underbrace{01234\cdots01234}_{k}01234}}\in\Omega_*(G).
\end{align*}
It is not easy to completely determine GLMY path complex for $C_5^{1,2}$. But with the strongly regular condition, the computation of path complex and cellular complex is much simpler.
\end{remark}

\begin{example}[The digraph $C_7^{1,3}$] %Let us look at the directed circulant graph $C_7^{1,3}$.
%\begin{figure}[H]
%	\centering
%	\input{diC713.pdf_tex}
%	%\caption[]{$\Supp(e_{0123})$}
%	\label{Fig:The circulant graphs $C_5{1,2}$ and $C_7^{1,2}$}
%\end{figure}
The sets $\mathcal{P}_{\adm,n}(C_7^{1,3})$, which are the same as $\mathcal{P}_n(C_7^{1,3})$, are given by
\begin{align*}
&\mathcal{P}_{\adm,0}(C_7^{1,3})=\mathcal{P}_0(C_7^{1,3})=V(C_7^{1,3});\\
&\mathcal{P}_{\adm,1}(C_7^{1,3})=\mathcal{P}_1(C_7^{1,3})=E(C_7^{1,3});\\
&\mathcal{P}_{\adm,2}(C_7^{1,3})=\mathcal{P}_2(C_7^{1,3})=\{e_{014}-e_{034},~e_{125}-e_{145},~e_{236}-e_{256},~e_{340}-e_{360},\\
&\qquad\qquad\qquad\qquad\qquad\qquad  e_{451}-e_{401},~e_{562}-e_{512},~e_{603}-e_{623}\};\\
&\mathcal{P}_{\adm,n}(C_7^{1,3})=\mathcal{P}_n(C_7^{1,3})=0,\quad n\geq3.
\end{align*}
One can check that $\left(C_7^{1,3},\cup_{n=0}^2\mathcal{P}_{n}\right)$ gives the path complex structure as well as CW structure of $C_7^{1,3}$. Furthermore, we have
\begin{align*}
&H_0^{\cell}(C_7^{1,3};\R)=H_0(C_7^{1,3};\R)\cong\R;\\
&H_1^{\cell}(C_7^{1,3};\R)=H_1(C_7^{1,3};\R)=\Span_{\R}\bigg\{[e_{03}+e_{36}+e_{60}],~\\
&\qquad\qquad\qquad\qquad\qquad\qquad\quad [e_{01}+e_{12}+e_{23}+e_{34}+e_{45}+e_{56}+e_{60}]\bigg\};\\
&H_2^{\cell}(C_7^{1,3};\R)=H_2(C_7^{1,3};\R)=\Span_{\R}\bigg\{[(e_{014}-e_{034})+(e_{125}-e_{145})\\
&\qquad\qquad\qquad\qquad\qquad\qquad\quad +(e_{236}-e_{256})+(e_{340}-e_{360})+(e_{451}-e_{401})\\
&\qquad\qquad\qquad\qquad\qquad\qquad\quad +(e_{562}-e_{512})+(e_{603}-e_{623})]\bigg\};\\
&H_n^{\cell}(C_7^{1,3};\R)=H_n(C_7^{1,3};\R)=0,\quad n\geq 3.
\end{align*}
\end{example}

%\subsection{The multipartite digraph}
%
%\begin{definition}A multipartite digraph $K_{n_1,\ldots,n_l}$ is defined as follows. Its set of vertices $V$ is a disjoint union of $l>1$ subsets $V_1,\ldots,V_l$ so that $|V_i|=n_i>1$, and the arrows are defined by
%$$x~\longrightarrow~y \quad\Longleftrightarrow\quad x\in V_i,~ y\in V_j \text{ for some }i<j.$$
%\end{definition}

\subsection{The Johnson digraph}
\begin{definition}Fix an integer $k$, $1\leq k\leq n$. The Johnson digraph $\overrightarrow{J}(n,k)$ is defined as a digraph whose vertices are all $k$-element subsets of the set $V=\{1,2,\ldots,n\}$, and the arrows are defined as follows: if $a$ and $b$ are two $k$-element subsets of $V$, then $a\rightarrow b$ in $\overrightarrow{J}(n,k)$ if $a$ and $b$ have the following form (unordered):
$$a=\{c_0,c_1,\ldots,c_{k-1}\},\quad b=\{c_0',c_1,\ldots,c_{k-1}\},\quad c_0'<c_0.$$
\end{definition}

\begin{example}[Johnson digraph $\overrightarrow{J}(4,2)$] For $n=4$, $k=2$, we get the following left digraph $\overrightarrow{J}(4,2)$. For convenience, we relabel its vertices as the right picture shows.
\begin{figure}[H]
	\centering
	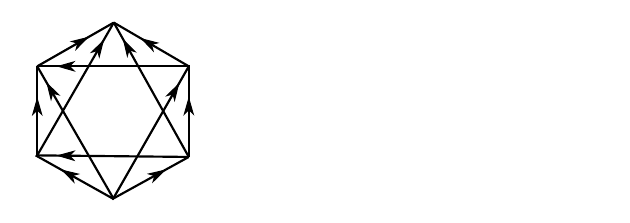
	%\caption[]{$\Supp(e_{0123})$}
	%\label{}
\end{figure}
The sets $\mathcal{P}_{\adm,n}=\mathcal{P}_{\adm,n}(\overrightarrow{J}(4,2))$, which are the same as $\mathcal{P}_n(\overrightarrow{J}(4,2))$, are given by
\begin{align*}
&\mathcal{P}_{\adm,0}=\mathcal{P}_0=V(\overrightarrow{J}(4,2));\\
&\mathcal{P}_{\adm,1}=\mathcal{P}_1=E(\overrightarrow{J}(4,2));\\
&\mathcal{P}_{\adm,2}=\mathcal{P}_2=\{e_{210},~e_{310},~e_{420},~e_{430},~e_{521},~e_{531},~e_{421}-e_{431},~e_{542},~e_{543},\\
&\qquad\qquad\qquad\quad e_{5i0}-e_{5j0},~ i,j=1,2,3,4, \text{ and }i\neq j \text{ (up to a sign)}\};\\
&\mathcal{P}_{\adm,3}=\mathcal{P}_3=\{e_{5310},~e_{5420}\}.
\end{align*}
There exist admissible relations among $\{e_{5i0}-e_{5j0},~ i,j=1,2,3,4, \text{ and }i\neq j \text{ (up to a sign)}\}$. Furthermore, we can compute the cellular homology of $\overrightarrow{J}(4,2)$, and obtain
$$\widetilde{H}_*^{\cell}(\overrightarrow{J}(4,2))=0.$$
The result could also follow from the contractibility of Johnson digraph $\overrightarrow{J}(4,2)$, where the deformation retraction is given by the composition of the following two digraph maps
\begin{align*}
&f_1(5)=f_1(4)=f_1(2)=2,\quad f_1(3)=f_1(1)=1,\quad f_1(0)=0;\\
&f_2(2)=f_2(1)=f_2(0)=0.
\end{align*}
\end{example}

\bigskip

\section{Discussion and Outlook}\label{5}

\subsection{The digraphs from geometry}

There are several ways to discrete a manifold $M$ and get a graph which captures some topological structure of $M$. For example, one can choose a simplicial complex $K$ which triangulates $M$, and its one-skeleton gives us a graph. Such a discretization has a wide range of applications, such as Cheeger-M\"{u}ller theorem, combinatorial quantum field theory.

CMY\cite{GMY1} prove that after twice barycentric subdivisions, the corresponding path homology of the resulting digraph is isomorphic to the simplicial homology of $K$.

Instead, we can consider the digraph from geometry in another way. For simplicity, let $(M,g)$ be an oriented complete Riemannian manifold with metric $g$ and injectivity radius bounded below by a positive constant $b$. First, choose a discrete collection of points $\{x_i\}=V\subset M$ satisfying
\begin{enumerate}
  \item $\dist_g(x_i,x_j)\geq\frac{b}{4}$,~ for any $x_i\neq x_j$;
  \item for every $x\in M$ there is some $x_i\in V$ such that $\dist_g(x_i,x)<\frac{b}{2}$.
\end{enumerate}
Next, define an (undirected) edge structure on $V$ by the relationship
$$x_i\sim x_j,\text{ if } \dist_g(x_i, x_j)<b \text{ and } i\neq j.$$
Then $G=(V,E)$ is a connected graph. Moreover, we consider the direction of each edge which is induced by the orientation of $M$: Given an undirected edge $x_i\sim x_j$, we say $x_i~\rightarrow~x_j$, if the tangent vector at $x_i$ to a distance minimizing geodesic connecting $x_i$ to $x_j$ is given by the orientation on $M$.

\begin{example}[The integer lattice $\Z^n$] Let $\R^n$ be the Euclidean space. Take $V=\mathbb{Z}^n$ and $b=\sqrt{2}$, then we get the standard nearest-neighbor digraph on the given lattice.
\end{example}

We will study the path homology and cellular homology of the digraph from geometry and compare them to the homology of the manifold in a separation paper.

\subsection{The relation to singular cubical homology}

In Remark \ref{regular}, we have the observation: if the digraph is acyclic,
%$$v_0\rightarrow v_1\rightarrow v_2\rightarrow\cdots\rightarrow v_n\rightarrow v_0,$$
%then the strongly regular condition is the same as the regular condition. It follows that
then our restricted path homology is also the same as GLMY path homology. In the sequel, we continue focusing on such special digraphs.

Note that the definition of our cellular homology is related to the singular cubical homology. It is natural to consider the relation between the two homologies. First, if Conjecture \ref{contractible} holds, then by homotopy invariance of both homologies, we have
$$\widetilde{H}_*^{\cell}(\Supp(P))=\widetilde{H}_*^c(\Supp(P))=0,\quad P\in\mathcal{P}_{\adm}(G).$$

Next, we will consider the following three typical digraphs:
\begin{itemize}
  \item $G=\{0\rightarrow 1\rightarrow 3, 0\rightarrow 2\rightarrow 3, 0\rightarrow 3\}$, i.e. Example \ref{basicex};
  \item the digraph $S_k$ in Example \ref{k-square} and Example \ref{summable};
  \item the exotic cube, i.e. Example \ref{exotic}.
\end{itemize}

\subsubsection{Example \ref{basicex} revisited}

In the digraph $G$ of Example \ref{basicex}, we have the following trivial linear relations
\begin{enumerate}
  \item $(e_{013},\Supp(e_{013}))+(-e_{013},\Supp(e_{013}))=0$;
  \item $(e_{013}-e_{023},\Supp(e_{013}-e_{023}))+(e_{023}-e_{013},\Supp(e_{013}-e_{023}))=0$;
  \item $(e_{013},\Supp(e_{013}))+(-e_{023},\Supp(e_{023}))=(e_{013}-e_{023},\Supp(e_{013}-e_{023}))$.\footnote{Although $e_{013}-e_{023}$ is not minimal in $G$, but it admits a singular cubic realization, we want to discuss the relations among these singular cubes.}
\end{enumerate}

For the trivial relations (1) (2), let us look at the following singular $3$-cubes $f_+$, $f_-$, $g_+$, $g_-$:
\begin{align*}
&f_+((0,0))=f_+((0,1))=0,\quad f_+((1,0))=1,\quad f_+((1,1))=3;\\
&f_-((0,0))=f_-((1,0))=0,\quad f_-((0,1))=1,\quad f_-((1,1))=3;\\
&g_+((0,0))=0,\quad g_+((1,0))=1,\quad g_+((0,1))=2,\quad g_+((1,1))=3;\\
&g_-((0,0))=0,\quad g_-((1,0))=2,\quad g_-((0,1))=1,\quad g_-((1,1))=3.
\end{align*}
It is easy to obtain
$$(f_+)_*(\omega_2)+(f_-)_*(\omega_2)=e_{013}+(-e_{013})=0,\quad (g_+)_*(\omega_2)+(g_-)_*(\omega_2)=0.$$
Furthermore, we can construct a non-degenerate $3$-cube $\tilde{f}$:
\begin{align*}
&\tilde{f}((0,0,0))=\tilde{f}((1,0,0))=\tilde{f}((0,0,1))=\tilde{f}((1,0,1))=0,\\
&\tilde{f}((0,1,0))=\tilde{f}((1,1,0))=\tilde{f}((0,1,1))=1,\quad \tilde{f}((1,1,1))=3.
\end{align*}
By definition, $\tilde{f}\circ F_{11}=f_+$, $\tilde{f}\circ F_{31}=f_-$, and other $\tilde{f}\circ F_{j\epsilon}$'s are degenerate. Thus, we have
$$\partial^c(\tilde{f})=f_++f_-.$$

Similarly, for the following non-degenerate $3$-cube $\tilde{g}$:
\begin{align*}
&\tilde{g}((0,0,0))=\tilde{g}((1,0,0))=\tilde{g}((0,0,1))=0,\quad \tilde{g}((1,0,1))=1\\
&\tilde{g}((0,1,0))=\tilde{g}((1,1,0))=\tilde{g}((0,1,1))=2,\quad \tilde{g}((1,1,1))=3.
\end{align*}
one can check $\partial^c(\tilde{g})=g_++g_-$.

\begin{remark}
Note that the above $f_+$ and $f_-$ ($g_+$ and $g_-$ respectively) are related by an automorphism of $I^{\boxdot2}$. Thus it is an interesting problem to study the action of the automorphism group of the cubes on the singular cubes. On the other hand, one can also consider other different singular cubical realizations of the triangle and discuss their relations.
\end{remark}

For the trivial relation (3), let us look at the following three singular $2$-cubes $f_1$, $f_2$, $f_3$:
\begin{align*}
&f_1((0,0))=0,~f_1((0,1))=1,~f((1,0))=f((1,1))=3;\\
&f_2((0,0))=f_2((1,0))=0,~ f_2((0,1))=2,~ f_2((1,1))=3;\quad f_3=g_+.
\end{align*}
Similarly, we have
$$(f_3)_*(\omega_2)=e_{013}-e_{023}=-(f_1)_*(\omega_2)+(f_2)_*(\omega_2).$$
And we can construct a non-degenerate $3$-cubes $h$ in $G$:
\begin{align*}
&h((0,0,0))=h((1,0,0))=h((0,0,1))=0,\quad h((1,0,1))=1\\
&h((0,1,0))=h((0,1,1))=2, \quad  h((1,1,0))=h((1,1,1))=3.
\end{align*}
It satisfies $\partial^c(h)=f_3+f_1-f_2.$

\vskip 0.2cm

\subsubsection{Example \ref{k-square} and Example \ref{summable} revisited}
The digraph $S_k$ in Example \ref{k-square} is contractible, since we have the deformation retraction
$$f:S_k\rightarrow \{S~\rightarrow~ 1\},$$
which are given by
$$f(S)=f(2)=f(3)=\cdots=f(k)=S,\quad f(1)=f(E)=1.$$
By homotopy invariance, we know that both kinds of homologies of $S_k$ are trivial. Note that there are many admissible relations among $\mathcal{P}_{\adm,2}(S_k)$. For simpliicty, we focus our discussions on the following two relations.
\begin{enumerate}
  \item For distinct $i,j,l$, $(e_{SiE}-e_{SjE})-(e_{SiE}-e_{SlE})+(e_{SjE}-e_{SlE})=0$;
  \item $(e_{S1E}-e_{S2E})+(e_{S3E}-e_{S4E})+(e_{S2E}-e_{S3E})+(e_{S4E}-e_{S1E})=0$.
\end{enumerate}

First each square in $S_k$ corresponds to a singular $2$-cube. For example, we consider the cubes
$$f_{ij}: I^{\boxdot 2}\rightarrow S_k,\quad i,j=1,\ldots,k,~i\neq j$$
which are given by
$$f_{ij}((0,0))=S,\quad f_{ij}((1,0))=i,\quad f_{ij}((0,1))=j,\quad f_{ij}((1,1))=E.$$
By definition, for distinct $i,j,l\in\{1,2,\ldots,k\}$, we have
\begin{itemize}
  \item $[(f_{ij})_*-(f_{il})_*+(f_{jl})_*](\omega_2)=0$, i.e. the above admissible relation (1);
  \item $\partial^c(f_{ij}-f_{il}+f_{jl})=0.$
\end{itemize}
%$$\partial^c(f_{ij}-f_{il}+f_{jl})=0.$$
Meanwhile $f_{ij}-f_{il}+f_{jl}$ is also a boundary of a non-degenerate singular $3$-cube $g:I^{\boxdot3}\rightarrow S_k$
\begin{align*}
&g((0,0,0))=g((1,0,0))=g((0,1,0))=g((0,0,1))=S,\\
&g((1,1,0))=i,\quad g((1,0,1))=j,\quad g((0,1,1))=l,\quad g((1,1,1))=E.
\end{align*}
One can check that $\partial^c(g)=f_{ij}-f_{il}+f_{jl}$ $\mod$ (degenerate cubes).\\

Similarly, the four summands in the above admissible relation (2) could be realized by singular $2$-cubes $f_{12}$, $f_{34}$, $f_{23}$ and $f_{41}$ respectively. By definition, we have
$$\partial^c(f_{12}+f_{34}+f_{23}+f_{41})=0.$$
Meanwhile, we can construct two non-degenerate singular $3$-cubes $g_1,g_2:I^{\boxdot3}\rightarrow S_k$:
\begin{align*}
&g_1((0,0,0))=g_1((1,0,0))=g_1((0,1,0))=g_1((0,0,1))=S,\\
&g_1((1,1,0))=1,\quad g_1((1,0,1))=2,\quad g_1((0,1,1))=3,\quad g_1((1,1,1))=E;\\
&g_2((0,0,0))=g_2((1,0,0))=g_2((0,1,0))=g_2((0,0,1))=S,\\
&g_2((1,1,0))=3,\quad g_2((1,0,1))=4,\quad g_2((0,1,1))=1,\quad g_2((1,1,1))=E.
\end{align*}
One can check that $\partial^c(g_1+g_2)=f_{12}+f_{34}+f_{23}+f_{41}$ $\mod$ (degenerate cubes).

\vskip 0.2cm

We learn from this example that the admissible relations may be realized as the boundary conditions in the singular cubical homology. We can also check this observation via Example \ref{summable}, where the digraph is
\begin{figure}[H]
	\centering
	
	%\caption[]{$\Supp(e_{0123})$}
	%\label{Fig:}
\end{figure}
All the admissible paths in
$$\mathcal{P}_{\adm,3}=\{e_{0136}-e_{0236}+e_{0256},~e_{0136}-e_{0146}+e_{0246}-e_{0236},~e_{0146}-e_{0246}+e_{0256}\},$$
could be realized by the following three singular 3-cubes %(for convenience, we embed $I^{\boxdot3}$ in $I^{\boxdot4}$)
%\begin{align*}
%&f: I^{\boxdot3}\boxdot\{0\}\rightarrow \Supp(e_{0136}-e_{0236}+e_{0256}),\\
%%&\qquad  f_*(\omega_3)=e_{0136}-e_{0236}+e_{0256};\\
%&g: I^{\boxdot2}\boxdot\{0\}\boxdot I\rightarrow\Supp(e_{0136}-e_{0146}+e_{0246}-e_{0236}),\\
%%&\qquad  g_*(\omega_3)=e_{0136}-e_{0146}+e_{0246}-e_{0236};\\
%&h: I\boxdot\{0\}\boxdot I^{\boxdot2}\rightarrow\Supp(e_{0146}-e_{0246}+e_{0256}),
%%&\qquad  h_*(\omega_3)=e_{0146}-e_{0246}+e_{0256}.
%\end{align*}
\begin{align*}
&f: I^{\boxdot3}\rightarrow \Supp(e_{0136}-e_{0236}+e_{0256}),\\
%&\qquad  f_*(\omega_3)=e_{0136}-e_{0236}+e_{0256};\\
&g: I^{\boxdot3}\rightarrow\Supp(e_{0136}-e_{0146}+e_{0246}-e_{0236}),\\
%&\qquad  g_*(\omega_3)=e_{0136}-e_{0146}+e_{0246}-e_{0236};\\
&h: I^{\boxdot3}\rightarrow\Supp(e_{0146}-e_{0246}+e_{0256}),
%&\qquad  h_*(\omega_3)=e_{0146}-e_{0246}+e_{0256}.
\end{align*}
which are explicitly given by
%\begin{align*}
%&f((0,0,0,0))=0,\quad f((1,0,0,0))=1,\quad f((0,1,0,0))=2,\quad f((0,0,1,0))=5,\\
%&f((1,1,0,0))=3,\quad f((1,0,1,0))=6,\quad f((0,1,1,0))=5,\quad f((1,1,1,0))=6;\\
%&g((0,0,0,0))=0,\quad g((1,0,0,0))=1,\quad g((0,1,0,0))=2,\quad g((0,0,0,1))=2,\\
%&g((1,1,0,0))=3,\quad g((1,0,0,1))=4,\quad g((0,1,0,1))=4,\quad g((1,1,0,1))=6;\\
%&h((0,0,0,0))=0,\quad h((1,0,0,0))=1,\quad h((0,0,0,1))=2,\quad h((0,0,1,0))=5,\\
%&h((1,0,0,1))=4,\quad h((1,0,1,0))=6,\quad h((0,0,1,1))=5,\quad h((1,0,1,1))=6.
%\end{align*}
\begin{align*}
&f((0,0,0))=0,\quad f((1,0,0))=1,\quad f((0,1,0))=2,\quad f((0,0,1))=5,\\
&f((1,1,0))=3,\quad f((1,0,1))=6,\quad f((0,1,1))=5,\quad f((1,1,1))=6;\\
&g((0,0,0))=0,\quad g((1,0,0))=1,\quad g((0,1,0))=2,\quad g((0,0,1))=1,\\
&g((1,1,0))=3,\quad g((1,0,1))=3,\quad g((0,1,1))=4,\quad g((1,1,1))=6;\\
&h((0,0,0))=0,\quad h((1,0,0))=2,\quad h((0,0,1))=5,\quad h((0,1,0))=1,\\
&h((1,1,0))=5,\quad h((1,0,1))=4,\quad h((0,1,1))=6,\quad h((1,1,1))=6.
\end{align*}
One can check that
\begin{align*}
&f_*(\omega_3)-g_*(\omega_3)-h_*(\omega_3)\\
=~&(e_{0136}-e_{0236}+e_{0256})-(e_{0136}-e_{0146}+e_{0246}-e_{0236})-(e_{0146}-e_{0246}+e_{0256})=0
\end{align*}
%where $\omega_3^1$, $\omega_3^2$ and $\omega_3^3$ are the lift of $\omega_3$ in $I^{\boxdot3}\boxdot\{0\}$, $I^{\boxdot2}\boxdot\{0\}\boxdot I$ and $I\boxdot\{0\}\boxdot I^{\boxdot2}$ respectively. One can check that
and also compute by definition
$$\partial^c(f-g+h)=0.$$
Meanwhile, we can construct $\phi:I^{\boxdot4}\rightarrow G$ as follows:
\begin{align*}
&\phi((0,0,0,0))=0,\quad \phi((1,0,0,0))=1,\quad \phi((0,1,0,0))=2,\quad \phi((0,0,1,0))=5,\\
&\phi((1,1,0,0))=3,\quad \phi((1,0,1,0))=6,\quad \phi((0,1,1,0))=5,\quad \phi((1,1,1,0))=6;\\
&\phi((0,0,0,0))=0,\quad \phi((1,0,0,0))=1,\quad \phi((0,1,0,0))=2,\quad \phi((0,0,0,1))=1,\\
&\phi((1,1,0,0))=3,\quad \phi((1,0,0,1))=3,\quad \phi((0,1,0,1))=4,\quad \phi((1,1,0,1))=6;\\
&\phi((0,0,0,0))=0,\quad \phi((0,1,0,0))=2,\quad \phi((0,0,1,0))=5,\quad \phi((0,0,0,1))=1,\\
&\phi((0,1,1,0))=5,\quad \phi((0,1,0,1))=4,\quad \phi((0,0,1,1))=6,\quad \phi((0,1,1,1))=6;\\
&\phi((1,1,1,1))=6.
\end{align*}
We write several repeated terms to make the following observations clearer
\begin{align*}
&f=\phi\circ F_{40}, \quad \text{via }~ F_{40}:I^{\boxdot3}\cong I^{\boxdot 3}\boxdot\{0\};\\
&g=\phi\circ F_{30}, \quad\text{via }~ F_{30}:I^{\boxdot3}\cong I^{\boxdot2}\boxdot\{0\}\boxdot I;\\
&h=\phi\circ F_{10}, \quad\text{via }~ F_{10}:I^{\boxdot3}\cong \{0\}\boxdot I^{\boxdot3}.
\end{align*}
%\begin{align*}
%&f=\phi\circ F_{40}=\phi\big|_{I^{\boxdot3}\boxdot\{0\}}, \quad \text{via }~ F_{40}:I^{\boxdot3}\cong I^{\boxdot 3}\boxdot\{0\},\\
%&g=\phi\circ F_{30}=\phi\big|_{I^{\boxdot2}\boxdot\{0\}\boxdot I}, \quad\text{via }~ F_{30}:I^{\boxdot3}\cong I^{\boxdot2}\boxdot\{0\}\boxdot I\\
%&h=\phi\circ F_{10}=\phi\big|_{\{0\}\boxdot I^{\boxdot3}}, \quad\text{via }~ F_{20}:I^{\boxdot3}\cong \{0\}\boxdot I^{\boxdot3}.
%\end{align*}
To compute $\partial^c\phi$, we also need to know
$$\phi\circ F_{20}~\text{ and }~ \phi\circ F_{j1},~j=1,2,3,4.$$
One can easily check one by one that these singular $3$-cubes are degenerate. Thus, by definition,
\begin{align*}
\partial^c{\phi}=~&\sum_{j=1}^4(-1)^j\left(\phi\circ F_{j0}-\phi\circ F_{j1}\right)\\
                =~&f-g-h \quad\mod\text{ (degenerate cubes)}.
\end{align*}

\vskip 0.2cm

\subsubsection{Example \ref{exotic} revisited}

In \cite{GJM}, a non-trivial 2-cycle $\phi$ (where they use the notation $\phi^{\Box}$) for the singular cubical homology of the exotic cube is constructed to be a linear combination of 8 non-generate singular $2$-cubes:
$$\phi=\sum_{i=1}^8\phi_i.$$
%$$\phi^{\Box}=\sum_{i=1}^8\phi_i^{\Box}.$$
In particular, the 8 singular $2$-cubes exactly give the singular cubical realizations of the 8 minimal 2-faces in $\mathcal{P}_{\adm,2}(\mbox{exotic cube})$ respectively. And
%$$\phi^{\Box}(\omega_2)=\sum_{i=1}^8\phi_i^{\Box}(\omega_2)$$
$$\phi_*(\omega_2)=\sum_{i=1}^8(\phi_i)_*(\omega_2)$$
gives us the non-trivial 2-cycle in the cellular homology which generates $H_2^{\cell}(\text{exotic cube};\R)$.

The singular cubical chains are so large that it is quite difficult to write down the chain complex and compute the corresponding homologies. But our cellular homologies are much easier. With the above observations, we have the following conjecture.

\begin{conjecture}\label{conj} If $G$ is an acyclic digraph, then
$$H_*^{\cell}(G;\R)\cong H_*^c(G;\R).$$
In particular, the admissible relations among the admissible pairs are the boundary conditions in the singular cubical homology.
\end{conjecture}

\bigskip

\begin{appendices}

\section{Appendix. The admissible minimal $3$-paths}\label{appendixsection}\label{AppendixA}

In this appendix section, we will list all the admissible minimal $3$-paths and the deformation retractions to the corresponding face components.

First, as we explain before, there are finitely many admissible pairs $(P,\Supp(P))$ for each $n\geq 0$. For $P\in\mathcal{P}_{\adm,n}(G)$, we have
\begin{itemize}
  \item $|V(P)|=|V(\Supp(P))|\leq |V(I^{\boxdot n})|=2^{n}$.
  \item Let $\NE(P)$ be the number of elementary path components in $P$. Then $f_*(\omega_n)=P$ implies
$$\NE(P)\leq \NE(\omega_n)=n!.$$
  \item Let $\NF(P)$ be the number of $(n-1)$-face components of $P$ in $\Supp(P)$. Then Lemma \ref{elementarylem} implies
$$\NF(P)\leq \NF(\omega_n)=2n.$$
  \item If we reverse all the edges in $\Supp(P)$, that is,
  $$\text{for any }a\rightarrow b\text{ in }\Supp(P),\qquad\rightsquigarrow\qquad b\rightarrow a,$$
  then obviously the resulting path $\tilde{P}$ is also admissible.
\end{itemize}
Thus, we can use these number constraints to analyze the set of admissible minimal $3$-paths up to the above global orientation of edges. Finally, we arrive at the following eight admissible minimal $3$-paths.

For simplicity, we will use the following notations
\begin{itemize}
  \item $K_3(a,b,c)$: the triangle $\{a\rightarrow b\rightarrow c, ~a\rightarrow c\}$;
  \item $S(a,b,c,d)$: the square $\{a\rightarrow b\rightarrow d, ~a\rightarrow c\rightarrow d\}$.
\end{itemize}

\begin{example}[$\NE(P)=1$] The admissible pair $(P,\Supp(P))$ with $\NE(P)=1$ is given by $P=e_{0123}$ and
\begin{figure}[H]
	\centering
	%% Creator: Inkscape 1.0.1 (3bc2e813f5, 2020-09-07), www.inkscape.org
%% Accompanies image file '0123.pdf' (pdf, eps, ps)
\begingroup%
  \makeatletter%
  \providecommand\color[2][]{%
    \errmessage{(Inkscape) Color is used for the text in Inkscape, but the package 'color.sty' is not loaded}%
    \renewcommand\color[2][]{}%
  }%
  \providecommand\transparent[1]{%
    \errmessage{(Inkscape) Transparency is used (non-zero) for the text in Inkscape, but the package 'transparent.sty' is not loaded}%
    \renewcommand\transparent[1]{}%
  }%
  \providecommand\rotatebox[2]{#2}%
  \newcommand*\fsize{\dimexpr\f@size pt\relax}%
  \newcommand*\lineheight[1]{\fontsize{\fsize}{#1\fsize}\selectfont}%
  \ifx\svgwidth\undefined%
    \setlength{\unitlength}{237.51158947bp}%
    \ifx\svgscale\undefined%
      \relax%
    \else%
      \setlength{\unitlength}{\unitlength * \real{\svgscale}}%
    \fi%
  \else%
    \setlength{\unitlength}{\svgwidth}%
  \fi%
  \global\let\svgwidth\undefined%
  \global\let\svgscale\undefined%
  \makeatother%
  \begin{picture}(1,0.21212938)%
    \lineheight{1}%
    \setlength\tabcolsep{0pt}%
    \put(0,0){\includegraphics[width=\unitlength,page=1]{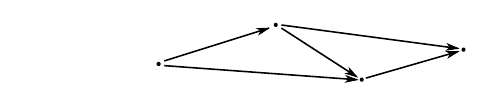}}%
    \put(-0.00093796,0.08014644){\makebox(0,0)[lt]{\lineheight{1.25}\smash{\begin{tabular}[t]{l}$\Supp(e_{0123})=$\end{tabular}}}}%
    \put(0.32766127,0.1153773){\makebox(0,0)[lt]{\lineheight{1.25}\smash{\begin{tabular}[t]{l}$0$\end{tabular}}}}%
    \put(0.53483801,0.1905926){\makebox(0,0)[lt]{\lineheight{1.25}\smash{\begin{tabular}[t]{l}$1$\end{tabular}}}}%
    \put(0.71206984,0.00278773){\makebox(0,0)[lt]{\lineheight{1.25}\smash{\begin{tabular}[t]{l}$2$\end{tabular}}}}%
    \put(0.95719778,0.10477412){\makebox(0,0)[lt]{\lineheight{1.25}\smash{\begin{tabular}[t]{l}$3$\end{tabular}}}}%
  \end{picture}%
\endgroup%

	%\caption[]{$\Supp(e_{0123})$}
	%\label{Fig:Q1}
\end{figure}
The singular cubical map (non-unique) $f:I^{\boxdot3}\rightarrow\Supp(P)$ could be given by
\begin{align*}
f((0,0,0))=0,\quad f((1,0,0))=1,\quad f((0,1,0))=2,\quad f((0,0,1))=2,\\
f((1,1,0))=2,\quad f((1,0,1))=3,\quad f((0,1,1))=2,\quad f((1,1,1))=3.
\end{align*}
We also have the following deformation retractions:
\begin{enumerate}
  \item $r_1:\Supp(P)\rightarrow \im\left(f\big|_{I^{\boxdot 2}\boxdot\{0\}}\right)=K_3(0,1,2)$,
  $$r_1(3)=r_1(2)=2,~r_1(1)=1,~r_1(0)=0;$$
  %\item $r_2:\Supp(P)\rightarrow \im\left(f\big|_{I\boxdot\{0\}\boxdot I}\right)$,
  %$r_{2}(3)=r_1(2)=2,~r_1(1)=1,~r_1(0)=0$;
  \item $r_2:\Supp(P)\rightarrow \im\left(f\big|_{\{0\}\boxdot I^{\boxdot 2}}\right)=\{0\rightarrow 2\}$,
  $$r_2(3)=r_2(2)=2,~r_2(1)=r_2(0)=0;$$
  \item $r_3:\Supp(P)\rightarrow \im\left(f\big|_{I^{\boxdot 2}\boxdot\{1\}}\right)=\{2\rightarrow 3\}$,
  $$r_3(3)=3,~r_3(2)=r_3(1)=r_3(0)=2;$$
  \item $r_4:\Supp(P)\rightarrow \im\left(f\big|_{I\boxdot\{1\}\boxdot I}\right)=\{2\rightarrow 3\}$,
  $$r_4(3)=3,~r_4(2)=r_4(1)=r_4(0)=2;$$
  \item $r_5:\Supp(P)\rightarrow \im\left(f\big|_{\{1\}\boxdot I^{\boxdot 2}}\right)=K_3(1,2,3)$,
  $$r_5(3)=3,~r_5(2)=2,~r_5(1)=r_5(0)=1.$$
\end{enumerate}
Note that $\im\left(f\big|_{I\boxdot\{0\}\boxdot I}\right)=S(0,1,2,3)$, there is no digraph map from $\Supp(P)$ onto $S(0,1,2,3)$. But it is easy to construct the digraph maps from $\Supp(P)$ onto
$$\{0\rightarrow1\},~\{0\rightarrow2\},~\{1\rightarrow3\},~\{2\rightarrow3\},~\{0\rightarrow1\rightarrow3\},~\{0\rightarrow2\rightarrow 3\}.$$
\end{example}

\begin{example}[$\NE(P)=2$] One of the admissible pairs $(P,\Supp(P))$ with $\NE(P)=2$ is given by $P=e_{0124}-e_{0134}$ and
\begin{figure}[H]
	\centering
	%% Creator: Inkscape 1.0.1 (3bc2e813f5, 2020-09-07), www.inkscape.org
%% Accompanies image file '0124.pdf' (pdf, eps, ps)
\begingroup%
  \makeatletter%
  \providecommand\color[2][]{%
    \errmessage{(Inkscape) Color is used for the text in Inkscape, but the package 'color.sty' is not loaded}%
    \renewcommand\color[2][]{}%
  }%
  \providecommand\transparent[1]{%
    \errmessage{(Inkscape) Transparency is used (non-zero) for the text in Inkscape, but the package 'transparent.sty' is not loaded}%
    \renewcommand\transparent[1]{}%
  }%
  \providecommand\rotatebox[2]{#2}%
  \newcommand*\fsize{\dimexpr\f@size pt\relax}%
  \newcommand*\lineheight[1]{\fontsize{\fsize}{#1\fsize}\selectfont}%
  \ifx\svgwidth\undefined%
    \setlength{\unitlength}{215.15112449bp}%
    \ifx\svgscale\undefined%
      \relax%
    \else%
      \setlength{\unitlength}{\unitlength * \real{\svgscale}}%
    \fi%
  \else%
    \setlength{\unitlength}{\svgwidth}%
  \fi%
  \global\let\svgwidth\undefined%
  \global\let\svgscale\undefined%
  \makeatother%
  \begin{picture}(1,0.34050514)%
    \lineheight{1}%
    \setlength\tabcolsep{0pt}%
    \put(0,0){\includegraphics[width=\unitlength,page=1]{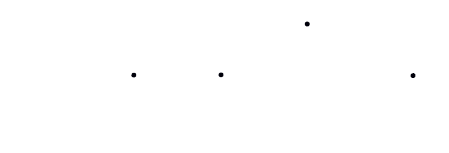}}%
    \put(-0.0015262,0.16452058){\makebox(0,0)[lt]{\lineheight{1.25}\smash{\begin{tabular}[t]{l}$\Supp(P)=$\end{tabular}}}}%
    \put(0.29275559,0.21062808){\makebox(0,0)[lt]{\lineheight{1.25}\smash{\begin{tabular}[t]{l}$0$\end{tabular}}}}%
    \put(0.54379348,0.16238082){\makebox(0,0)[lt]{\lineheight{1.25}\smash{\begin{tabular}[t]{l}$1$\end{tabular}}}}%
    \put(0.66221435,0.32140433){\makebox(0,0)[lt]{\lineheight{1.25}\smash{\begin{tabular}[t]{l}$2$\end{tabular}}}}%
    \put(0.67377717,0.00266865){\makebox(0,0)[lt]{\lineheight{1.25}\smash{\begin{tabular}[t]{l}$3$\end{tabular}}}}%
    \put(0,0){\includegraphics[width=\unitlength,page=2]{0124.pdf}}%
    \put(0.95795054,0.16495191){\makebox(0,0)[lt]{\lineheight{1.25}\smash{\begin{tabular}[t]{l}$4$\end{tabular}}}}%
  \end{picture}%
\endgroup%

	%\caption[]{$\Supp(e_{0123})$}
	%\label{Fig:Q1}
\end{figure}
Another one is given by reversing all the edges.

The singular cubical map (non-unique) $f:I^{\boxdot3}\rightarrow\Supp(P)$ could be given by
\begin{align*}
f((0,0,0))=0,\quad f((1,0,0))=1,\quad f((0,1,0))=2,\quad f((0,0,1))=3,\\
f((1,1,0))=2,\quad f((1,0,1))=3,\quad f((0,1,1))=4,\quad f((1,1,1))=4.
\end{align*}
We also have the following deformation retractions:
\begin{enumerate}
  \item $r_1:\Supp(P)\rightarrow \im\left(f\big|_{I^{\boxdot 2}\boxdot\{0\}}\right)=K_3(0,1,2)$,
  $$r_1(4)=r_1(2)=2,~r_1(3)=r_1(1)=1,~r_1(0)=0;$$
  \item $r_2:\Supp(P)\rightarrow \im\left(f\big|_{I\boxdot\{0\}\boxdot I}\right)=K_3(0,1,3)$,
  $$r_2(4)=r_2(3)=3,~r_2(2)=r_2(1)=1,~r_2(0)=0;$$
  \item $r_3:\Supp(P)\rightarrow \im\left(f\big|_{\{0\}\boxdot I^{\boxdot 2}}\right)=S(0,2,3,4)$,
  $$r_3(4)=4,~r_3(3)=3,~r_3(2)=2,~r_3(1)=r_3(0)=0;$$
  \item $r_4:\Supp(P)\rightarrow \im\left(f\big|_{I^{\boxdot 2}\boxdot\{1\}}\right)=\{3\rightarrow 4\}$,
  $$r_4(4)=r_4(2)=4,~r_4(3)=r_4(1)=r_4(0)=3;$$
  \item $r_5:\Supp(P)\rightarrow \im\left(f\big|_{I\boxdot\{1\}\boxdot I}\right)=\{2\rightarrow 4\}$,
  $$r_5(4)=r_5(3)=4,~r_5(2)=r_5(1)=r_5(0)=2;$$
  \item $r_6:\Supp(P)\rightarrow \im\left(f\big|_{\{1\}\boxdot I^{\boxdot 2}}\right)=S(1,2,3,4)$,
  $$r_6(4)=4,~r_6(3)=3,~r_6(2)=2,~r_6(1)=r_6(0)=1.$$
\end{enumerate}
\end{example}

\begin{example}[$\NE(P)=3$, $\NF(P)=4$] The admissible pair $(P,\Supp(P))$ with $\NE(P)=3$, $\NF(P)=4$ is given by
$$P=e_{0124}-e_{0134}+e_{0234},$$
\begin{figure}[H]
	\centering
	%% Creator: Inkscape 1.0.1 (3bc2e813f5, 2020-09-07), www.inkscape.org
%% Accompanies image file '0234.pdf' (pdf, eps, ps)
\begingroup%
  \makeatletter%
  \providecommand\color[2][]{%
    \errmessage{(Inkscape) Color is used for the text in Inkscape, but the package 'color.sty' is not loaded}%
    \renewcommand\color[2][]{}%
  }%
  \providecommand\transparent[1]{%
    \errmessage{(Inkscape) Transparency is used (non-zero) for the text in Inkscape, but the package 'transparent.sty' is not loaded}%
    \renewcommand\transparent[1]{}%
  }%
  \providecommand\rotatebox[2]{#2}%
  \newcommand*\fsize{\dimexpr\f@size pt\relax}%
  \newcommand*\lineheight[1]{\fontsize{\fsize}{#1\fsize}\selectfont}%
  \ifx\svgwidth\undefined%
    \setlength{\unitlength}{228.79696583bp}%
    \ifx\svgscale\undefined%
      \relax%
    \else%
      \setlength{\unitlength}{\unitlength * \real{\svgscale}}%
    \fi%
  \else%
    \setlength{\unitlength}{\svgwidth}%
  \fi%
  \global\let\svgwidth\undefined%
  \global\let\svgscale\undefined%
  \makeatother%
  \begin{picture}(1,0.23123019)%
    \lineheight{1}%
    \setlength\tabcolsep{0pt}%
    \put(0,0){\includegraphics[width=\unitlength,page=1]{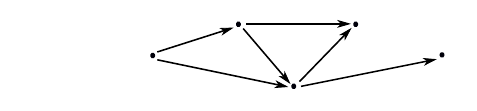}}%
    \put(0.25529844,0.0998988){\makebox(0,0)[lt]{\lineheight{1.25}\smash{\begin{tabular}[t]{l}$0$\end{tabular}}}}%
    \put(0.48328289,0.21315275){\makebox(0,0)[lt]{\lineheight{1.25}\smash{\begin{tabular}[t]{l}$1$\end{tabular}}}}%
    \put(0.60022828,0.00214886){\makebox(0,0)[lt]{\lineheight{1.25}\smash{\begin{tabular}[t]{l}$2$\end{tabular}}}}%
    \put(0.7397654,0.21462924){\makebox(0,0)[lt]{\lineheight{1.25}\smash{\begin{tabular}[t]{l}$3$\end{tabular}}}}%
    \put(0.96907889,0.1029996){\makebox(0,0)[lt]{\lineheight{1.25}\smash{\begin{tabular}[t]{l}$4$\end{tabular}}}}%
    \put(0,0){\includegraphics[width=\unitlength,page=2]{0234.pdf}}%
    \put(-0.00121548,0.10278668){\makebox(0,0)[lt]{\lineheight{1.25}\smash{\begin{tabular}[t]{l}$\Supp(P)=$\end{tabular}}}}%
  \end{picture}%
\endgroup%

	%\caption[]{$\Supp(e_{0123})$}
	%\label{Fig:Q1}
\end{figure}
The singular cubical map (non-unique) $f:I^{\boxdot3}\rightarrow\Supp(P)$ could be given by
\begin{align*}
f((0,0,0))=0,\quad f((1,0,0))=1,\quad f((0,1,0))=1,\quad f((0,0,1))=2,\\
f((1,1,0))=3,\quad f((1,0,1))=2,\quad f((0,1,1))=3,\quad f((1,1,1))=4.
\end{align*}
We also have the following deformation retractions:
\begin{enumerate}
  \item $r_1:\Supp(P)\rightarrow \im\left(f\big|_{I^{\boxdot 2}\boxdot\{0\}}\right)=\{0\rightarrow 1\rightarrow 3\}$,
  $$r_1(4)=r_1(3)=3,~r_1(2)=r_1(1)=1,~r_1(0)=0;$$
  \item $r_2:\Supp(P)\rightarrow \im\left(f\big|_{I\boxdot\{0\}\boxdot I}\right)=K_3(0,1,2)$,
  $$r_2(4)=r_2(3)=r_2(2)=2,~r_2(1)=1,~r_2(0)=0;$$
  %\item $r_3:\Supp(P)\rightarrow \im\left(f\big|_{\{0\}\boxdot I^{\boxdot 2}}\right)=S(0,1,2,3)$,
  %$$r_3(4)=r_3(3)=3,~r_3(2)=2,~r_3(1)=1,~r_3(0)=0;$$
  \item $r_4:\Supp(P)\rightarrow \im\left(f\big|_{I^{\boxdot 2}\boxdot\{1\}}\right)=K_3(2,3,4)$,
  $$r_4(4)=4,~r_4(3)=3,~r_4(2)=r_4(1)=r_4(0)=2;$$
  \item $r_5:\Supp(P)\rightarrow \im\left(f\big|_{I\boxdot\{1\}\boxdot I}\right)=\{1\rightarrow 3\rightarrow 4\}$,
  $$r_5(4)=4,~r_5(3)=r_5(2)=3,~r_5(1)=r_5(0)=1;$$
  %\item $r_6:\Supp(P)\rightarrow \im\left(f\big|_{\{1\}\boxdot I^{\boxdot 2}}\right)=S(1,2,3,4)$,
  %$$r_6(4)=4,~r_6(3)=3,~r_6(2)=2~,r_6(1)=r_6(0)=1.$$
\end{enumerate}
Note that $\im\left(f\big|_{\{0\}\boxdot I^{\boxdot 2}}\right)=S(0,1,2,3)$ and $\im\left(f\big|_{\{1\}\boxdot I^{\boxdot 2}}\right)=S(1,2,3,4)$, there are no digraph maps from $\Supp(P)$ onto these two faces, but there exist digraph maps onto their components, such as $r_1$ and $r_5$.
\end{example}

\begin{example}[$\NE(P)=3$, $\NF(P)=5$, i.e. Example \ref{NE32}] The admissible pair $(P,\Supp(P))$ with $\NE(P)=3$, $\NF(P)=5$ is given by
$$P=e_{0135}-e_{0235}+e_{0245},$$
\begin{figure}[H]
	\centering
	%% Creator: Inkscape 1.0.1 (3bc2e813f5, 2020-09-07), www.inkscape.org
%% Accompanies image file 'EX315old.pdf' (pdf, eps, ps)
\begingroup%
  \makeatletter%
  \providecommand\color[2][]{%
    \errmessage{(Inkscape) Color is used for the text in Inkscape, but the package 'color.sty' is not loaded}%
    \renewcommand\color[2][]{}%
  }%
  \providecommand\transparent[1]{%
    \errmessage{(Inkscape) Transparency is used (non-zero) for the text in Inkscape, but the package 'transparent.sty' is not loaded}%
    \renewcommand\transparent[1]{}%
  }%
  \providecommand\rotatebox[2]{#2}%
  \newcommand*\fsize{\dimexpr\f@size pt\relax}%
  \newcommand*\lineheight[1]{\fontsize{\fsize}{#1\fsize}\selectfont}%
  \ifx\svgwidth\undefined%
    \setlength{\unitlength}{234.09038634bp}%
    \ifx\svgscale\undefined%
      \relax%
    \else%
      \setlength{\unitlength}{\unitlength * \real{\svgscale}}%
    \fi%
  \else%
    \setlength{\unitlength}{\svgwidth}%
  \fi%
  \global\let\svgwidth\undefined%
  \global\let\svgscale\undefined%
  \makeatother%
  \begin{picture}(1,0.25197849)%
    \lineheight{1}%
    \setlength\tabcolsep{0pt}%
    \put(0,0){\includegraphics[width=\unitlength,page=1]{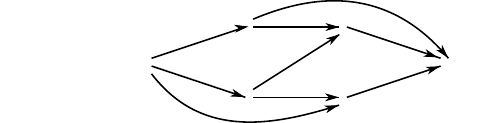}}%
    \put(0.25068877,0.10523922){\makebox(0,0)[lt]{\lineheight{1.25}\smash{\begin{tabular}[t]{l}$0$\end{tabular}}}}%
    \put(0.46185204,0.21954819){\makebox(0,0)[lt]{\lineheight{1.25}\smash{\begin{tabular}[t]{l}$1$\end{tabular}}}}%
    \put(0.47987608,0.08080941){\makebox(0,0)[lt]{\lineheight{1.25}\smash{\begin{tabular}[t]{l}$2$\end{tabular}}}}%
    \put(0.7017894,0.14904051){\makebox(0,0)[lt]{\lineheight{1.25}\smash{\begin{tabular}[t]{l}$3$\end{tabular}}}}%
    \put(0.70061249,0.00487309){\makebox(0,0)[lt]{\lineheight{1.25}\smash{\begin{tabular}[t]{l}$4$\end{tabular}}}}%
    \put(0.94352494,0.10575315){\makebox(0,0)[lt]{\lineheight{1.25}\smash{\begin{tabular}[t]{l}$5$\end{tabular}}}}%
    \put(0,0){\includegraphics[width=\unitlength,page=2]{EX315old.pdf}}%
    \put(-0.00297418,0.10774043){\makebox(0,0)[lt]{\lineheight{1.25}\smash{\begin{tabular}[t]{l}$\Supp(P)=$\end{tabular}}}}%
  \end{picture}%
\endgroup%

	%\caption[]{$\Supp(e_{0123})$}
	%\label{Fig:Q1}
\end{figure}
The singular cubical map (non-unique) $f:I^{\boxdot3}\rightarrow\Supp(P)$ could be given by
\begin{align*}
f((0,0,0))=0,\quad f((1,0,0))=1,\quad f((0,1,0))=2,\quad f((0,0,1))=4,\\
f((1,1,0))=3,\quad f((1,0,1))=5,\quad f((0,1,1))=4,\quad f((1,1,1))=5.
\end{align*}
We also have the following deformation retractions:
\begin{enumerate}
  \item $r_1:\Supp(P)\rightarrow \im\left(f\big|_{I^{\boxdot 2}\boxdot\{0\}}\right)=S(0,1,2,3)$,
  $$r_1(5)=r_1(3)=3,~r_1(4)=r_1(2)=2,~r_1(1)=1,~r_1(0)=0;$$
  \item $r_2:\Supp(P)\rightarrow \im\left(f\big|_{I\boxdot\{0\}\boxdot I}\right)=S(0,1,4,5)$,
  $$r_2(5)=5,~r_2(4)=4,~r_2(3)=r_2(1)=1,~r_2(2)=r_2(0)=0;$$
  \item $r_3:\Supp(P)\rightarrow \im\left(f\big|_{\{0\}\boxdot I^{\boxdot 2}}\right)=K_3(0,2,4)$,
  $$r_3(5)=r_3(4)=4,~r_3(3)=r_3(2)=2,~r_3(1)=r_3(0)=0;$$
  \item $r_4:\Supp(P)\rightarrow \im\left(f\big|_{I^{\boxdot 2}\boxdot\{1\}}\right)=\{4\rightarrow 5\}$,
  $$r_4(5)=r_4(3)=r_4(1)=5,~r_4(4)=r_4(2)=r_4(0)=4;$$
  \item $r_5:\Supp(P)\rightarrow \im\left(f\big|_{I\boxdot\{1\}\boxdot I}\right)=S(2,3,4,5)$,
  $$r_5(5)=5,~r_5(4)=4,~r_5(3)=r_5(1)=3,~r_5(2)=r_5(0)=2;$$
  \item $r_6:\Supp(P)\rightarrow \im\left(f\big|_{\{1\}\boxdot I^{\boxdot 2}}\right)=K_3(1,3,5)$,
  $$r_6(5)=r_6(4)=5,~r_6(3)=r_6(2)=3,~r_6(1)=r_6(0)=1.$$
\end{enumerate}
\end{example}

\begin{example}[$\NE(P)=4$, $\NF(P)=4$] The admissible pair $(P,\Supp(P))$ with $\NE(P)=4$, $\NF(P)=4$ is given by
$$P=e_{0135}-e_{0235}+e_{0245}-e_{0145},$$
\begin{figure}[H]
	\centering
	%% Creator: Inkscape 1.0.1 (3bc2e813f5, 2020-09-07), www.inkscape.org
%% Accompanies image file 'EX012345.pdf' (pdf, eps, ps)
\begingroup%
  \makeatletter%
  \providecommand\color[2][]{%
    \errmessage{(Inkscape) Color is used for the text in Inkscape, but the package 'color.sty' is not loaded}%
    \renewcommand\color[2][]{}%
  }%
  \providecommand\transparent[1]{%
    \errmessage{(Inkscape) Transparency is used (non-zero) for the text in Inkscape, but the package 'transparent.sty' is not loaded}%
    \renewcommand\transparent[1]{}%
  }%
  \providecommand\rotatebox[2]{#2}%
  \newcommand*\fsize{\dimexpr\f@size pt\relax}%
  \newcommand*\lineheight[1]{\fontsize{\fsize}{#1\fsize}\selectfont}%
  \ifx\svgwidth\undefined%
    \setlength{\unitlength}{223.10847005bp}%
    \ifx\svgscale\undefined%
      \relax%
    \else%
      \setlength{\unitlength}{\unitlength * \real{\svgscale}}%
    \fi%
  \else%
    \setlength{\unitlength}{\svgwidth}%
  \fi%
  \global\let\svgwidth\undefined%
  \global\let\svgscale\undefined%
  \makeatother%
  \begin{picture}(1,0.31815132)%
    \lineheight{1}%
    \setlength\tabcolsep{0pt}%
    \put(0,0){\includegraphics[width=\unitlength,page=1]{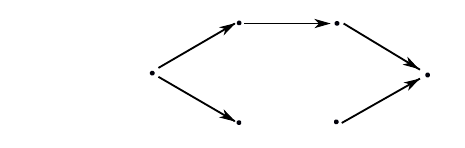}}%
    \put(-0.0021868,0.14789372){\makebox(0,0)[lt]{\lineheight{1.25}\smash{\begin{tabular}[t]{l}$\Supp(P)=$\end{tabular}}}}%
    \put(0.26636067,0.1485136){\makebox(0,0)[lt]{\lineheight{1.25}\smash{\begin{tabular}[t]{l}$0$\end{tabular}}}}%
    \put(0.48668306,0.2906922){\makebox(0,0)[lt]{\lineheight{1.25}\smash{\begin{tabular}[t]{l}$1$\end{tabular}}}}%
    \put(0.7011585,0.29138943){\makebox(0,0)[lt]{\lineheight{1.25}\smash{\begin{tabular}[t]{l}$3$\end{tabular}}}}%
    \put(0.94783225,0.14375306){\makebox(0,0)[lt]{\lineheight{1.25}\smash{\begin{tabular}[t]{l}$5$\end{tabular}}}}%
    \put(0,0){\includegraphics[width=\unitlength,page=2]{EX012345.pdf}}%
    \put(0.48417832,0.00346404){\makebox(0,0)[lt]{\lineheight{1.25}\smash{\begin{tabular}[t]{l}$2$\end{tabular}}}}%
    \put(0.69946788,0.00511006){\makebox(0,0)[lt]{\lineheight{1.25}\smash{\begin{tabular}[t]{l}$4$\end{tabular}}}}%
    \put(0,0){\includegraphics[width=\unitlength,page=3]{EX012345.pdf}}%
  \end{picture}%
\endgroup%

	%\caption[]{$\Supp(e_{0123})$}
	%\label{Fig:Q1}
\end{figure}
The singular cubical map (non-unique) $f:I^{\boxdot3}\rightarrow\Supp(P)$ could be given by
\begin{align*}
f((0,0,0))=0,\quad f((1,0,0))=1,\quad f((0,1,0))=2,\quad f((0,0,1))=2,\\
f((1,1,0))=3,\quad f((1,0,1))=4,\quad f((0,1,1))=4,\quad f((1,1,1))=5.
\end{align*}
We also have the following deformation retractions:
\begin{enumerate}
  \item $r_1:\Supp(P)\rightarrow \im\left(f\big|_{I^{\boxdot 2}\boxdot\{0\}}\right)=S(0,1,2,3)$,
  $$r_1(5)=3,~r_1(4)=2,~r_1(3)=3,~r_1(2)=2,~r_1(1)=1,~r_1(0)=0;$$
  \item $r_2:\Supp(P)\rightarrow \im\left(f\big|_{I\boxdot\{0\}\boxdot I}\right)=S(0,1,2,4)$,
  $$r_2(5)=4,~r_2(4)=4,~r_2(3)=1,~r_2(2)=2,~r_2(1)=1,~r_2(0)=0;$$
  \item $r_3:\Supp(P)\rightarrow \im\left(f\big|_{\{0\}\boxdot I^{\boxdot 2}}\right)=\{0\rightarrow 2\rightarrow 4\}$,
  $$r_3(5)=4,~r_3(4)=4,~r_3(3)=2,~r_3(2)=2,~r_3(1)=2,~r_3(0)=0;$$
  \item $r_4:\Supp(P)\rightarrow \im\left(f\big|_{I^{\boxdot 2}\boxdot\{1\}}\right)=\{2\rightarrow 4\rightarrow 5\}$,
  $$r_4(5)=5,~r_4(4)=4,~r_4(3)=4,~r_4(2)=2,~r_4(1)=2,~r_4(0)=2;$$
  \item $r_5:\Supp(P)\rightarrow \im\left(f\big|_{I\boxdot\{1\}\boxdot I}\right)=S(2,3,4,5)$,
  $$r_5(5)=5,~r_5(4)=4,~r_5(3)=3,~r_5(2)=2,~r_5(1)=2,~r_5(0)=2;$$
  \item $r_6:\Supp(P)\rightarrow \im\left(f\big|_{\{1\}\boxdot I^{\boxdot 2}}\right)=S(1,3,4,5)$,
  $$r_6(5)=5,~r_6(4)=4,~r_6(3)=3,~r_6(2)=1,~r_6(1)=1,~r_6(0)=1.$$
\end{enumerate}
The (3-step) homotopy $F$ between $i\circ r_1$ (similarly, $i\circ r_2$, $i\circ r_5$, $i\circ r_6$) and $\id_{\Supp(P)}$ can be found in Subsection 3.4 of \cite{TY}. Alternatively, if we relax the condition $r\big|_{H}=\id_H$ in the definition of retraction, we can also construct the following convenient digraph maps:
\begin{enumerate}
  \item $s_1:\Supp(P)\rightarrow \im\left(f\big|_{I^{\boxdot 2}\boxdot\{0\}}\right)=S(0,1,2,3)$,
  $$s_1(5)=3,~s_1(4)=2,~s_1(3)=1,~s_1(2)=0,~s_1(1)=0,~s_1(0)=0;$$
  \item $s_2:\Supp(P)\rightarrow \im\left(f\big|_{I\boxdot\{0\}\boxdot I}\right)=S(0,1,2,4)$,
  $$s_2(5)=4,~s_2(4)=2,~s_2(3)=1,~s_2(2)=0,~s_2(1)=0,~s_2(0)=0;$$
  \item $s_3:\Supp(P)\rightarrow \im\left(f\big|_{\{0\}\boxdot I^{\boxdot 2}}\right)=\{0\rightarrow 2\rightarrow 4\}$,
  $$s_3(5)=4,~s_3(4)=2,~s_3(3)=2,~s_3(2)=2,~s_3(1)=0,~s_3(0)=0;$$
  \item $s_4:\Supp(P)\rightarrow \im\left(f\big|_{I^{\boxdot 2}\boxdot\{1\}}\right)=\{2\rightarrow 4\rightarrow 5\}$,
  $$s_4(5)=5,~s_4(4)=4,~s_4(3)=5,~s_4(2)=4,~s_4(1)=4,~s_4(0)=2;$$
  \item $s_5:\Supp(P)\rightarrow \im\left(f\big|_{I\boxdot\{1\}\boxdot I}\right)=S(2,3,4,5)$,
  $$s_5(5)=5,~s_5(4)=5,~s_5(3)=5,~s_5(2)=3,~s_5(1)=4,~s_5(0)=2;$$
  \item $s_6:\Supp(P)\rightarrow \im\left(f\big|_{\{1\}\boxdot I^{\boxdot 2}}\right)=S(1,3,4,5)$,
  $$s_6(5)=5,~s_6(4)=5,~s_6(3)=5,~s_6(2)=4,~s_6(1)=3,~s_6(0)=1.$$
\end{enumerate}
These maps could play the key step in the homotopy to prove that $r_S:\Supp(P)\rightarrow\{0\}$ and $r_E:\Supp(P)\rightarrow\{5\}$ are deformation retractions.
\end{example}

\begin{example}[$\NE(P)=4$, $\NF(P)=6$, i.e. Example \ref{NE46}] The admissible pair $(P,\Supp(P))$ with $\NE(P)=4$, $\NF(P)=6$ is given by
$$P=e_{0136}-e_{0236}+e_{0256}-e_{0146},$$
\begin{figure}[H]
	\centering
	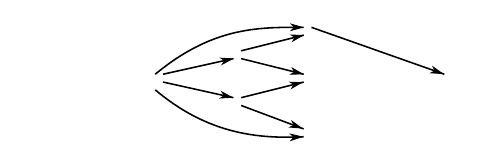
	%\caption[]{$\Supp(e_{0123})$}
	%\label{Fig:Q1}
\end{figure}
The singular cubical map (non-unique) $f:I^{\boxdot3}\rightarrow\Supp(P)$ could be given by
\begin{align*}
f((0,0,0))=0,\quad f((1,0,0))=1,\quad f((0,1,0))=0,\quad f((0,0,1))=2,\\
f((1,1,0))=4,\quad f((1,0,1))=3,\quad f((0,1,1))=5,\quad f((1,1,1))=6.
\end{align*}
We also have the following deformation retractions:
\begin{enumerate}
  \item $r_1:\Supp(P)\rightarrow \im\left(f\big|_{I^{\boxdot 2}\boxdot\{0\}}\right)=K_3(0,1,4)$,
  $$r_1(6)=4,~r_1(5)=0,~r_1(4)=4,~r_1(3)=1,~r_1(2)=0,~r_1(1)=1,~r_1(0)=0;$$
  \item $r_2:\Supp(P)\rightarrow \im\left(f\big|_{I\boxdot\{0\}\boxdot I}\right)=S(0,1,2,3)$,
  $$r_2(6)=3,~r_2(5)=2,~r_2(4)=1,~r_2(3)=3,~r_2(2)=2,~r_2(1)=1,~r_2(0)=0;$$
  \item $r_3:\Supp(P)\rightarrow \im\left(f\big|_{\{0\}\boxdot I^{\boxdot 2}}\right)=K_3(0,2,5)$,
  $$r_3(6)=5,~r_3(5)=5,~r_3(4)=0,~r_3(3)=2,~r_3(2)=2,~r_3(1)=0,~r_3(0)=0;$$
  \item $r_4:\Supp(P)\rightarrow \im\left(f\big|_{I^{\boxdot 2}\boxdot\{1\}}\right)=S(2,3,5,6)$,
  $$r_4(6)=6,~r_4(5)=5,~r_4(4)=3,~r_4(3)=3,~r_4(2)=2,~r_4(1)=3,~r_4(0)=2;$$
  \item $r_5:\Supp(P)\rightarrow \im\left(f\big|_{I\boxdot\{1\}\boxdot I}\right)=S(0,4,5,6)$,
  $$r_5(6)=6,~r_5(5)=5,~r_5(4)=4,~r_5(3)=6,~r_5(2)=5,~r_5(1)=4,~r_5(0)=0;$$
  \item $r_6:\Supp(P)\rightarrow \im\left(f\big|_{\{1\}\boxdot I^{\boxdot 2}}\right)=S(1,3,4,6)$,
  $$r_6(6)=6,~r_6(5)=3,~r_6(4)=4,~r_6(3)=3,~r_6(2)=3,~r_6(1)=1,~r_6(0)=1.$$
\end{enumerate}
The (3-step) homotopy $F$ between $i\circ r_4$ (similarly, $i\circ r_6$) and $\id_{\Supp(P)}$ can be found in Subsection 3.4 of \cite{TY}.
\end{example}

\begin{example}[$\NE(P)=5$, $\NF(P)=6$, i.e. Example \ref{NE5}] The admissible pair $(P,\Supp(P))$ with $\NE(P)=5$, $\NF(P)=6$ is given by
$$P=e_{0136}-e_{0156}+e_{0456}+e_{0246}-e_{0236},$$
\begin{figure}[H]
	\centering
	
	%\caption[]{$\Supp(e_{0123})$}
	%\label{Fig:Q1}
\end{figure}
The singular cubical map (non-unique) $f:I^{\boxdot3}\rightarrow\Supp(P)$ could be given by
\begin{align*}
f((0,0,0))=0,\quad f((1,0,0))=1,\quad f((0,1,0))=2,\quad f((0,0,1))=4,\\
f((1,1,0))=3,\quad f((1,0,1))=5,\quad f((0,1,1))=4,\quad f((1,1,1))=6.
\end{align*}
We also have the following obvious deformation retractions:
\begin{enumerate}
  %\item $r_1:\Supp(P)\rightarrow \im\left(f\big|_{I^{\boxdot 2}\boxdot\{0\}}\right)=S(0,1,2,3)$,
  %$$r_1(6)=3,~r_1(5)=r_1(4)=2,~r_1(3)=1,~r_1(2)=0,~r_1(1)=0,~r_1(0)=0;$$
  %$$r_1(6)=3,~r_1(5)=3,~r_1(4)=2,~r_1(3)=3,~r_1(2)=2,~r_1(1)=1,~r_1(0)=0;$$
  \item $r_2:\Supp(P)\rightarrow \im\left(f\big|_{I\boxdot\{0\}\boxdot I}\right)=S(0,1,4,5)$,
  $$r_2(6)=5,~r_2(5)=5,~r_2(4)=4,~r_2(3)=1,~r_2(2)=0,~r_2(1)=1,~r_2(0)=0;$$
  \item $r_3:\Supp(P)\rightarrow \im\left(f\big|_{\{0\}\boxdot I^{\boxdot 2}}\right)=K_3(0,2,4)$,
  $$r_3(6)=4,~r_3(5)=4,~r_3(4)=4,~r_3(3)=2,~r_3(2)=2,~r_3(1)=0,~r_3(0)=0;$$
  \item $r_4:\Supp(P)\rightarrow \im\left(f\big|_{I^{\boxdot 2}\boxdot\{1\}}\right)=K_3(4,5,6)$,
  $$r_4(6)=6,~r_4(5)=5,~r_4(4)=4,~r_4(3)=6,~r_4(2)=4,~r_4(1)=5,~r_4(0)=4;$$
  \item $r_5:\Supp(P)\rightarrow \im\left(f\big|_{I\boxdot\{1\}\boxdot I}\right)=S(2,3,4,6)$,
  $$r_5(6)=6,~r_5(5)=6,~r_5(4)=4,~r_5(3)=3,~r_5(2)=2,~r_5(1)=3,~r_5(0)=2;$$
  %\item $r_6:\Supp(P)\rightarrow \im\left(f\big|_{\{1\}\boxdot I^{\boxdot 2}}\right)=S(1,3,5,6)$,
  %$$r_6(6)=6,~r_6(5)=5,~r_6(4)=3,~r_6(3)=3,~r_6(2)=3,~r_6(1)=1,~r_6(0)=1.$$
\end{enumerate}
%It is easy to see that $i\circ r_1\simeq \id_{\Supp(P)}$ can not be realized by a one-step homotopy. Alternatively, one can also consider the digraph map $s_1:\Supp(P)\rightarrow \im\left(f\big|_{I^{\boxdot 2}\boxdot\{0\}}\right)=S(0,1,2,3)$:
%$$s_1(6)=3,~s_1(5)=s_1(4)=2,~s_1(3)=1,~s_1(2)=0,~s_1(1)=0,~s_1(0)=0.$$
%It is easy to check that $i\circ s_1\simeq \id_{\Supp(P)}$, $s_1\circ i\simeq\id_{S(0,1,2,3)}$. One can apply a similar argument to the homopopy equivalent between $S(1,3,4,6)$ and $\Supp(P)$.
\end{example}

\begin{remark} There is another minimal $3$-paths with $5$ elementary path components up to the global orientation of edges and isomorphism of digraphs. That is,
$$P=e_{0147}-e_{0247}+e_{0257}-e_{0357}+e_{0367}.$$
%\begin{align*}
%&P_1=e_{0147}-e_{0247}+e_{0257}-e_{0357}+e_{0367},\\
%&P_2=e_{0147}-e_{0247}+e_{0267}-e_{0157}+e_{0357}.
%\end{align*}
Its supporting digraph is given by the following one\footnote{This digraph is isomorphic to $\Supp(P)$ in Subsection 3.4 of \cite{TY}, where
$$P=e_{0137}-e_{0237}+e_{0267}-e_{0157}+e_{0457}.$$}.
\begin{figure}[H]
	\centering
	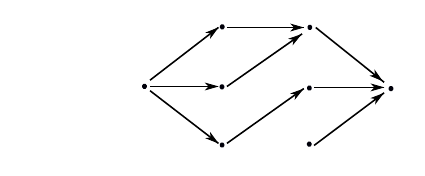
	%\caption[]{$\Supp(e_{0123})$}
	%\label{Fig:Q1}
\end{figure}

One can see that $\Supp(P)$ has $7$ face components, thus $P\notin\mathcal{P}_{\adm,3}$.
\end{remark}

\begin{example}[$\NE(P)=6$, $\NF(P)=6$] The admissible pair $(P,\Supp(P))$ with $\NE(P)=6$, $\NF(P)=6$ is given by
$$P=e_{0147}-e_{0157}+e_{0357}-e_{0367}+e_{0267}-e_{0167},$$
and the supporting digraph $\Supp(P)$ is nothing but the 3-cube $I^{\boxdot3}$. The corresponding six retraction digraph maps are the standard projections.
\end{example}

\end{appendices}

\end{document}